\documentclass{article}
\usepackage[pagebackref, colorlinks=true,linkcolor=red,citecolor=blue]{hyperref}
\usepackage{amsfonts,amssymb,amsmath} 

\usepackage{color}

\usepackage{latexsym}

\usepackage{amsthm}

\font\tengoth=eufm10 at 10pt
\font\sevengoth=eufm7 at 6pt
\newfam\gothfam
\textfont\gothfam=\tengoth
\scriptfont\gothfam=\sevengoth

\newcommand{\mlabel}[1]{\marginpar{#1}\label{#1}}

\newcommand{\fS}{{\mathfrak S}}

\newcommand{\g}{{\mathfrak g}}

\newcommand{\fb}{{\mathfrak b}}

\newcommand{\fe}{{\mathfrak e}}

\newcommand{\fg}{{\mathfrak g}}
\newcommand{\fh}{{\mathfrak h}}

\newcommand{\fk}{{\mathfrak k}}
\newcommand{\fl}{{\mathfrak l}}
\newcommand{\fm}{{\mathfrak m}}
\newcommand{\fn}{{\mathfrak n}}

\newcommand{\fp}{{\mathfrak p}}
\newcommand{\fr}{{\mathfrak r}}
\newcommand{\fs}{{\mathfrak s}}
\newcommand{\ft}{{\mathfrak t}}
\newcommand{\fu}{{\mathfrak u}}

\newcommand{\fz}{{\mathfrak z}}

\renewcommand\sp{\mathfrak {sp}} 
\newcommand\hsp{\mathfrak {hsp}} 
 
\newcommand\heis{\mathfrak {heis}}

\renewcommand{\:}{\colon}
\newcommand{\1}{\mathbf{1}}

\newcommand{\cA}{\mathcal{A}}

\newcommand{\cC}{\mathcal{C}}

\newcommand{\cH}{\mathcal{H}}

\newcommand{\cL}{\mathcal{L}}

\newcommand{\cO}{\mathcal{O}}

\newcommand{\cV}{\mathcal{V}}
\newcommand{\cW}{\mathcal{W}}

\newcommand{\eset}{\emptyset}

\newcommand{\dd}{{\tt d}}

\newcommand{\trile}{\trianglelefteq}
\newcommand{\subeq}{\subseteq}
\newcommand{\supeq}{\supseteq}

\newcommand{\into}{\hookrightarrow}

\newcommand{\shalf}{{\textstyle{\frac{1}{2}}}}

\def\onto{\to\mskip-14mu\to}

\newcommand{\N}{{\mathbb N}}
\newcommand{\Z}{{\mathbb Z}}
\newcommand{\R}{{\mathbb R}}
\newcommand{\C}{{\mathbb C}}

\newcommand{\bP}{{\mathbb P}}
\newcommand{\Q}{{\mathbb Q}}

\newcommand{\bE}{{\mathbb E}}

\newcommand{\bS}{{\mathbb S}}

\renewcommand{\hat}{\widehat}

\renewcommand{\tilde}{\widetilde}

\renewcommand{\L}{\mathop{\bf L{}}\nolimits}

\newcommand{\Aff}{\mathop{{\rm Aff}}\nolimits}

\newcommand{\GL}{\mathop{{\rm GL}}\nolimits}
\newcommand{\SL}{\mathop{{\rm SL}}\nolimits}

\newcommand{\PGL}{\mathop{{\rm PGL}}\nolimits}

\newcommand{\PSL}{\mathop{{\rm PSL}}\nolimits}
\newcommand{\SO}{\mathop{{\rm SO}}\nolimits}
\newcommand{\SU}{\mathop{{\rm SU}}\nolimits}

\newcommand{\U}{\mathop{\rm U{}}\nolimits}

\newcommand{\Sp}{\mathop{{\rm Sp}}\nolimits}
\newcommand{\Sym}{\mathop{{\rm Sym}}\nolimits}

\newcommand{\Ham}{\mathop{{\rm Ham}}\nolimits}

\newcommand{\aff}{\mathop{{\mathfrak{aff}}}\nolimits}

\newcommand{\mot} {\mathop{{\mathfrak{mot}}}\nolimits}

\newcommand{\fsl} {\mathop{{\mathfrak{sl} }}\nolimits}

\newcommand{\su}  {\mathop{{\mathfrak{su} }}\nolimits}
\newcommand{\so}  {\mathop{{\mathfrak{so} }}\nolimits}

\newcommand{\Fix}{\mathop{{\rm Fix}}\nolimits}

\newcommand{\ad}{\mathop{{\rm ad}}\nolimits}
\newcommand{\Ad}{\mathop{{\rm Ad}}\nolimits}

\newcommand{\epi}{\mathop{{\rm epi}}\nolimits}
\newcommand{\tr}{\mathop{{\rm tr}}\nolimits}

\newcommand{\Heis}{\mathop{{\rm Heis}}\nolimits}

\newcommand{\Aut}{\mathop{{\rm Aut}}\nolimits}

\newcommand{\Diff}{\mathop{{\rm Diff}}\nolimits}

\newcommand{\id}{\mathop{{\rm id}}\nolimits}

\newcommand{\rad}{\mathop{{\rm rad}}\nolimits}
\renewcommand{\dim}{\mathop{{\rm dim}}\nolimits}

\newcommand{\supp}{\mathop{{\rm supp}}\nolimits}

\newcommand{\Inn}{\mathop{{\rm Inn}}\nolimits}

\newcommand{\cone}{\mathop{{\rm cone}}\nolimits}
\newcommand{\conv}{\mathop{{\rm conv}}\nolimits}

\newcommand{\vol}{\mathop{{\rm vol}}\nolimits}

\newcommand{\Spann}{\mathop{{\rm span}}\nolimits}
\newcommand{\ev}{\mathop{{\rm ev}}\nolimits}

\newcommand{\PSO}{\mathop{{\rm PSO}}\nolimits}

\newcommand{\Rarrow}{\Rightarrow}
\newcommand{\nin}{\noindent} 
\newcommand{\oline}{\overline}

\newcommand{\la}{\langle}
\newcommand{\ra}{\rangle}

\newcommand{\Mot}{{\rm Mot}}

\newcommand{\res}{\vert}

\newcommand{\ssssarr}{\hbox to 15pt{\rightarrowfill}}
\newcommand{\sssarr}{\hbox to 20pt{\rightarrowfill}}
\newcommand{\ssarr}{\hbox to 30pt{\rightarrowfill}}
\newcommand{\sarr}{\hbox to 40pt{\rightarrowfill}}
\newcommand{\arr}{\hbox to 60pt{\rightarrowfill}}
\newcommand{\sssslarr}{\hbox to 15pt{\leftarrowfill}}
\newcommand{\ssslarr}{\hbox to 20pt{\leftarrowfill}}
\newcommand{\sslarr}{\hbox to 30pt{\leftarrowfill}}
\newcommand{\slarr}{\hbox to 40pt{\leftarrowfill}}
\newcommand{\larr}{\hbox to 60pt{\leftarrowfill}}

\newcommand{\Arr}{\hbox to 80pt{\rightarrowfill}}

\def\theoremname{Theorem}
\def\propositionname{Proposition}
\def\corollaryname{Corollary}
\def\lemmaname{Lemma}
\def\remarkname{Remark}
\def\conjecturename{Conjecture} 

\def\definitionname{Definition}
\def\exercisename{Exercise}
\def\examplename{Example}
\def\examplesname{Examples}
\def\problemname{Problem}
\def\problemsname{Problems}

\def\satzname{Satz} 
\def\koroname{Korollar}
\def\folgname{Folgerung}
\def\bemerkname{Bemerkung}
\def\aufgname{Aufgabe}

\def\beisname{Beispiel}
\def\beissname{Beispiele}
\def\bewname{Beweis}

\def\@thmcounter#1{\noexpand\arabic{#1}}
\def\@thmcountersep{}
\def\@begintheorem#1#2{\it \trivlist \item[\hskip 
\labelsep{\bf #1\ #2.\quad}]}
\def\@opargbegintheorem#1#2#3{\it \trivlist
      \item[\hskip \labelsep{\bf #1\ #2.\quad{\rm #3}}]}
\makeatother
\newtheorem{theor}{\theoremname}[section]
\newtheorem{propo}[theor]{\propositionname}
\newtheorem{coro}[theor]{\corollaryname}
\newtheorem{lemm}[theor]{\lemmaname}

\newenvironment{thm}{\begin{theor}\it}{\end{theor}}

\newenvironment{prop}{\begin{propo}\it}{\end{propo}}

\newenvironment{cor}{\begin{coro}\it}{\end{coro}}

\newenvironment{lem}{\begin{lemm}\it}{\end{lemm}}

\newenvironment{Lemma}{\begin{lemm}\it}{\end{lemm}}

\newtheorem{rema}[theor]{\remarkname}

\newenvironment{rem}{\begin{rema}\rm}{\end{rema}}

\newtheorem{stepnow}[theor]{}

\newtheorem{defin}[theor]{\definitionname} 

\newenvironment{defn}{\begin{defin}\rm}{\end{defin}}

\newtheorem{exerc}{\exercisename}[section]

\newtheorem{exa}[theor]{\examplename}

\newenvironment{ex}{\begin{exa}\rm}{\end{exa}}

\newtheorem{exas}[theor]{\examplesname}

\newtheorem{conj}[theor]{\conjecturename}

\newtheorem{pro}[theor]{\problemname}

\newtheorem{prs}[theor]{\problemsname}

\newtheorem{aufg}{\aufgname}[section]

\newenvironment{prf}{\begin{proof}}{\end{proof}}
 
\newcommand{\pmat}[1]{\begin{pmatrix} #1 \end{pmatrix}}

{\hfill\qed\end{trivlist}}

\newenvironment{beweis*}{\begin{trivlist}\item[\hskip%
\labelsep{\bf\bewname.\quad}]}%
{\end{trivlist}}

\newtheorem{satzn}[theor]{\satzname}

\newtheorem{koro}[theor]{\koroname}

\newtheorem{folg}[theor]{\folgname}

\newtheorem{bem}[theor]{\bemerkname}

\newtheorem{aufgn}[theor]{\aufgname}

\newtheorem{beis}[theor]{\beisname}

\newtheorem{beiss}[theor]{\beissname}

\renewcommand{\div}{\mathop{{\rm div}}\nolimits}
\newcommand{\comp}{\mathop{{\rm comp}}\nolimits}
\renewcommand{\phi}{\varphi} 
\renewcommand{\aff}{\mathop{{\rm aff}}\nolimits}
\newcommand{\co}{\mathop{{\rm co}}\nolimits}
\newcommand{\add}{\mathop{{\rm add}}\nolimits}

\newcommand{\be}{{\bf{e}}}

\newcommand{\bc}{{\bf{c}}}

\newtheorem{theora}{\theoremname}
\newenvironment{thma}{\begin{theora}\it}{\end{theora}}

\renewcommand{\mlabel}{\label} 

\title{A classification of coadjoint orbits \\
carrying Gibbs ensembles} 
\author{Karl-Hermann Neeb} 

\begin{document} 

\maketitle

\abstract{A coadjoint orbit $\cO_\lambda \subeq \g^*$ of a Lie group $G$
  is said to carry a Gibbs ensemble if the set of all $x \in \g$, for which the function  
  $\alpha \mapsto e^{-\alpha(x)}$ on the orbit is integrable with respect to the
  Liouville measure, has non-empty interior $\Omega_\lambda$.
  We describe a classification of all coadjoint orbits of finite-dimensional
  Lie algebras with this property. In the context of Souriau's
  Lie group thermodynamics, the subset $\Omega_\lambda$
  is the geometric temperature, a parameter space for a family
  of Gibbs measures on the coadjoint orbit. The
  corresponding Fenchel--Legendre transform maps
  $\Omega_\lambda/\fz(\g)$ diffeomorphically 
  onto the interior of the convex hull of the coadjoint orbit
  $\cO_\lambda$.   This provides an interesting perspective on the
  underlying information geometry.

We also show that already
  the integrability of $e^{-\alpha(x)}$ for one $x \in \g$ implies
  that $\Omega_\lambda \not=\eset$ and that, for general Hamiltonian
  actions, the existence of Gibbs measures implies that the range
  of the momentum maps consists of coadjoint orbits $\cO_\lambda$ as above.

  \nin  {\bf Keywords:} Hamiltonian action,
  Lie group thermodynamics, Gibbs measure, admissible Lie algebra,
  tempered coadjoint orbit, Liouville measure,\\  
  MSC: Primary 37J37; Secondary 22F30, 53D20, 58F05, 70H33,   82B05}

\tableofcontents

\section{Introduction} 
\mlabel{sec:1}

Let $G$ be a connected (finite-dimensional) Lie group 
with Lie algebra $\g$. 
The conjugation action of $G$ on itself induces on $\g$ the 
{\it adjoint action} 
$\Ad \: G \to \Aut(\g)$ and by dualization we obtain on
the dual space $\g^*$ the {\it coadjoint action} 
\[ \Ad^* \: G\to \GL(\g^*)\quad \mbox{  with  } \quad 
\Ad^*(g) \lambda := \lambda \circ \Ad(g)^{-1}.\] 
We call a subset of $\g$, resp., $\g^*$ {\it invariant} if it is invariant 
under $\Ad(G)$, resp., $\Ad^*(G)$. 
Let \break $\sigma \: G \times M \to M$ be a (strongly) Hamiltonian action
of the Lie group $G$ on the symplectic manifold 
$(M,\omega)$ and 
\[ \Psi \: M \to \g^* \]
the corresponding
equivariant momentum map.\begin{footnote}{One
    also  studies Hamiltonian actions for which all vector fields
    $\dot\sigma(x)$ come from Hamiltonian functions, but no equivariant
    momentum map $M \to \g^*$ exists. If $M$ is connected, this can be
    overcome by replacing the Lie algebra $\g$ by a suitable central 
    extension $\hat\g$. Taking this into account, it is no loss of generality
    to assume, as we do throughout,
    the existence of an equivariant momentum map, resp.,
    that the action is {\it strongly Hamiltonian}. We refer to
  Section~\ref{sec:non-strong} for a discussion of this issue in our context.}
\end{footnote} Then 
$H_x(m) := \Psi(m)(x)$ is the Hamiltonian function of $x \in \g$,
i.e., $\dd H_x = -i_{\dot\sigma(x)}\omega$ holds for the vector field
$\dot\sigma(x)$ of the derived action $\dot\sigma \:  \g\to \cV(M)$ 
(cf.\ \cite{GS84}). 
We write $\lambda_M$ for the Liouville measure on $M$,
specified by the volume form $\frac{\omega^n}{(2\pi)^n n!}$,
where $\dim M = 2n$. 
The open subset
\[ \Omega := \Big\{ x \in \g \:  \int_M e^{- H_x(m)}
  \, d\lambda_M(m) < \infty \Big\}^\circ  \]
is called the corresponding {\it geometric temperature}.
This is an open subset of $\g$ and the 
Laplace transform of the push-forward measure $\mu := \Psi_*\lambda_M$
on $\g^*$ defines on $\Omega$ an analytic convex function: 
\begin{equation}
  \label{eq:Z}
 Z(x) := \cL(\mu)(x) = \int_{\g^*} e^{-\alpha(x)} \, d\mu(\alpha)
 = \int_M e^{- H_x(m)} \, d\lambda_M(m).
\end{equation}
In Statistical Mechanics \eqref{eq:Z} corresponds to the {\it partition function},
hence the notation~$Z(x)$. 
The family of the probability measures 
$\lambda_x = \frac{e^{-H_x}}{Z(x)} \lambda_M$ is called 
the {\it Gibbs ensemble of the dynamical group $G$} acting on $M$.
The specific form of the density of Gibbs measures can be characterized
among all measures with smooth density and the same expectation value
in $\g^*$ by the maximality of their entropy 
(cf.~Remark~\ref{rem:entropy} and Theorem~\ref{thm:b.3}).
Therefore the Gibbs measures are natural models of equilibrium states
in thermodynamical systems. \\

Generalized temperatures
of a Hamiltonian action of a Lie group were introduced by 
J.-M.~Souriau in~\cite{So66, So75} and elaborated in
\cite[Ch.~IV]{So97}, as {\it Lie group   thermodynamics}. 
The idea was, that the momentum map $\Psi \: M \to \g^*$ 
of a Hamiltonian action generalizes the case where $\g$ is one-dimensional,
where $\Psi$ corresponds to the energy function of an isolated system. 
In Statistical Mechanics, the probability density of a state is given
in terms of the energy $E$ by the {\it Boltzmann distribution}
\[ P_\beta(E) = \frac{1}{Z(\beta)} e^{-\beta E},\]
where $\beta > 0$ corresponds to the inverse temperature, and
the {\it partition function} $Z(\beta)$ is a normalizing factor.
Souriau now replaces the ``inverse temperature'' $\beta = \frac{1}{kT}$
by a Lie algebra element~$x$, 
so that we obtain Gibbs measures $\lambda_x$ as above.

The building blocks for Hamiltonian actions are the transitive
ones (cf.\ Subsection~\ref{subsec:6.5}).
Then the momentum map $\Psi$ is a covering map from $M$ onto a coadjoint 
orbit $\cO_\lambda := \Ad^*(G)\lambda \subeq \g^*$.
One of our main results is a classification of those coadjoint orbits
for which the corresponding
geometric temperature $\Omega_\lambda$ is non-empty, i.e.,
for which the Laplace transform of the Liouville measure
$\mu_\lambda$ on $\cO_\lambda$ is finite on an open subset of $\g$.
\begin{footnote}{In many interesting situations $\cO_\lambda$ is simply connected
    and $\Psi$ is a diffeomorphism, but this is not always the case.
    The nilpotent coadjoint orbits in $\fsl_2(\R)^*$ are examples
    with $\pi_1(\cO_\lambda) \cong \Z$. Although $\Omega_\lambda \not=\eset$
    in this case, for the action of
    $\tilde\SL_2(\R)$ on its simply connected covering $\tilde\cO_\lambda$,
    all functions $e^{-H_x}$ have infinite integral.
}\end{footnote}

To this end, we may factorize the ideal $\cO_\lambda^\bot
= \{ x \in \g \: (\forall \alpha \in \cO_\lambda) \ \alpha(x) = 0\}\trile \g$
and thereafter assume that $\cO_\lambda$ spans $\g^*$. This
entails in particular that
$\dim \fz(\g) \leq 1$ because central elements define constant
Hamiltonian functions on $\cO_\lambda$. 
The first key observation is that, if $\cO_\lambda$ spans $\g^*$ and
\begin{equation}
  \label{eq:dlambda-intro}
  D_{\mu_\lambda} := \{ x \in \g \: \cL(\mu_\lambda)(x)< \infty\} \not=\eset,
\end{equation}
then the Lie algebra $\g$ is {\it admissible} 
(Theorem~\ref{thm:5.9}), i.e.,
contains a generating closed convex $\Ad(G)$-invariant 
subset not containing affine lines.

Admissible Lie algebras have a well-developed structure theory, exposed
in detail in the monograph \cite{Ne00}. Key facts are:
\begin{itemize}
\item A simple Lie algebra is admissible if and only if it is compact
  or hermitian, i.e., non-compact with non-trivial invariant convex cones
  (cf.\ \cite{Vi80}).  
\item Reductive Lie algebras are admissible if and only if their simple
  ideals are compact or hermitian.
\item For a symplectic vector space $(V,\Omega)$, the Jacobi--Lie algebra 
  $\hsp(V,\Omega) \cong \heis(V,\Omega) \rtimes \sp(V,\Omega)$
  of polynomials of degree $\leq 2$ on $V$, with respect to the Poisson
  bracket, is admissible (cf.\ \cite[App.~A.IV]{Ne00}).
\item Non-reductive admissible Lie algebras with at most one-dimensional
  center are semidirect sums
  $\g = \heis(V,\Omega)\rtimes_\sigma \fl$, where $\fl$ is reductive admissible
  with a homomorphism $\sigma \: \fl \to \sp(V,\Omega)$, satisfying
  certain positivity properties; see Subsection~\ref{subsec:1.4} for
  details. 
\end{itemize}

An important structural feature of admissible Lie algebras is
that they contain a compactly embedded Cartan subalgebra $\ft$
(cf.\ \cite{HH89}), 
\begin{footnote}{We call a subalgebra $\fb \subeq \g$ {\it compactly embedded}
    if the subgroup of $\Aut(\g)$ generated by $e^{\ad \fb}$ has compact closure.}
\end{footnote}
so that we obtain a
root decomposition
\[ \g_\C = \ft_\C \oplus \bigoplus_{\alpha \in \Delta} \g_\C^\alpha
  \quad\mbox{ with } \quad 
  \g_\C^\alpha = \{ z \in \g_\C \: (\forall x \in \ft)\
  [x,z] = \alpha(x)z\} 
  \quad \mbox{ and } \quad \Delta \subeq i \ft^*\]
(\cite[Thm.~VII.2.2]{Ne00}). 
In addition, there exists a unique maximal compactly embedded subalgebra
$\fk \subeq \g$, containing $\ft$
(\cite[Prop.~VII.2.5]{Ne00}). It specifies a subset
$\Delta_k := \{ \alpha \in \Delta\: \g_\C^\alpha \subeq \fk_\C\}$
of {\it compact roots}, and the corresponding
reflections generate a Weyl group~$\cW_{\fk}$, acting on~$\ft$ and~$\Delta$.
A positive system $\Delta^+ \subeq \Delta$ of roots is said to be
{\it adapted}, if the set $\Delta_p^+ := \Delta^+ \setminus \Delta_k$
of positive non-compact roots is invariant under the Weyl group~$\cW_{\fk}$
(\cite[Def.~VII.2.6, Prop.~VII.2.12]{Ne00}).
For $z = x + iy \in \g_\C$, we put
$z^* := -x + iy$ and 
associate to any such system two $\cW_{\fk}$-invariant convex cones in $\ft$:
\begin{equation}
  \label{eq:maxcon-into}
C_{\rm min} := \oline\cone(\{ i [x_\alpha, x_\alpha^*] \: 
x_\alpha \in \g_\C^\alpha, \alpha \in \Delta_p^+\})\subeq \ft,  
\end{equation} 
and 
\begin{equation}
  \label{eq:maxcone2-into}
  C_{\rm max} :=
  \{ x \in \ft \: (\forall \alpha \in \Delta_p^+) \ i\alpha(x) \geq 0\}
\end{equation}
(\cite[Def.~VII.3.6]{Ne00}). 
On the level of $\g$, they correspond to the cones
\[ W_{\rm max} := \{ y \in \g \:  p_\ft(\Ad(G)y) \subeq C_{\rm max} \}
\quad\mbox{ and }  \quad   W_{\rm min} := \{ y \in \g \:
  p_\ft(\Ad(G)y) \subeq C_{\rm min} \},\]
where $p_\ft \: \g \to \ft$ is the projection with kernel
$[\ft,\g]$ (\cite[Prop.~VIII.3.7]{Ne00}).\begin{footnote}
  {The terminology is motivated by the case of simple hermitian
    Lie algebras, where $W_{\rm min}$ is a minimal generating
    invariant cone and $W_{\rm max}$ is maximal.}\end{footnote}
If $C_{\rm min} \subeq C_{\rm max}$, then
$W_{\rm min} \subeq W_{\rm max}$ by definition
(cf.\ \cite[Thm.~VIII.3.8]{Ne00}). 
We are now ready to formulate our first main result.

\begin{thma} {\rm(Classification Theorem)} 
   Let $\cO_\lambda \subeq \g$ be a coadjoint orbit spanning $\g^*$.
   Then $\Omega_\lambda\not=\eset$  if and only if
     $\g$ is admissible   and there exists an adapted positive system
     $\Delta^+$ with $C_{\rm min}$ pointed and contained in $C_{\rm max}$
     such that $\lambda \in W_{\rm min}^\star := \{ \beta \in \g^* \:
     \beta(W_{\rm min}) \subeq [0,\infty)\}$.
\end{thma}

This result is contained in Theorem~\ref{conj:4.4}.  
Our strategy to obtain this classification is as follows:
If $D_{\mu_\lambda} \not=\eset$ (cf.\ \eqref{eq:dlambda-intro}), then quite general arguments
show that $\g$ is admissible and that
$\lambda \in W_{\rm min}^\star$ for an adapted  positive system $\Delta^+$
as above (Theorem~\ref{thm:5.9}).

The converse is harder. The main ingredients are: 
\begin{itemize}
\item The coadjoint orbit $\cO_\lambda$ of $\hsp(V,\Omega)$
  corresponding to the affine symplectic action on $(V,\Omega)$
  satisfies $\Omega_\lambda \not=\eset$. This can be seen by direct
  evaluation of  Gaussian integrals.
\item If $\lambda \in C_{\rm min}^\star$, then $\cO_\lambda$ is a so-called
  {\it admissible orbit}, i.e., closed, and its convex hull contains
  no affine lines (\cite[Def.~VII.3.14]{Ne00}). For these orbits, there exist explicit formulas
  for the Laplace transform $\cL(\mu_\lambda)$, based on stationary
  phase methods (Duistermaat--Heckman formulas), that have been obtained
  in \cite{Ne96a}. They imply that $\Omega_\lambda \not=\eset$ in
  this case (Subsection~\ref{sec:5.3}).
\item If $\g$ is not reductive, then $\cO_\lambda$ is a symplectic product
  of an orbit corresponding to an affine action on a symplectic
  vector space and an orbit of a reductive Lie algebra.
  Since the affine case has been dealt with explicitly, this reduces 
  our problem to reductive, and hence to simple Lie algebras
  (Subsection~\ref{subsec:6.4}). 
\item If $\g$ is a compact simple Lie algebra, then $W_{\rm min} = \{0\}$ and
  all coadjoint orbits satisfy $\Omega_\lambda = \g$ because $\mu_\lambda$
  is a finite measure.   
\item The most difficult case are orbits of simple hermitian Lie algebras
  that are admissible. Then we have a Jordan decomposition
  $\lambda = \lambda_s + \lambda_n$ with $\lambda_s, \lambda_n
  \in W_{\rm min}^\star$,
  $\lambda_n$ nilpotent and $\lambda_s$ semisimple  (cf.\ \cite{NO22}).
  \begin{footnote}{Here we use that the Cartan--Killing form
      $\kappa$ on $\g$ induces a $G$-equivariant linear isomorphism
      $\g \to \g^*, x \mapsto \kappa(x,\cdot)$. Accordingly, we translate
      the Jordan decomposition from elements of $\g$ to elements of $\g^*$.    
}  \end{footnote}
Here $\cO_{\lambda_s}$ is admissible, a case we already dealt with,
  and the Liouville measure on the nilpotent orbit
  $\cO_{\lambda_n}$ can be treated with methods from \cite{Rao72},
  which imply that it is tempered. Since it is contained in a pointed
  cone, $\Omega_{\lambda_n} \not=\eset$ follows from Borcher's
  Theorem on tempered distributions
  (cf.\ Proposition~\ref{prop:dom-lapl-temp}).
  The Liouville measure $\mu_\lambda$ is a ``fibered product'' of
  $\mu_{\lambda_s}$ and a nilpotent Liouville measure of the centralizer
  $\fl$ of $\lambda_s$ (\cite{Rao72}).
  To deal with this situation, we prove a convexity
  theorem for the projection $p \:  \g \to \fl$ to show that
  $\Omega_\lambda\not=\eset$.
\end{itemize}

The strategy described above further shows that,
for $\lambda \in W_{\rm min}^\star$, the geometric temperature 
$\Omega_\lambda$ is the open convex cone~$W_{\rm  max}^\circ$.
This does not tell us anything about the finiteness of
$\cL(\mu_\lambda)$ in boundary points of this cone, but we also have:

\begin{thma} {\rm(Domain  Theorem)}  \mlabel{thm:dom}
  Suppose that $\cO_\lambda$ spans $\g^*$.
 If $\g$ is admissible with compactly embedded Cartan subalgebra $\ft$ 
 and $\Delta^+$ is adapted with $C_{\rm min}$ pointed and contained in
 $C_{\rm max}$, then $\lambda \in W_{\rm min}^\star$ implies that 
  $W_{\rm max}^\circ = D_{\mu_\lambda} = \Omega_\lambda.$
  In particular, the domain $D_{\mu_\lambda}$ of the Laplace transform
  $\cL(\mu_\lambda)$ is open. 
\end{thma}

The central argument for this theorem is the observation
that $\cL(\mu_\lambda)(x) < \infty$ leads to an invariant probability 
measure on the dual of the Lie subalgebra $\fz_\g(x) = \ker(\ad x)$ whose support is generating.
To show that this
can only happen for $x \in W_{\rm max}^\circ$, we use the following rather
general tool (Theorem~\ref{thm:measlemb}): 

\begin{thma} {\rm(Compactness Theorem)} \mlabel{thm:measlem}
Let $V$ be a finite dimensional real vector space. 
\begin{itemize}
\item[\rm(a)] If $\mu$ is a finite positive Borel measure on $V$ whose
  support spans $V$, then its stabilizer group
  $\GL(V)^\mu := \{g \in \GL(V) \: g_*\mu = \mu\}$ is closed and
  has the property that  all its elements are elliptic, i.e., generate relatively
  compact subgroups of $\GL(V)$.
\item[\rm(b)] If $H \subeq \GL(V)$ is a closed subgroup, such that
  all elements of $H$ are elliptic, then $H$ is compact. 
\end{itemize}
\end{thma}

For a coadjoint $\cO_\lambda$ with Liouville measure $\mu_\lambda$
and $Z_\lambda = \cL(\mu_\lambda)$, we have in the context of
Theorem~\ref{thm:dom} the analytic function 
\begin{equation}
  \label{eq:qdef-intro}
  Q \: \Omega_\lambda = W_{\rm max}^\circ  \to \g^*, \quad
  Q(x) := -\dd \log Z_\lambda(x)= \frac{1}{Z_\lambda(x)}
  \int_{\g^*} \alpha e^{-\alpha(x)}\, d\mu_\lambda(\alpha). 
\end{equation}
It associates to $x$ the expectation value of the probability measure
\[ d\lambda_x(\alpha)
  = \frac{e^{-\alpha(x)}}{Z_\lambda(x)}\, d\mu_\lambda(\alpha) \]  on
$\cO_\lambda$, hence $Q(x)$ is contained in its closed convex hull, but we actually
have much finer information.
The Domain Theorem implies that the smooth convex function $Z_\lambda$
on $\Omega_\lambda$ has a closed epigraph.
One can now derive from Fenchel's Convexity Theorem
(\cite[Thm.~V.3.31]{Ne00}, \cite[Thm.~1.16]{Ne19})
that $Q$ factors through a diffeomorphism
\[ \oline Q \: \Omega_\lambda/\fz(\g) =
  W_{\rm max}^\circ/\fz(\g) \to \conv(\cO_\lambda)^\circ \]
onto the relative interior of the convex hull of $\cO_\lambda$
(Theorem~\ref{conj:4.4}). Here the main point is the determination 
of the range of this map. That it is a diffeomorphism onto an open
subset follows from rather general facts on Laplace transforms. \\

The structure of this paper is as follows. 
In Section~\ref{app:a} we collect the relevant general material
on convex functions and Laplace transforms of measures.
In Section~\ref{app:coad-fin-vol}
we prove the Compactness Theorem. 
Section~\ref{sec:2} contains material on admissible Lie algebras,
supplemented by new results relating to invariant measures on $\g^*$ and their
Laplace transforms. For instance, Theorem~\ref{thm:5.9}
shows that, if $\mu$ is an invariant measure on $\g^*$,
whose support spans $\g^*$, and $D_\mu\not=\eset$, 
then $\g$ is admissible and $\supp(\mu) \subeq W_{\rm  min}^\star$
for an adated positive system.
In Section~\ref{sec:3} we briefly recall the concepts related to 
symplectic Gibbs ensembles.  
In Section~\ref{sec:5} we initialize the proof of the Classification Theorem with the
observation that $\g$ needs to be admissible and that
$\lambda \in W_{\rm  min}^\star$ is necessary for $\Omega_\lambda \not=\eset$
(a consequence of Theorem~\ref{thm:5.9} for Liouville measures of coadjoint orbits).
We then inspect the action on a symplectic vector space and
on admissible coadjoint orbits.
In Section~\ref{sec:6} we first treat nilpotent coadjoint
orbits in simple Lie algebras,
then mixed orbits, and finally split the problem
into the affine action of $\Heis(V,\Omega) \rtimes \Sp(V,\Omega)$
on $(V,\Omega)$ and the case of reductive Lie algebras.
In Section~\ref{sec:disint} we show that the measure 
$\mu$ always disintegrates into Liouville measures
on coadjoint orbits (Theorem~\ref{thm:disint}).
Finally, we discuss in Section~\ref{sec:non-strong} how to translate
our results to the context of non-strongly Hamiltonian actions,
where the momentum map is covariant with respect to a suitable
affine action of $G$ on $\g^*$. 

We conclude with a brief discussion of interesting perspectives
in Section~\ref{sec:7}. 
In particular, it would be interesting to develop a closer connection
between Gibbs ensembles on coadjoint orbits and
Gibbs states of the $C^*$-algebra $B(\cH)$, $\cH$ a complex Hilbert space.
They should be closely related to the KMS states
studied in \cite{Si23} for unitary highest weight representations
$(U,\cH)$. Then the operators $e^{-i \partial U(x)}$, $x \in W_{\rm max}^\circ$,
are trace class, so that $(U(\exp tx))_{t \in \R}$ is a unitary one-parameter
group with a unique Gibbs state for any inverse temperature $\beta > 0$.
On the ``classical side'', in $\g^*$, we find, by the Domain Theorem,
the same parameter space $W_{\rm max}^\circ$
for the Gibbs ensemble on~$\cO_\lambda$. This shows that,
for finite-dimensional Lie algebras, 
Gibbs ensembles on $\g^*$ and Gibbs states in unitary representations
share the same geometric environment. 

It is also interesting to connect all this with
information geometry. In this context, the key structure
is the {\it Fisher--Rao metric} on $\Omega$ (cf.\ \cite{Fr91}).
\begin{footnote}{It is called {\it geometric capacity} by Souriau and
{\it heat capacity} by Barbaresco.}  
\end{footnote} 
 It is given by the second
differential 
\begin{equation}
  \label{eq:d2Z}
  (\dd^2 \log Z)(x)(v,w)
  = \bE_{\lambda_x}[(H_v - \oline H_v)(H_w - \oline H_w)] \geq 0, 
  \quad \mbox{ where }  \quad
  \oline H_v := \bE_{\lambda_x}[H_v].
\end{equation}
This is positive definite if the convex hull of the support of
$\mu = \Phi_*\lambda_M$ has interior points, because then no non-zero
function $H_v$ is constant (cf.~Proposition~\ref{prop:I.9}(iii)).
This part of Souriau's work was taken up by
Barbaresco in \cite{Ba16}, who observed that the metric defined by
Souriau in \cite{So75}
coincides with the Fisher--Rao metric in the context
of statistical manifolds in 
information geometry (see also \cite[\S 4.3]{Neu22}
\cite{Ko61} and \cite{Sh07} for metrics defined by Hessians of convex functions
on domains in vector spaces). 
Souriau's concepts have
been translated to modern terminology and explored further
by Marle in \cite{Ma20a, Ma20b, Ma21}; see also the interesting
discussion in \cite[\S 5]{Bo19}. For the link with the
thermodynamics of continua, we refer to \cite{dS16}.

Souriau discusses in \cite{So97} the Galilei group
$\R^4 \rtimes \Mot_3(\R)$ 
and the Poincar\'e group $\R^{1,3} \rtimes \SO_{1,3}(\R)_e$.
In both cases (relativistic and non-relativistic),
he finds that no coadjoint orbit 
with non-trivial geometric temperature exists, so that it is
necessary to restrict to subgroups. We refer to Souriau's book
for an interesting discussion of the physical interpretations
of this fact, f.i., for the Galilei group, the non-existence of Gibbs states
is related to the universe being expanding and not stationary.
In both cases, the subgroup $\R^4  \times \SO_3(\R)$ has admissible
central extensions, to which our results apply.
In \cite[(17.136)]{So97}, the subgroup $H = \R \times \SO_3(\R)$
of the Poincar\'e group is discussed in connection with a relativistic
ideal gas.

In \cite{BDNP23} it was shown that, in a hermitian simple Lie algebra,
the minimal nilpotent orbit has non-empty geometric temperature,
and that, for the nilpotent orbit of
$\g = \fsl_2(\R)$, the 
Fisher--Rao metric turns the Gibbs cone $Q(\Omega_\lambda)
= \conv(\cO_\lambda)^\circ$ into a Riemannian symmetric space. \\

\nin {\bf Non-transitive actions:} In the present paper we determine
all coadjoint orbits for which the domain of the Laplace transform
of the Liouville measure is non-empty. In general, 
Souriau's Lie group thermodynamics leads to an $\Ad^*(G)$-invariant 
measure $\mu$ on $\g^*$ whose support spans $\g^*$ and
for which $\cL(\mu)$ is finite in some point of $\g$.
Then Theorem~\ref{thm:5.9} shows that
$\Psi(M) \subeq W_{\rm min}^\star$ for an adapted positive
system $\Delta^+$ with $C_{\rm min}$ pointed and contained in~$C_{\rm max}$.
In Theorem~\ref{thm:disint}, we show that
there exists a measurable subset $S \subeq \Psi(M)$ 
and a measure $\nu$ on $S$, for which
\begin{equation}
  \label{eq:mu-disint}
 \mu = \int_S \mu_\lambda\,  d\nu(\lambda),
  \quad \mbox{ and thus } \quad
  \cL(\mu)(x) = \int_S \cL(\mu_\lambda)(x)\, d\nu(\lambda).
\end{equation}
Since $\cL(\mu_\lambda)(x) < \infty$ for all $x \in W_{\rm max}^\circ$
by the Domain Theorem~\ref{thm:dom}, the finiteness properties of
$\cL(\mu)$ only depend on the measure $\nu$ on the cross section.
We show in Section~\ref{sec:disint}
that, if $\cC \subeq W_{\rm min}^\star$ is open and
contains no affine lines and $G$ contains a lattice $\Gamma$,
i.e., $\Gamma$ is discrete with $\vol(G/\Gamma) < \infty$,
then the restriction of Lebesgue measure $\lambda_{\g^*}$ to $\cC$
occurs as $\mu$ for $M \subeq T^*(\Gamma \backslash G)$,  and \eqref{eq:mu-disint} provides a
``Plancherel decomposition'' of $\lambda_{\g^*}\res_{\cC}$
into Liouville measures on coadjoint orbits.
If $\g$ is abelian, then all coadjoint orbits are trivial
and the Liouville measures $\mu_\lambda$ are point measures,
so that $\cL(\mu) = \cL(\nu)$. \\


\nin{\bf Acknowledgment:} {We are grateful to
  Tobias Simon for reading a first draft of this paper and
  for numerous comments that helped to improve the exposition.
  We also thank Toshiyuki Kobayashi for pointing out
  Rao's paper \cite{Rao72} that we used to deal with nilpotent orbits.
  We thank Pierre Bieliavsky for inspiring discussions on
  information geometry and for pointing out \cite{BDNP23}
  and \cite{Neu22}. We are also indepted to
  Nicolo Drago for discussions on geometric KMS states
  and Weinstein's paper \cite{We97}, which started the whole project.
  Further thanks go to Yves Cornulier for an inspiring  email
  exchange on linear torsion groups that led to a very effective
  proof of the Compactness Theorem~\ref{thm:measlem}.
  We thank F.~Barbaresco for supplying us generously with references
  on the connection between Souriau's Lie group thermodynamics, 
  information geometry, and Koszul's work on transitive affine actions
  on convex domains.   Finally we thank  Yoshiki Oshima for pointing out
  du Cloux's paper \cite{dCl91}
  as a means to show that Liouville measures of coadjoint orbits
  of reductive Lie algebras are tempered.
  In the end, we managed to bypass du Cloux's elaborate machinery 
  of Schwartz functions on semialgebraic varieties 
  by using a suitable Convexity Theorem~\ref{thm:conv-ell} for orbit
  projections. This was enough to obtain the desired
  finiteness of the Laplace transforms and even to show that
  $\mu_\lambda$ is tempered if $\Omega_\lambda\not=\eset$
  (Theorem~\ref{thm:temp}). 
}

\section{Convex sets and functions} 
\mlabel{app:a}

In this section we collect some some basic facts on 
convex sets, convex functions, and Laplace transforms
of positive measures.

Let $V$ be a finite-dimensional real vector space and $V^*$ be its dual space. 
We write $\la \alpha, v \ra = \alpha(v)$ for the natural pairing 
$V^* \times V \to \R$. For a subset $C \subeq V^*$, we consider the 
{\it dual cone} 
\begin{equation}
  \label{eq:dualcone}
 C^\star := \{ v \in V \: (\forall \alpha \in C)\ \alpha(v) \geq 0\} \quad 
 \mbox{ and also }  \quad  B(C) := \{ v \in V \: \inf\la C, v \ra > - \infty\}
\end{equation}
(cf.\ \cite[\S V.1]{Ne00}). 
Both are convex cones and $C^\star$ is closed.  
For a convex subset $C\subeq V$, we define 
its {\it recession cone} 
\begin{equation}
  \label{eq:limcone}
 \lim(C) := \{ x \in V \: C + x \subeq C \}  
\quad \mbox{ and } \quad H(C) := \lim(C) \cap -\lim(C)
= \{ x \in V \: C + x = C \}.
\end{equation}
Then $\lim(C)$ is a convex cone and $H(C)$ a linear subspace.
We write $C^\circ$ for the interior of $C$ in the affine subspace
$\aff(C)$ generated by $C$. Note that $C^\circ \not=\eset$ whenever $C \not=\eset$. 

\begin{Lemma} \label{lem:limcone} {\rm(\cite[Lemma~2.9]{Ne10},
    \cite[Prop.~V.1.6]{Ne00})} 
If $\eset\not=C \subeq V$ 
is an open or closed convex subset, then the following assertions hold: 
\begin{description}
\setlength\itemsep{0em}
  \item[\rm(i)] $\lim(C) = \lim(\oline C)$ is a closed convex cone. 
  \item[\rm(ii)] $v \in \lim(C)$ if and only if 
$t_j c_j \to v$ for a net with 
$t_j \geq 0$, $t_j \to 0$ and $c_j \in C$. 
  \item[\rm(iii)] If $c \in C$ and $d \in V$ satisfy 
$c + \R_+ d \subeq C$, then $d \in \lim(C)$. 
\item[\rm(iv)] $H(C) = \{0\}$ if and only if $C$ contains no affine lines. 
\item[\rm(v)] $B(C)^\star = \lim(C)$ and $B(C)^\bot = H(C)$. 
\end{description}
\end{Lemma}

A function $f \: V \to \R \cup \{\infty\}$ is said to be {\it convex}  
if its epigraph 
\[ \epi(f) := \{(x,t) \in V \times \R\: f(x) \leq t \} \] 
is convex, and {\it lower semicontinuous} if its epigraph is closed
(cf.~\cite[Lemma~V.3.1]{Ne00}). For a convex
function $f \: D \to \R\cup \{\infty\}$ ($D \subeq V$ convex),
there is a unique
convex function $\oline f$ whose epigraph $\epi(\oline f)$
is the closure~$\oline{\epi(f)}$ (\cite[Prop.~V.3.7]{Ne00}).
If, conversely, $f$ is a closed convex function and 
$D_f := f^{-1}(\R)$, then
$f\res_{D_f^\circ}$ is continuous and its closure
coincides with $f$ (\cite[Prop.~V.3.2]{Ne00}).

\begin{lem} \mlabel{lem:a.2}  Suppose that
  $f$ is a lower semicontinuous convex function.
  If $f$ is bounded on a ray $v + \R_+ h \subeq D_f$, then
  \[ h \in \lim(D_f) \quad \mbox{ and }  \quad
    f(x + th) \leq f(x) \quad \mbox{ for all}\quad
    x \in D_f, t \geq 0.\]   
\end{lem}

\begin{prf} Our assumption implies the existence of $c \in \R$ for which 
  $(v + th, c) \in \epi(f)$ for all $t \geq 0$.
  This implies that $(h,0) \in \lim(\epi(f))$
  (Lemma~\ref{lem:limcone}(iii)).
  We conclude that, for all $x \in D_f$, we have
  \[ (x,f(x)) + \R_+ (h,0) \subeq \epi(f),\]
  which means that $f(x + th) \leq f(x)$ for all $t \geq 0$.
\end{prf}

\begin{lem} \mlabel{lem:2.3}
  Let $V$ be a finite-dimensional real vector space and
  $\mu$ a positive Borel measure on $V^*$ whose support spans~$V^*$.
We consider its  Laplace transform
  \[ \cL(\mu) \: D_\mu := \Big\{ v \in V \:  \int_{V^*}
    e^{-\alpha(v)}\, d\mu(\alpha) < \infty \Big\} \to \R,
    \quad \cL(\mu)(v) := \int_{V^*}    e^{-\alpha(v)}\, d\mu(\alpha). \]
  Then the following assertions hold:
  \begin{itemize}
  \item[\rm(a)] If $x \in D_\mu$ and $y \in \R$ are such that
    \begin{equation}
      \label{eq:decr}
\cL(\mu)(x + ty) \leq \cL(\mu)(x) \quad \mbox{ for all } \quad
t \geq 0,
    \end{equation}
    then $y \in \supp(\mu)^\star$.
  \item[\rm(b)] Let $y \in V$. If there exists some
    $x \in D_\mu$ with 
    \begin{equation}
    \label{eq:non-red}
    \cL(\mu)(x + ty) = \cL(\mu)(x) \quad \mbox{ for all }  \quad t\in \R,
  \end{equation}
  then $y =  0$.
\end{itemize}  
\end{lem}

\begin{prf} (a) Since the convex function $\cL(\mu)$ on $D_\mu$
  has a closed epigraph, the condition under (a)
  implies that $(y,0) \in \lim(\epi(\cL(\mu)))$ (Lemma~\ref{lem:limcone}(c)).
  The Monotone Convergence Theorem implies for $d \in \R$ that 
  \begin{align*}
\lim_{t \to \infty} e^{t d} {\cal L}(\mu)(x + ty) 
&= \lim_{t \to \infty} e^{t d} \int_{V^*} e^{-\alpha(x + ty)} \, d\mu(\alpha) 
= \lim_{t \to \infty} \int_{V^*} e^{t(d-\alpha(y))} \ e^{-\alpha(x)}\, d\mu(\alpha) \\
&=
  \begin{cases}
     0 & \text{ for } d < \inf \la \supp(\mu),y \ra \\ 
\int_{\alpha(y) = d} e^{-\alpha(v)}\, d\mu(\alpha) 
& \text{ for } d = \inf \la \supp(\mu), y\ra \\ 
\infty & \text{ for } d > \inf \la \supp(\mu),y \ra 
  \end{cases}
  \end{align*}
  (\cite[Rem.~V.4.12]{Ne00}).
  In view of \eqref{eq:decr}, this limit is $0$ for all $d < 0$,
  so that we must  have
  \[ \inf \la \supp(\mu),y\ra \geq 0, \quad \mbox{ i.e.,} \quad
    y \in \supp(\mu)^\star.\] 

  \nin (b) Applying (a) to $y$ and $-y$, it follows that
  $y \in \supp(\mu)^\star \cap -\supp(\mu)^\star = \supp(\mu)^\bot$. Since $\supp(\mu)$ spans $V^*$,
  we obtain $y= 0$.
\end{prf}

We continue with the setting of Lemma~\ref{lem:2.3}. 
For $x \in D_\mu$ and  $x^*(\alpha) = \alpha(x)$, the measure 
\begin{equation}
  \label{eq:mux}
  \mu_x := e^{-\log \cL(\mu)(x) - x^*} \cdot \mu = \frac{e^{-x^*} \mu}{\cL(\mu)(x)}
\end{equation}
is a probability measure on $V^*$. If $D_\mu$ has interior points in $V$
and $x \in D_\mu^\circ$, then the smoothness
of the Laplace transform on the open convex set~$D_\mu^\circ$  
implies that the expectation value of this measure exists and is given by 
\begin{equation}
  \label{eq:exp-val}
Q(x) := \frac{1}{\cL(\mu)(x)} \int_{V^*} \alpha e^{-\alpha(x)}\, d\mu(\alpha) 
= -\dd(\log\cL(\mu))(x) 
\end{equation}
(\cite[Prop.~V.4.6]{Ne00}). It is contained in
\begin{equation}
  \label{eq:cmudef}
  C_\mu := \oline\conv(\supp(\mu)) \subeq V^*.  
\end{equation}

\begin{prop}
  \mlabel{prop:I.9}
\begin{itemize}
\item[\rm(i)] The functions $\cL(\mu)$ and $\log(\cL(\mu))$ are 
convex and lower semicontinuous. 
If $C_\mu$ has interior points in $\g^*$,
then $\cL(\mu)$ and $\log\cL(\mu)$ are 
strictly convex on $D_\mu$.
\item[\rm(ii)] The function $\cL(\mu)$ is analytic on $D_\mu^\circ$ and
  has a 
holomorphic extension to the tube domain $D_\mu^\circ + i V$.
\item[\rm(iii)] Let $N_\mu := (C_\mu - C_\mu)^\bot$ be the linear subspace of all elements 
$x \in V$ for which $x^*$ is constant on~$\supp(\mu)$. 
Then $N_\mu + D_\mu = D_\mu$, the function $Q = - \dd(\log\cL(\mu))$ is constant on 
the $N_\mu$-cosets and factors through a function 
\[ \oline Q \: D_\mu/N_\mu \to C_\mu \subeq V^*.\] 
Its restriction to the relative
interior $D_\mu^\circ/N_\mu$ is a diffeomorphism onto a 
relatively open subset of $C_\mu$ in the affine subspace generated by $C_\mu$. 
If $C_\mu$ has interior points in $V^*$, then the bilinear form 
$\dd^2 (\log \cL(\mu))(x)$ is positive definite for all $x \in D_\mu^\circ$. 
\end{itemize}
\end{prop}

\begin{prf} (i) follows from \cite[Prop.~V.4.3, Cor.~V.4.4]{Ne00},
  and (ii) from \cite[Prop.~V.4.6]{Ne00}. 

\nin (iii) For $z \in N_\mu$ and $x \in D_\mu$, we have 
\[ \cL(\mu)(x + z) = e^{-z^*} \cL(\mu)(x) \quad \mbox{ and }  \quad 
\log \cL(\mu)(x + z) = - z^* + \log \cL(\mu)(x).\] 
This implies  $Q(x + z) = Q(x)$. For $x \in D_\mu^\circ$ and $y \in V$,
the argument in the proof of \cite[Prop.~V.4.6(iii)]{Ne00} shows that 
\[ \dd^2(\log\cL(\mu))(x)(y,y) \geq 0, \]
with equality if and only if $y \in N_\mu$, which is equivalent
to the linear function $v^*$ being $\mu$-almost everywhere constant 
(cf.~\eqref{eq:d2Z}). For $\oline y := y + N_\mu \in V/N_\mu$, we thus obtain 
\[ \la \dd \oline Q(\oline x)(\oline y), y \ra 
= \la \dd Q(x)(y), y \ra 
= - \dd^2(\log\cL(\mu))(x)(y,y) < 0 \quad \mbox{ if }\quad \oline y \not=0.\]  
This implies that 
$\dd Q(\oline x) \: V/N_\mu \to \aff(C_\mu)$ 
is injective, hence invertible 
because
\[ \dim(V/N_\mu) = \dim N_\mu^\bot = \dim(\aff(C_\mu)).\]
Therefore $\oline Q \: D_\mu^\circ /N_\mu \to C_\mu$ 
has open range in the affine~$\aff(C_\mu)$, 
and $\oline Q$ is a local diffeomorphism. 
To see that it is injective, we argue as in \cite[Lemma~1.3]{Ne19}
with $f := \log\cL(\mu)$. For $x,x +y \in D_\mu^\circ$ we have 
\[ \la \oline Q(\oline x + \oline y) - \oline Q(\oline x), y\ra 
= -\int_0^1 \dd^2 f(x + ty)(y,y)\, dt.\] 
If $\oline y\not= 0$, then $y \not\in N_\mu$, so that the right hand side 
is non-zero. Hence $\oline Q$ is injective. 
\end{prf}

With $N_\mu$ as in Proposition~\ref{prop:I.9}(iii), we now have: 
\begin{thm} \mlabel{thm:conv-lapl}
  {\rm(Convexity Theorem for Laplace Transforms)} 
  If $D_\mu \not=\eset$ is open, hence equal to $\Omega_\mu$,
  then   $\oline Q$ maps $\Omega_\mu/N_\mu$ diffeomorphically onto $C_\mu^\circ$. 
\end{thm}

\begin{prf} This follows from \cite[Thm.~V.4.9]{Ne00}
  because the domain $D_\mu$ of the closed convex function $\log\cL(\mu)$
  has no boundary points by assumption, hence satisfies the required
  essential smoothness condition by \cite[Lemma~V.3.18(v)]{Ne00}.
\end{prf}

Part (a) of the next proposition follows from
\cite[Thm.~II.1.7]{Bo96}, dealing more generally with tempered distributions.
We include the rather direct proof for the special case of tempered measures
and also add a very useful converse that can be used to verify
temperedness of measures. 

\begin{prop} \mlabel{prop:dom-lapl-temp} {\rm(Laplace transforms
    and temperedness)} 
  Let $V$ be a finite-dimensional real vector space
  and $\mu$ a positive Borel measure on $V^*$
  for which $C_\mu$ contains no affine lines,  i.e.,
  $B(C_\mu)$ has interior points {\rm(\cite[Prop.~V.1.16]{Ne00})}.
  Then the following assertions hold:
  \begin{itemize}
  \item[\rm(a)] If $\mu$ is tempered, then   $B(C_\mu)^\circ \subeq D_\mu$
    and there exists a $k \in \N$, such that, for every $z \in B(C_\mu)^\circ$
    \[ \limsup_{t \to 0+} \cL(\mu)(tz) t^k < \infty.\] 
  \item[\rm(b)] If there exists an $x \in B(C_\mu)^\circ$
    and $k \in \N$, such that
    \[ \limsup_{t \to 0+} \cL(\mu)(tx) t^k < \infty,\]
    then $\mu$ is tempered.
  \end{itemize}
\end{prop}

\begin{prf} We enlarge $V$ to the space
  $\tilde V = V \times \R$ and consider
  $\mu$ as a measure on $V^* \times \{1\} \subeq \tilde V^*$.
  Then
  \[ \cL(\mu)(x,t) = e^{-t} \cL(\mu)(x) \]
and $\aff(C_\mu) \subeq V^* \times \{1\}$ is an affine hyperplane
  not containing~$0$. This implies that
  \[ C := \cone(C_\mu) = \R_+ C_\mu \cup \big(\lim(C_\mu) \times \{0\}\big) \]
 is a pointed convex cone (\cite[Prop.~V.1.15]{Ne00})  and
 \[ B(C_\mu) = C_\mu^\star + \R (0,1) = C^\star + \R(0,1).\]

\nin (a) We have to show that $(C^\star)^\circ \subeq D_\mu$.

Let $z = (x,c) \in (C^\star)^\circ$. Then 
$C_1 := \{ \alpha \in C \:  \alpha(z) = 1 \}$ 
is a compact base of the cone $C$. We choose a norm $\|\cdot\|$ 
on $\tilde V$, such that its unit ball $B$ contains $C_1$, so that
\begin{equation}
  \label{eq:norm-esti}
  \alpha(z) \geq \|\alpha\| \quad \mbox{ for all } \quad
  \alpha \in C \supeq C_\mu.
\end{equation}

Since $\mu$ is tempered, by definition,
there exists a $k\in \N$ such that 
$\int_{V^*}  \frac{d\mu(\alpha)}{(1 + \|\alpha\|^2)^k} < \infty.$
For the Laplace transform of $\mu$ we now obtain
\begin{align*}
  \cL(\mu)(z)
  &= \int_{C_\mu} e^{-\alpha(z)}\, d\mu(\alpha) 
    \leq
    \int_{V^*} e^{-\|\alpha\|}\, d\mu(\alpha) 
    =     \int_{V^*} \underbrace{e^{-\|\alpha\|}(1 + \|\alpha\|^2)^k}_{\text{bounded}}
    \frac{d\mu(\alpha)}{(1 + \|\alpha\|^2)^k} < \infty.
\end{align*}
As $z \in (C^\star)^\circ$ was arbitrary,
this proves that $(C^\star)^\circ \subeq D_{\mu}$.

For $t > 0$, we further obtain
\begin{align*}
  \cL(\mu)(tz)
  &= \int_{C_\mu} e^{-t\alpha(z)}\, d\mu(\alpha) 
    \leq    \int_{V^*} e^{-t\|\alpha\|}\, d\mu(\alpha) 
    =     \int_{V^*}
    \underbrace{e^{-t\|\alpha\|}(1 + \|t\alpha\|^2)^k}_{\text{bounded}}
    \frac{(1 + \|\alpha\|^2)^k}{(1 + \|t\alpha\|^2)^k}
    \frac{d\mu(\alpha)}{(1 + \|\alpha\|^2)^k}. 
\end{align*}
As
\[ t^{2k} \frac{(1 + \|\alpha\|^2)^k}{(1 + t^2\|\alpha\|^2)^k}
  =   \frac{(t^{2k} + \|t \alpha\|^2)^k}{(1 + \|t\alpha\|^2)^k}
  \leq 1 \quad \mbox{ for } \quad 0 < t \leq 1,\]
it follows that $\limsup_{t \to 0+} \cL(\mu)(tz) t^{2k} < \infty$.

\nin (b) In view of the construction preceding the proof of (a),
we may w.l.o.g.\ assume that $\supp(\mu)$ is contained in a pointed
closed convex cone $C$ and that $x \in (C^\star)^\circ$. 
Our assumption implies the existence
of $c > 0$ and $\delta > 0$ such that
\[ \cL(\mu)(tx) \leq c t^{-k} \quad \mbox{ for }  \quad  0 < t \leq \delta.\]
For the measure $\mu_x := (x^*)_* \mu$ on $\R$, we have
$\cL(\mu_x)(t) = \cL(\mu)(tx),$ 
so that \cite[Prop.~4]{FNO25} implies that the measure $\mu_x$ on $\R$
is tempered, hence that there exists an $m \in \N$ with
\[   \int_{C} \frac{d\mu(\alpha)}{(1 + \alpha(x)^2)^m} 
  =   \int_\R \frac{d\mu_x(\alpha)}{(1 + \alpha^2)^m}
  < \infty.\] 

We choose a norm $\|\cdot\|$ on $V^*$ such that
$\|x^*\| \leq 1$, so that $|\alpha(x)| \leq \|\alpha\|$ for
$\alpha \in V^*$.
Then we have 
\begin{equation*}
\int_{V^*} \frac{d\mu(\alpha)}{(1 + \|\alpha\|^2)^m}
\leq \int_{C} \frac{d\mu(\alpha)}{(1 + \alpha(x)^2)^m}
= \int_\R \frac{d\mu_x(\alpha)}{(1 + \alpha^2)^m} < \infty.
\qedhere
\end{equation*}
\end{prf}

\subsection*{Entropy and Gibbs measures}

\begin{defn} \mlabel{def:entropy} Let $\lambda_M$
  be a positive Borel measure on the manifold~$M$,
  let $V$ be a finite-dimensional real   vector space, 
  and $\Psi \: M \to V^*$ be a smooth map. We write
  $\mu := \Psi_*\lambda_M$ for the push-forward measure on~$V^*$.
  
  \nin (a) A {\it related Gibbs measure} is a
  measure of the form
  \[ d \lambda_x(m)  = e^{-z(x) - \Psi(m)(x)}\, d\lambda_M(m) \quad
    \mbox{ with } \quad 
    z(x) =  \log \int_M e^{-\Psi(m)(x)}\, d\lambda_M(m).\]
We write  $\mu_x := \Psi_*\lambda_x$ for the corresponding
  probability measure on $V^*$. 

\nin (b) The {\it entropy} of the probability measure $\lambda_x$ with 
respect to the density function 
\[ p_x = e^{-z(x) - \la \Psi(\cdot),x \ra} \]
is defined by 
\begin{align}\label{eq:entropie-s}
 s(x) 
&:= - \int_M \log(p_x) \cdot p_x\, d\lambda_M 
= - \int_M \log(p_x) \, d\lambda_x \notag \\
&=  \int_{V^*} \alpha(x) + z(x)  \, d\mu_x(\alpha) 
=   Q(x)(x) + z(x). 
\end{align}
\end{defn}

\begin{thm}\mlabel{thm:b.3}
  {\rm(\cite[Thm.~(16.200)]{So97})} 
Let $\lambda_x$ be a Gibbs measure on $M$ related to the continuous map
$\Psi \: M \to V^*$ and the measure $\lambda_M$ on $M$.
Suppose that the expectation value 
\[ Q(x) = \int_M\Psi\, d\lambda_x = \int_{V^*}\alpha\, d\mu_x(\alpha) 
  \quad \mbox{ of } \quad \mu_x = \Psi_*\lambda_x \]
 exists.  Then the  $\lambda_M$-entropy  $s(x)$ exists and equals 
 \begin{equation}
   \label{eq:entropie-eq}
   s(x) = z(x) + Q(x)(x).
 \end{equation}
All other probability measures which are completely continuous with 
respect to $\lambda_M$ and with the same expectation value
$Q(x)$ have an entropy strictly less than~$s(x)$. 
\end{thm}

\section{Invariant probability measures for linear groups} 
\mlabel{app:coad-fin-vol}

Our starting point in this section is the 
Poincar\'e Recurrence Theorem~\ref{thm:PRT}.
We shall use it to derive that, if a connected Lie group
$G \subeq \GL(V)$ preserves a probability measure $\mu$ on $V$,
whose support spans $V$, then
the closure of $G$ is compact.
As a consequence, coadjoint orbits whose Liouville measure
is finite arise only from compact groups. But we shall see below,
that there are stronger conclusions concerning the
openness of the domain of Laplace transforms of invariant
(not necessarily finite) measures on $\g^*$. In particular, we
shall see, in the context of geometric temperatures
in Lie algebras, that $D_\mu \subeq \comp(\g)^\circ$ whenever
the measure $\mu$ spans $\g^*$.\\


\begin{thm} \mlabel{thm:PRT} {\rm(Poincar\'e Recurrence Theorem)} 
Let $(X,\Sigma, \mu)$ be a finite measure space and $f \: X \to X$ 
be a measure preserving Borel automorphism. Then, for any 
$E \in \Sigma$, the sets 
\[ E_+(f) := \{ x \in E \: (\exists N \in \N_0)(\forall n > N)\  f^n(x) \not\in E \}
  = E \setminus \bigcup_{N \in \N} \bigcap_{n > N} f^{-n}(E) \]
and $E_-(f) := E_+(f^{-1})$ have measure zero.
\end{thm}

This means that almost every point $x \in E$ 
returns to $E$ in the sense that there exists a strictly 
increasing sequence $(n_k)_{k \in \N}$ of natural numbers 
with $f^{n_k}(x) \in E$, and that the same holds for $f^{-1}$. 

\begin{prf} For the sake of completeness, we include a sketch of the
  simple proof (cf.~\cite[\S 1.29]{Na13}).
  As $f^{-1}$ also satisfies the assumption,
  it suffices to show that $\mu(E_+(f)) = 0$. 
  We consider the measurable subset
  \[ F := \{ x \in E \:  (\forall k \geq 1) \, f^k(x) \not\in E\}
    = E \setminus \bigcup_{k > 0} f^{-k}(E).\]
Then it is easily seen that  
the sequence $(f^n(F))_{n  \in \Z}$ is  pairwise disjoint.
Therefore the invariance and the finiteness of the measure imply that
$\mu(F) = 0$, so that $\bigcup_{k \geq 0} f^k(F) \supeq E_+(f)$ 
is also a $\mu$-null set .
\end{prf}

In the following lemma we shall use the multiplicative Jordan
decompositions $g= g_e g_h g_u$ of $g \in \GL(V)$, $V$ a finite-dimensional
real vector space. These are uniquely determined commuting factors, where
$g_e$ is elliptic (diagonalizable over $\C$ with eigenvalues of
absolute value $1$), $g_h$ is hyperbolic (diagonalizable
with positive eigenvalues), and $g_u$ is unipotent, i.e.,
$(g_u - \1)^N = 0$ for some $N \in \N$.

\begin{lem} \mlabel{lem:gotoinf} 
  Let $V$ be a finite-dimensional real vector space
  and $g\in \GL(V)$. We write $g = g_e g_h g_u$ for its multiplicative
  Jordan decomposition into elliptic, hyperbolic and unipotent factor.
Then the following assertions hold: 
\begin{itemize}
\item[\rm(a)]  If $v \in V$, then one of the sequences $g^nv$ or $g^{-n}v$
  eventually leaves every compact subset of $V$ if and only if
$v$ is not fixed by $g_h g_u$. 
\item[\rm(b)] If $\mu$ is a finite $g$-invariant Borel measure on~$V$, 
then $\supp(\mu) \subeq \Fix(g_h g_u)$. 
\end{itemize}
\end{lem}

\begin{prf} Since $g_e^\Z$ has compact closure in $\GL(V)$, there exists a
 $g_e$-invariant norm on $V$. 

  \nin  (a) Suppose that $v \in V$ is not fixed by $g_h g_u$, the trigonalizable
  Jordan component of $g$.
  Let 
  \[ V_\lambda(g_h) = \ker(g_h - \lambda \1) \]
  be the eigenspaces of the hyperbolic factor $g_h$ and recall that all
  eigenvalues are positive.

  \nin {\bf Step 1:} We consider $v \in V$
  that is not fixed by $g_u$ and the linear subspace
  \[ W := \Spann \{ g_u^n.v \:  n \in \N_0\} \subeq V, \]
  for which our assumption implies $\dim W > 1$.  
  Since $g_u - \1$ is nilpotent and non-zero on $W$, the Jordan Normal
  Form implies that $\dim (g_u-\1)^k W = \dim W - k$ for $k \leq \dim W$,
  so that $\oline W := W/(g_u-\1)^2 W$ is $2$-dimensional.
  The image $\oline v$ of $v$ in this space satisfies
  \[ (\oline g_u -\1) \oline v \not=0 \quad \mbox{ and } \quad 
    (\oline g_u -\1)^2 \oline v  =0,\]
  so that
  \[ \oline g_u^n.\oline v = \oline v + n (\1 - \oline g_u) \oline v
    \quad \mbox{ for } \quad n \in \Z.\]
  As this sequence is unbounded in both directions in $\oline W$,
  the same holds for the sequence $g_u^n.v$ in $V$.

\nin  {\bf Step 2:} If there exists an eigenvalue
  $\lambda > 1$, then $v$ has a non-zero component $v_\lambda$
  in this eigenspace, which is a generalized eigenspace of
  $g_h g_u$. Then 
  \[ \|g^n.v_\lambda \| = \lambda^n \|g_u^n.v_\lambda\|,\]
  and if $g_u$ does not fix $v_\lambda$, then Step 1 implies that
  $\|g_u^n.v_\lambda\| \to \infty$; otherwise 
  $g_u^n.v_\lambda =v_\lambda$ for all $n \in \Z$.
  In both cases $\lambda > 1$ implies that 
  $\|g^n.v_\lambda \| \to \infty$. 

  If there exists an eigenvalue $\lambda < 1$ of $g_h$,
  then the same argument   applies to $g^{-1} = g_e^{-1} g_h^{-1} g_u^{-1}$
  and shows that $\|g^{-n}.v_\lambda\| \to \infty$.

  \nin  {\bf Step 3:} In view of Steps 1 and 2, a necessary condition
  for neither $g^n.v$ nor $g^{-n}.v$ to tend to infinity is that,
  on the cyclic subspace generated by $v$, we have 
  $g_h = \1$, i.e., all its eigenvalues are $1$, and that
  $g_u = \1$ as well. This means that $g_h g_u.v = v$.
  If, conversely, this condition is satisfied,
  then the sequence $g^n.v = g_e^n.v$ is bounded. This completes
  the proof of (a).
  
\nin (b)  If $v \in V$ with $(g_h g_u).v \not=v$, 
then either $g^nv \to \infty$ or $g^{-n} v\to \infty$ by (a). 
We conclude that, for every compact subset $C \subeq V\setminus \Fix(g_h g_u)$, 
no point $v \in C$ is recurrent for $g$ and $g^{-1}$. 
By the Poincar\'e Recurrence Theorem (Theorem~\ref{thm:PRT}),
the set of all $v \in C$ with $g^n.v \to \infty$ has measure zero,
and so does the set of all $v \in C$ with $g^{-n}.v \to \infty$.
This shows that $\mu(C) = 0$, and hence that 
$\mu(V \setminus \Fix(g)) = 0$ 
because the open set $V\setminus \Fix(g)$ is a countable union of compact 
subsets.
We conclude that $\supp(\mu) \subeq \Fix(g_hg_u)$. 
\end{prf}

\begin{thm} {\rm(Compactness Theorem)} \mlabel{thm:measlemb}
Let $V$ be a finite dimensional real vector space. 
\begin{itemize}
\item[\rm(a)] If $\mu$ is a finite positive Borel measure on $V$ whose
  support spans $V$, then its stabilizer group
  $\GL(V)^\mu := \{g \in \GL(V) \: g_*\mu = \mu\}$ is closed and
  has the property that  all its elements are elliptic, i.e., generate relatively
  compact subgroups of $\GL(V)$.
\item[\rm(b)] If $G\subeq \GL(V)$ is a closed subgroup, such that
  all elements of $G$ are elliptic, then $G$ is compact. 
\end{itemize}
\end{thm}

\begin{prf} (a) For $\xi \in C_c(V)$ the function
  \[ \GL(V) \to \R, \quad g \mapsto
    \int_V \xi(v)\, d(g_*\mu)(v) 
    \int_V \xi(gv)\, d\mu(v) \]
  is continuous, so that the stabilizer $\GL(V)^\mu$ is a
  closed subgroup of $\GL(V)$.  \begin{footnote}
    {By \cite[Thm.~3.2.4]{Zi84} and an embedding of $\GL(V)$ into
      $\PGL(V \oplus \R)$ one can even show that this group is algebraic.}
  \end{footnote}
  By Lemma~\ref{lem:gotoinf}(b), all elements
$g \in \GL(V)^\mu$ are elliptic, i.e., $g = g_e$.

\nin (b) As $[\g,\rad(\g)]$ consists of nilpotent elements 
(\cite[\S 5.4.2]{HN12}), its exponential 
image consists of unipotent elements, hence is trivial. 
Therefore $\rad(\g)$ is central in $\g$, which means that $\g$ is reductive. 
The Cartan decomposition shows that,
any non-compact simple real Lie algebra contains non-zero
$\ad$-diagonalizable elements, and their exponential image
is hyperbolic. As this is excluded,
all simple ideals of $\g$ are compact,
and this entails that $\g$ is a compact Lie algebra.
We now have $\g = \fz(\g) \oplus [\g,\g]$ with $[\g,\g]$
compact semisimple. Then the Lie group $\la \exp [\g,\g]\ra$
is compact (\cite[Thm.~12.1.17]{HN12}). That $\exp(\fz(\g))$ also has compact 
closure  follows from the fact that, for each $x \in \fz(\g)$,
$\exp(\R x) = \exp([0,1]x) \exp(\Z x)$
has compact closure because $\exp(x)$ is elliptic.

This implies that the identity component $G_e$ is compact. Moreover,
for every $g \in G$, the closed subgroup $\oline{g^\Z} \subeq G$
is compact, hence has at most
finitely many connected components. Therefore $\pi_0(G) := G/G_e$ is a
torsion group.

As $G_e$ is compact, it is also Zariski closed,
so that its normalizer $N\subeq \GL(V)$
is a real algebraic group containing $G$.
In $N$, the identity component $G_e$ is a normal algebraic
subgroup, so that $H := N/G_e$ is an affine algebraic group,
hence has a realization as an algebraic subgroup of some
$\GL_d(\R)$. In the Lie group topology of $H$, the image of
$G$ is discrete and isomorphic to $\pi_0(G)$, hence a discrete
torsion group. Therefore the Corollary in \cite{Wa74}
implies that $\pi_0(G)$ is finite. This proves that $G$ is compact.

Instead of Wang's paper, we can also use
\cite[Lemma~2]{Le76}, asserting that every torsion subgroup of a connected
Lie group is contained in a maximal compact subgroup. 
It implies that the image of $G \cap N_e$ has compact closure in~$H$,
but since it is also discrete, it is finite. As $N$ is algebraic,
the group $\pi_0(N)$ is finite (\cite[Prop.~2.3]{BHC62}), so that 
$G \cap N_e$ has finite index in $G$, and therefore $G$ is compact. 
\end{prf}

We thank Yves Cornulier for the reduction argument in the preceding
proof, using algebraic groups and for pointing out the following
  example of a linear group $\Gamma$ which is not closed
  and not relatively compact, although all of its elements are elliptic.

\begin{ex}   We conside the group
  \[ G := \C^2 \rtimes \SU_2(\C) \subeq \Aff(\C^2) \subeq \GL_3(\C).\]
  Then every element $(v,u) \in G$ with $u \not=\1$ is conjugate to
  an element of $\SU_2(\C)$ because $u$ has no non-zero fixed points,
  so that any $w \in \C^2$ with $uw-w = v$ conjugates $(v,u)$ to
  $(v + w -uw,u) = (0,u)$. Therefore the complement
  of the normal abelian subgroup $A := \C^2 \times \{\1\}$
  of $G$ consists of elliptic   elements.

  Next we recall that the Lie algebra $\g = \C^2 \rtimes \su_2(\C)$
  is generated by two elements $a,b$ (\cite[Thm.~6]{Ku51}).
  In fact, let $x,y \in \su_2(\C)$ be two generators and consider
  elements of the form $a = (0,x), b = (v,y) \in \g$. Since $\ad x$
  has on $\su_2(\C)$ different eigenvalues than on $\C^2$, it easily
  follows that $a$ and $b$ generate the perfect Lie algebra~$\g$. 
  Kuranishi shows that these elements 
  can be chosen in such a way that the projections of 
  $g := \exp(a)$ and $h := \exp(b)$ to $\SU_2(\C)$ generate a
  free subgroup (\cite[Thm.~8]{Ku51})
  and that the group $\Gamma$ generated by $g$ and $h$
  is dense  in $G$. Freeness of the projection
  to $\SU_2(\C)$ then implies that 
  $\Gamma \cap A =  \{e\}$. Therefore $\Gamma$
  consists of elliptic elements, but its closure $G$ does not. 
\end{ex}

We now describe an alternative argument for the compactness of the
stabilizer of a probability measure in $\GL(V)$,
using Shalom's variant of F\"urstenberg's Lemma
(cf.\ \cite[p.~171]{Sh98}, \cite[Lemma~3]{Fu76}),
which deals with measures on projective spaces. 

\begin{lem} {\rm(F\"urstenberg--Shalom Lemma)} \mlabel{lem:fuerst}
  Let $k$  be a locally compact, non-discrete field
  and $H \subeq \GL_n(k)$ be an algebraic subgroup, 
  $\mu$ a probability measure on the projective space
  $\bP_{n-1}(k) = \bP(k^n)$, and $H^\mu$ the stabilizer group
  of $\mu$ in $H$. Then there exist finitely many linear subspaces
    $V_1,\ldots, V_\ell \subeq k^n$ such that
    \[ \mu([V_1] \cup \cdots \cup [V_\ell]) = 1,\]
    and an algebraic normal cocompact subgroup $H_S \subeq H^\mu$ which
    fixes every point in
    \[ S := [V_1] \cup \cdots \cup [V_\ell].\] 
\end{lem}

Shalom concludes from this lemma, that, if $H \subeq \GL_n(k)$
is semisimple algebraic and 
$G \subeq H$ amenable and Zariski dense in $H$, then $G$ has compact closure.
In our context, it provides the following
more direct, but less informative,
proof of the combination of (a) and (b) in
the Compactness Theorem:

\begin{prf} Let $\mu$ be a probability measure on $V$
  whose support spans $V$.
  We have to show that, in the algebraic group $H := \GL(V)$,
  the stabilizer $H^\mu$ of $\mu$ is compact. 
To this end, we consider the enlarged space $\tilde V := V \times \R$
  and embed $V$ as the affine subspace $A := V \times \{1\}$.
  Then $[A] \subeq \bP(\tilde V)$ is a dense open subset and we consider
  $\mu$ as a probability measure on~$A$. Further, 
\[ H = \GL(V) \into \PGL(\tilde V),\quad   g \mapsto [g \oplus 1] \]  
is a closed embedding onto an algebraic subgroup.
Let $V_1, \ldots, V_\ell$ be as in Lemma~\ref{lem:fuerst}.
Then $\mu$ is
supported in the union of the  affine subspaces $V_j \cap A$ of $A \cong V$.
Our assumption now implies that the affine spaces $V_j \cap A$
generate $V$ as a linear space. Therefore the pointwise stabilizer of
this union in $\GL(V)$ is trivial, and thus
F\"urstenberg's Lemma, as stated in \cite[p.~171]{Sh98},
implies that the stabilizer $H^\mu$ of $\mu$ is compact.
\end{prf}

\subsection*{Applications to coadjoint orbits}

\begin{cor} \mlabel{cor:measlem}
  Let $G$ be a finite-dimensional Lie group
  with Lie algebra $\g$ and $\mu$ an $\Ad^*(G)$-invariant
  Borel measure on $\g^*$ whose support spans $\g^*$.
  Then, for every $x \in \g$ with $\cL(\mu)(x) < \infty$, we
  have $\ker(\ad x) \subeq \comp(\g)$ and $x \in \comp(\g)^\circ$
  {\rm(cf.~Definition~\ref{def:2.1}(b))}. 
\end{cor}

\begin{prf} If $H_x(\alpha) = \alpha(x)$ denotes the evaluation
  functional on $\g^*$, then our assumption implies that
  $e^{-H_x} \mu$ is a finite positive Borel measure on $\g^*$
  invariant under the action of the group $\Ad(G^x)$.
  Theorem~\ref{thm:measlemb}
  thus implies that $\Ad(G^x)$ is relatively compact,
  so that   $\fz_\g(x) = \ker(\ad x) = \L(G^x)$ is compactly embedded,
  hence contained in $\comp(\g)$.
  That this is equivalent to $x \in \comp(\g)^\circ$
  follows from \cite[Lemma~VII.1.7(c)]{Ne00}. 
\end{prf}

\begin{cor} 
\mlabel{cor:d.1} 
Let $\cO_\lambda \subeq \g^*$ be a coadjoint orbit 
spanning $\g^*$. Then the following are equivalent: 
\begin{itemize}
\item[\rm(a)] The Liouville measure $\mu_\lambda$ is finite. 
\item[\rm(b)] $\g$ is  a compact Lie algebra. 
\item[\rm(c)] $\cO_\lambda$ is compact. 
\end{itemize}
\end{cor}

\begin{prf} (b) $\Rarrow$ (c): For a compact Lie algebra $\g$, 
the adjoint group is compact, so that all coadjoint orbits are compact. 
  
\nin (c) $\Rarrow$ (a) follows from the fact that the Liouville measure 
is finite on compact subsets. 
  
\nin (a) $\Rarrow$ (b): This is the non-trivial part. It follows 
from Corollary~\ref{cor:measlem}. 
\end{prf}

\begin{cor} \mlabel{cor:fin-vol}
  If $\cO_\lambda \subeq \g^*$ is  a coadjoint orbit
  of finite Liouville measure, then the quotient
  $\g/\cO_\lambda^\bot$ is a compact Lie algebra.   
\end{cor}

\begin{prf} If $\mu_\lambda$ is finite, then
  Corollary~\ref{cor:d.1} applies to the quotient Lie algebra
  $\g/\cO_\lambda^\bot$, whose dual is spanned by $\cO_\lambda.$ 
\end{prf}

\section{Admissible Lie algebras} 
\mlabel{sec:2}

Let $G$ be a connected Lie group with Lie algebra $\g$.
Subsection~\ref{subsec:1.1} introduces admissible Lie algebras 
The key tool to describe the fine structure of admissible Lie algebras 
is the root decomposition with respect to a compactly embedded 
Cartan subalgebra (Subsection~\ref{subsec:1.2}). 
In Subsection~\ref{subsec:1.3} we briefly recall from
\cite{Ne96b} and \cite{Ne00} how invariant 
convex functions relate to the root decomposition. 
The structure of admissible Lie algebras is 
described in Subsection~\ref{subsec:1.4}.
We conclude this section with
the proof of Theorem~\ref{thm:5.9} in Subsection~\ref{subsec:4.5}.
It draws from $D_\mu \not=\eset$ for an invariant measure~$\mu$ 
on~$\g^*$, whose support spans $\g^*$, the conclusion that $\g$
is admissible and $\supp(\mu) \subeq W_{\rm min}^\star$ for a suitable
positive system.

\subsection{From invariant convex functions to admissible Lie algebras}
\mlabel{subsec:1.1}

\begin{defn} \mlabel{def:2.1}
  (a) A Lie algebra $\g$ is said to be 
{\it admissible} if it contains a non-empty open invariant convex 
subset not containing affine lines. 

\nin (b) An element $x \in \g$ is said to be {\it elliptic}, or {\it compact}, 
if the one-parameter subgroup $e^{\R \ad x} \subeq \Aut(\g)$ has compact closure, 
i.e., if $\ad x$ is semisimple with purely imaginary spectrum. 
We write $\comp(\g)$ for the set of compact elements of $\g$. 

\nin (c) A subalgebra $\fs \subeq \g$ is said to be {\it compactly embedded} 
if the subgroup generated by $e^{\ad \fs} \subeq \Aut(\fg)$ has compact closure. 
\end{defn}

\begin{rem} \mlabel{rem:hermitian} (a)
  A simple Lie algebra $\g$ is admissible if and only if it either 
  is compact or hermitian, i.e.,
  a maximal compactly embedded subalgebra $\fk \subeq \g$ 
  has non-trivial center (cf. \cite[Prop.~VII.2.14]{Ne00}).
  For compact Lie algebras, admissibility
  follows from the existence of an invariant norm, so that the
  balls are invariant and contain no affine lines.
  For hermitian Lie algebras, admissibility
  follows from the existence a pointed   generating invariant cone.
  This is a consequence of the Kostant--Vinberg Theorem  on
  the existence of invariant cones in representations
  (cf.\ \cite{Vi80}).
  We refer to \cite[Thm.~VII.25]{HN93} for a rather direct argument.
  Here is a list of the simple hermitian Lie algebras:
  \[ \su_{p,q}(\C),  \qquad
    \so_{2,d}(\R),\ d > 2, \qquad
    \sp_{2n}(\R), \qquad
    \so^*(2n), \qquad \fe_{6(-14)}, \quad
    \fe_{7(-25)}.\]

\nin (b) A reductive Lie algebra $\g$ is admissible if and only if all its 
simple ideals are admissible (\cite[Lemma~VII.3.3]{Ne00}). 
\end{rem}

Let $\eset \not= \Omega \subeq \g$ an $\Ad(G)$-invariant convex subset and 
$f \: \Omega \to \R$ a convex function which is invariant under the 
adjoint action, i.e., constant on adjoint orbits.
Then the subset 
 \begin{equation}
   \label{eq:nf}
\fn_f := \{ x \in \g \: 
x + \Omega = \Omega,\ (\forall y \in \Omega)\, f(x + y) = f(y)\} 
 \end{equation}
 is an ideal of $\g$ because $f$ is $\Ad(G)$-invariant and
 $\Ad(G)$-invariant linear subspaces of $\g$ are ideals.
 The function $f$ is constant on the cosets 
$x + \fn_f$. Hence $f$ factors through a convex function on the convex 
subset $\Omega/\fn_f$ in the quotient Lie algebra $\g/\fn_f$. 
We call $f$ {\it reduced} if  $\fn_f = \{0\}$. So the following 
proposition asserts that the existence of 
reduced convex functions implies that $\g$ is admissible. 

\begin{prop} \mlabel{prop:1.2} Suppose that $\fn_f = \{0\}$.
  \begin{itemize}
  \item[\rm(a)] If $\Omega$ is open, then the following assertions hold:
    \begin{itemize}
    \item     $\g$ is admissible,
  \item For $c \in \R$, the open subset
    $\Omega_c := \{ x \in \Omega \: f(x) < c\}$ contains no affine lines. 
  \item $\Omega \subeq \comp(\g)$ {\rm(cf.~Definition~\ref{def:2.1})}.
    \end{itemize}
  \item[\rm(b)] Suppose that $f$ is closed, i.e., $\epi(f)$ is closed in
    $\g \oplus \R$. Then, for each $c \in \R$, the subset 
    $D_c := \{ f \leq c\}$ is closed and convex,  not containing
    affine lines.
  \item[\rm(c)] If $f$ is closed and $\g = \Spann D_f$,
    then $\g$ is admissible. 
  \end{itemize}
\end{prop}

\begin{prf}
Let $c \in \R$ be such that 
the open subset $\Omega_c := \{ x \in \Omega \: f(x) < c\}$ 
is non-empty. As $f$ is continuous and invariant and 
$\Omega$ is invariant, the subset $\Omega_c$ is an open 
convex invariant subset of $\g$. 
If $x + \R y \subeq \Omega_c$ is an affine line, then 
$f$ is bounded from above on this line, hence constant,
as a bounded convex function. 
Lemma~\ref{lem:limcone} implies that $\Omega_c + \R y = \Omega_c$,
hence $y \in H(\Omega)$ by Lemma~\ref{lem:limcone}(iii). 
We further obtain $(y,0) \in H(\oline{\epi(f)})$ (see \eqref{eq:limcone}),
so that $f$ is bounded, hence constant, 
on all affine lines $z + \R y$, $z \in \Omega$. Therefore 
$y \in \fn_f = \{0\}$, and we conclude that
$\Omega_c$ contains no affine lines. Therefore $\g$ is admissible.

For $c \in \R$, the inclusion $\Omega_c \subeq \comp(\g)$ 
now follows from \cite[Prop.~VII.3.4(e)]{Ne00}, so that 
$\Omega = \bigcup_{c \in \R} \Omega_c \subeq \comp(\g).$


\nin (b) Suppose that $D_c \not=\eset$. Note that this subset
is closed and $\Ad(G)$-invariant. Any affine line 
$x + \R y \subeq D_c$ leads with the same argument as under (a) to
$f$ being constant on all lines $z + \R y$, $z \in D_c$,
and we conclude as above that $y = 0$.

\nin (c) The assumption that $D_f$ spans $\g$
implies that, either $D_f$ has interior points in $\g$ or
in a proper affine hyperplane $\aff(D_f)$. Restricting
$f$ to the relative interior $D_f^\circ \subeq \aff(D_f)$, we obtain
a continuous function. Hence, for any $x_0 \in D_f^\circ$
and $c > f(x_0)$, the sublevel set $D_c$ contains a neighborhood
of $x_0$ in $\aff(D_f)$, hence also spans $\g$.
Since the invariant closed convex subset $D_c$ contains no affine lines,
\cite[Lemma~VII.3.1]{Ne00} and the definition of
admissibility imply that $\g$ is admissible. 
\end{prf}

\subsection{Root decomposition} 
\mlabel{subsec:1.2}

If the subset $\comp(\g)$ of compact elements in the Lie algebra
$\g$ has interior points, such as
in the context of Proposition~\ref{prop:1.2}, 
then \cite[Thm.~VII.1.8]{Ne00} implies the existence of a {\it compactly 
embedded Cartan subalgebra} $\ft \subeq \g$, i.e., 
$\ft$ is abelian, compactly embedded and 
coincides with its own centralizer: 
\[ \ft  = \fz_\g(\ft) := \{ x \in \g \: [x,\ft] = \{0\}\}.\] 
Then we have the root decomposition 
\[ \g_\C = \ft_\C \oplus \bigoplus_{\alpha \in \Delta} 
\g_\C^\alpha, \quad \mbox{ where } \quad 
 \g_\C^\alpha := \{ x \in \g_\C \: (\forall h \in \ft_\C)\, [h,x]= \alpha(h)
x\} \] 
and 
\[ \alpha(\ft) \subeq i \R \quad \mbox{ for every  root } 
\quad \alpha \in \Delta :=  
\{ \alpha \in \ft_\C^* \setminus \{0\} \: 
\g_\C^\alpha\not= \{0\}\}.\]
For $x + iy \in \g_\C$ we put $(x+ iy)^* := -x + iy$, so that 
\[ \g = \{ x \in \g_\C \: x^* = -x\}.\] 
We then have $x_\alpha^* \in \g_\C^{-\alpha}$ for 
$x_\alpha \in \g_\C^\alpha$. 
We call a root $\alpha \in \Delta$ 
\begin{itemize}
\item {\it compact}, if 
there exists an $x_\alpha \in \g_\C^\alpha$ with 
$\alpha([x_\alpha, x_\alpha^*]) > 0$, and 
\item {\it non-compact}, if 
there exists a non-zero $x_\alpha \in \g_\C^\alpha$ with 
$\alpha([x_\alpha, x_\alpha^*]) \leq 0$.
\item {\it solvable}, if it occurs in the root space decomposition
  of $\fr_\C$, where $\fr \trile \g$ is the maximal solvable ideal of $\g$.
\item {\it semisimple}, if it occurs in the root space decomposition
  of $(\g/\fr)_\C$. 
\end{itemize}
We write 
$\Delta_k$, $\Delta_p$, $\Delta_r,$ resp., $\Delta_s  \subeq \Delta$
for the subset of compact, non-compact, solvable,
semisimple roots (cf.\ \cite[Thm.~VII.2.2]{Ne00}). Then $\Delta_k \subeq \Delta_s$ and we
also write $\Delta_{p,s} := \Delta_p \setminus \Delta_s$ for the
semisimple non-compact roots.

If $\alpha$ is compact, then 
$\dim \g_\C^\alpha = 1$ and there exists a unique element 
$\alpha^\vee \in i \ft \cap [\g_\C^\alpha, \g_\C^{-\alpha}]$ 
with $\alpha(\alpha^\vee) = 2$. 
The linear endomorphism 
\[ r_\alpha \: \ft \to \ft,\quad  r_\alpha(x) 
:= x - \alpha(x)  \alpha^\vee 
= x  + (i \alpha)(x) i\alpha^\vee \] 
is called the corresponding reflection and 
\[ \cW_{\fk} := \cW(\fk,\ft) := \la r_\alpha \: \alpha \in \Delta_k \ra 
\subeq \GL(\ft) \] 
is called the {\it Weyl group}.

A subset $\Delta^+\subeq \Delta$ is called a {\it positive system} 
if there exists an $x_0 \in \ft$ with $\alpha(x_0) \not=0$ for every 
$\alpha \in \Delta$ and 
\[ \Delta^+  = \{ \alpha \in \Delta \: i\alpha(x_0) > 0 \}.\] 
A positive system is said to be {\it adapted} 
if $\Delta_p^+ := \Delta^+ \cap \Delta_p$ is invariant under 
$\cW_{\fk}$ (cf.~\cite[Prop.~VII.2.12]{Ne00}). 
Any such system specifies two $\cW_{\fk}$-invariant 
convex cones in $\ft$, which are relevant for invariant 
convex sets and functions (\cite[Def.~VII.3.6]{Ne00}): 
\begin{equation}
  \label{eq:maxcon}
C_{\rm min} :=  C_{\rm min}(\Delta_p^+) 
:= \oline\cone(\{ i [x_\alpha, x_\alpha^*] \: 
x_\alpha \in \g_\C^\alpha, \alpha \in \Delta_p^+\})\subeq \ft 
\end{equation} 
and 
\begin{equation}
  \label{eq:maxcone2}
C_{\rm max} := C_{\rm max}(\Delta_p^+) 
:= \{ x \in \ft \: (\forall \alpha \in \Delta_p^+) \ i\alpha(x) \geq 0\}.
\end{equation}

We collect the key results concerning invariant cones
in the following theorem: 

\begin{thm} \mlabel{thm:invcon-adm} Let  $\g$ be admissible
  {\rm(Definition~\ref{def:2.1}(a))} and
  $\Delta^+ \subeq \Delta$ be an adapted positive system
  with $C_{\rm min} \subeq C_{\rm max}$. Then the following assertions hold:
  \begin{itemize}
  \item[\rm(a)] $W_{\rm max} = \oline{\Ad(G)C_{\rm max}}$ is a closed convex invariant cone with $W_{\rm max}^\circ = \Ad(G) C_{\rm max}^\circ \subeq \comp(\g)$. 
  \item[\rm(b)] For $x \in C_{\rm max}^\circ$,   we have 
  \[\oline{\conv(\Ad(G)x)} =
    \{ y \in \g \: p_\ft(\Ad(G)y)\subeq \conv(\cW_{\fk} x) + C_{\rm min}\}
    \subeq W_{\rm max}^\circ,\]
  where $p_\ft \: \g \to \ft$ is the projection with kernel
$[\ft,\g]$. 
  \item[\rm(c)]
For $x \in W_{\rm max}^\circ$, we have 
\[   W_{\rm min} := \{ y \in \g \:
  p_\ft(\Ad(G)y) \subeq C_{\rm min} \}
  = \lim \big(\oline{\conv(\Ad(G)x)}\big) \subeq W_{\rm max}.\]
In particular, this cone does not depend on $x$. 
\item[\rm(d)] $W_{\rm max} \cap \ft = C_{\rm max}$ and
  $W_{\rm min} \cap \ft = C_{\rm min}$. 
  \end{itemize}
\end{thm}

\begin{prf} (a) follows from Prop.~VIII.3.7 and Lemma~VIII.3.9 in \cite{Ne00}. 

\nin (b) \cite[Thm.~VIII.3.18]{Ne00}; \nin (c) \cite[Lemma~VIII.3.27]{Ne00};
  (d) \cite[Lemma~VIII.3.22, 27]{Ne00};  
\end{prf}

\subsection{Invariant convex functions} 
\mlabel{subsec:1.3}

We now refine the conclusions from Proposition~\ref{prop:1.2}
by using the cones $W_{\rm min}$ and $W_{\rm max}$. 
We show that reduced invariant convex functions live on domains in 
$W_{\rm max}^\circ$ for some adapted positive system,
and that these functions are decreasing in the direction 
of the corresponding cone~$W_{\rm min}$. 

\begin{prop} \mlabel{prop:inv-fct}
  Let $\Omega \subeq \g$ be an open convex subset and
  $f \:  \Omega \to \R$ an  invariant convex function with 
  $\fn_f = \{0\}$. 
  Then $\g$ is admissible and
  contains a compactly embedded Cartan subalgebra~$\ft$, 
  and there exists an  adapted positive system $\Delta^+$
  with $C_{\rm min} \subeq C_{\rm max}$, such that 
  \begin{itemize}
  \item[\rm(a)] $\Omega \subeq W_{\rm max}^\circ$,  and 
  \item[\rm(b)]  $f(x + y) \leq f(x)$ for $x \in \Omega$ and $y\in W_{\rm min}$.
  \end{itemize}
 The set $\Delta_p^+$ of positive non-compact roots is uniquely
  determined by~$f$.
\end{prop}

\begin{prf} First, $\Omega \subeq \comp(\g)$
  follows from Proposition~\ref{prop:1.2}. The existence of interior
  points in $\comp(\g)$ implies the existence of a compactly
  embedded Cartan subalgebra (\cite[Thm.~VII.1.8]{Ne00}). 
Next we derive from \cite[Thm.~VII.3.8]{Ne00} the existence of a uniquely 
determined adapted positive system $\Delta^+$, such that, for every $c \in \R$ 
and $\Omega_c = \{ x \in \Omega \: f(x) < c\}$, we have 
\begin{equation}
  \label{eq:wmin-incl}
 W_{\rm min} \subeq \lim(\Omega_c) \quad \mbox{ and } \quad 
 \Omega_c \subeq W_{\rm max}.
\end{equation}

In fact, the Sandwich Theorem~\cite[Thm.~VII.3.8]{Ne00} shows
that
\[ \Omega_c\cap \ft\subseteq C_{\rm max} \quad \mbox{ and } \quad
  C_{\rm min}  \subeq \lim(\Omega_c\cap \ft).\]
Then we use Theorem~\ref{thm:invcon-adm} and
the definition of $W_{\rm min/max}$ to get \eqref{eq:wmin-incl}.
As a consequence of 
$W_{\rm min} \subseteq \lim(\Omega_c)$,
we get $x+\R^+ y\subeq \Omega_c$ for $x\in \Omega_c$ and $y\in W_{\rm min}$.
We thus obtain with Lemma~\ref{lem:a.2} that
\begin{equation}
  \label{eq:fdecr}
\Omega = \bigcup_{c \in \R} \Omega_c \subeq W_{\rm max}^\circ  
\quad \mbox{ and } \quad f(x + y) \leq f(x) \quad \mbox{ for } \quad 
x \in \Omega,\   y\in W_{\rm min} 
\end{equation}
(\cite[Thm.~VII.3.8]{Ne00}).
The uniqueness of $\Delta_p^+$ follows from the fact that
the open convex cone $W_{\rm max}$ with
$W_{\rm max} \cap \ft = C_{\rm max}$  determines
$\Delta_p^+$ as the subset $\{ \alpha \in \Delta_p \: i\alpha(C_{\rm max})
\subeq [0,\infty)\}$. \end{prf}

\subsection{Structure of admissible Lie algebras} 
\mlabel{subsec:1.4}

The structure of admissible Lie algebras is particularly well understood
in terms of a decomposition that goes back to K.~Spindler
(cf.~\cite{Sp88} and the notes to \S VII.2 in \cite{Ne00}).
The following theorem follows from
\cite[Thms.~VIII.2.7, VIII.2.26, Prop.~VIII.2.9]{Ne00}: 

\begin{thm} \mlabel{thm:spind}
A Lie algebra $\g$ with compactly embedded Cartan subalgebra 
$\ft$ is admissible if and only if it has a $\ft$-invariant semidirect 
decomposition $\g = \fu \rtimes \fl,$ 
where $\fu = \fz(\g) \oplus V$ is $2$-step nilpotent with 
\begin{itemize}
\item[\rm(S1)] $V = [\fl, \fu] = [\ft,\fu]$ and $[V,V] \subeq \fz(\g)$. 
\item[\rm(S2)] $\fl$ is reductive admissible
  with $\fz(\g) \cap \fl = \{0\}$. 
\item[\rm(S3)] There exists an adapted positive system $\Delta^+$ with
  $C_{\rm min} \subeq C_{\rm max}$ and a
  linear functional $\lambda_\fz \in \fz(\g)^*$ such that, 
for every non-zero $x_\alpha \in \fg_\C^\alpha = \fu_\C^\alpha$, $\alpha \in \Delta_r^+$, 
we have $\lambda_\fz(i[x_\alpha, x_\alpha^*]) > 0.$
\end{itemize}
\end{thm}

We call the decomposition from the preceding
theorem a {\it Spindler decomposition} of $\g$.  Then 
\begin{equation}
  \label{eq:omegaf}
  \Omega(v,w) := \lambda_\fz([v,w])
\end{equation}
defines on $V$ a symplectic form, which, in view of (S3), satisfies
\[ \Omega([x,v], v) > 0 \quad \mbox{ for } \quad 
 x \in C_{\rm max}^\circ, 0 \not= v \in V.\] 
Note that $H_x(v) := \frac{1}{2}\Omega([x,v], v)$ 
is the Hamiltonian function corresponding to the Hamiltonian flow 
on $(V,\Omega)$ generated by~$\ad x$ (\cite[Prop.~A.IV.15]{Ne00}).  
For details we refer to Section~VIII.2 in \cite{Ne00}, and 
in particular to \cite[Thm.~VIII.2.7]{Ne00}; see also \cite{NO22}. 

\subsection{From finiteness of Laplace transforms to  admissibility} 
\mlabel{subsec:4.5} 

The following theorem implies in particular that, whenever we have
a momentum map of a Hamiltonian action whose image spans $\g^*$, 
and the Laplace transform
of the corresponding measure $\Psi_*\lambda_M$ is finite in
one point, then $\g$ is admissible.

\begin{thm} \mlabel{thm:5.9} 
  Let $\g$ be a finite-dimensional Lie algebra
  and $\mu$ a positive $\Ad^*(G)$-invariant Borel 
  measure on $\g^*$ whose support spans~$\g^*$.
  If there exists an $x \in \g$ with $\cL(\mu)(x) < \infty$, then
  \begin{itemize}
  \item[\rm(a)]  $\cL(\mu)$ is reduced in the sense of {\rm Subsection~\ref{subsec:1.1}}.
  \item[\rm(b)]  $\g$ is admissible.   
  \item[\rm(c)] There
    exists a compactly embedded Cartan subalgebra $\ft \subeq \g$
    and an adapted positive system $\Delta^+$ for which $C_{\rm min}$
    is pointed and contained in $C_{\rm max}$
      {\rm(cf.~\eqref{eq:maxcon-into}, \eqref{eq:maxcone2-into})}, 
    \[ D_\mu \subeq W_{\rm max}^\circ, \quad \mbox{ and } \quad
      \supp(\mu)\subeq W_{\rm min}^\star.\] 
  \end{itemize}
\end{thm}

\begin{prf} (a) Lemma~\ref{lem:2.3}(b) implies that $\cL(\mu)$ is reduced.

\nin (b),(c):   {\bf Step 1:} First Corollary~\ref{cor:measlem}
  implies that
  \[ D_\mu = \{ x \in \g \: \cL(\mu)(x) < \infty\}
    \subeq \comp(\g)^\circ,\]
  so that
  $\comp(\g)$ has interior points, and thus
  $\g$ possesses a compactly embedded Cartan subalgebra~$\ft$.
  Then 
\[ \comp(\g)^\circ = \Ad(G).(\ft \cap \comp(\g)^\circ)\]
(\cite[Thm.~VII.1.8(i)]{Ne00}). In particular, we have 
$D_\mu \cap \ft \not=\eset$. 

  \nin {\bf Step 2:} Next we show that $\g$ has  {\it cone potential}, i.e.,
  for $0 \not= x_\alpha \in \g_\C^\alpha$ with $\alpha \in \Delta_p$, we have
  $[x_\alpha, x_\alpha^*] \not=0$. We assume that this is not the case.
  We pick $h \in \ft \cap D_\mu$ and consider the $3$-dimensional
  subspace
  \[ \fb := \R h +    \R (x_\alpha - x_\alpha^*) +
    \R i(x_\alpha + x_\alpha^*) \subeq \g.\]
  As $[h, x_\alpha] = \alpha(h) x_\alpha \in i \R x_\alpha$, it follows
  that $\fb$ is a Lie subalgebra.
  Further,  $h \in \comp(\g)^\circ$ by Step 1, so that
  the Lie algebra $\ker(\ad h)$ is compact,
  hence cannot contain the non-compact Lie algebra $\ft + \fb$.
  Therefore $\alpha(h) \not=0$, and thus $\fb$ is   
 isomorphic to the Lie algebra
 $\mot_2(\R)$ of the motion group of the euclidean plane.
 We  write it as $\fb = \R^2 \rtimes \R h$ with
 $h = \pmat{ 0 & 1 \\ -1 & 0}$ and
 $\fb^*$ is spanned by the dual basis $\be_1^*, \be_2^*$ and $h^*$. 
  Let $\mu_\fb$ denote the projection of the measure $\mu$ under the
  restriction map $\g^* \to \fb^*$. Its support spans $\fb^*$.
 The non-trivial 
  coadjoint orbits in $\fb^*$ are   cylinders
  \[ \cO_r := \{ a \be_1^* +  b \be_2^* + c h^* \:
    a^2 + b^2 = r^2 \}, \quad r > 0, \]
  with the axis $\R h^*$, and $\R h^*$ consists
  of fixed points because $h^*([\fb,\fb]) = \{0\}$.
  On any non-trivial orbit $\cO_r$, $r > 0$, 
 the invariant   measure $\mu_r$ satisfies $\cL(\mu_r)(h) = \infty$ because
  its projection to the axis $\R h^*$ is translation invariant.

  We decompose $\mu_\fb$ as sum $\mu_\fb^0 + \mu_\fb^1$, where
  $\mu_\fb^0$ is supported in $\ft_\fb^*$ and $\mu_\fb^1$ on its complement.
The measure $\mu_\fb^1$ has a canonical desintegration
  \[ \mu_\fb^1 = \int_{(0,\infty)}  \mu_r\, d\nu(r) \]
  for some positive measure $\nu$ on $(0,\infty)$, so that the finiteness of 
  \[  \cL(\mu_\fb^1)(h) =  \int_{(0,\infty)}  \cL(\mu_r)(h)\, d\nu(r)
    = \infty \cdot \nu((0,\infty))\]
  implies that $\nu = 0$. Therefore $\mu_\fb$ is supported on
  $\ft_\fb^*$, contradicting that its support spans $\fb^*$.
  This contradiction now implies that $[x_\alpha, x_\alpha^*] \not=0$.

  \nin {\bf Step 3:} We have just seen that $\g$ has cone potential,
  and this implies that it is root reduced, in the sense that the
  subspace $[\ft,\g]$ contains no non-zero ideal
  (\cite[Prop.~VII.2.25]{Ne00}). We now consider the, by Step 1  non-empty, 
  convex subset
  \[   D_\mu \cap \ft \subeq \comp(\g)^\circ.\]
  As $D_\mu$ is $\Ad(G)$-invariant, this set is invariant under the finite 
  Weyl group  $\cW_\fk$, hence contains a fixed point
  $z_0$, i.e.,
  an element in $\fz(\fk)$ (\cite[Lemma~VII.2.11(i)]{Ne00}).
  So $z_0 \in D_\mu \cap  \fz(\fk)$. 
  Since the centralizer of $z_0$ is compact, no non-compact root
  vanishes on $z_0$, and
  \[ \Delta_p^+ := \{ \alpha \in \Delta_p \: i\alpha(z_0) > 0 \} \]
  is a $\cW_\fk$-invariant positive system of non-compact roots.
  Picking a regular element $x_0 \in \ft$ so close to $z_0$
  that, for $\alpha \in \Delta_k$ and $\beta \in \Delta_p$,
  we have $|\alpha(x_0)| < |\beta(x_0)|$, the subset 
  \[ \Delta^+ := \{ \alpha \in \Delta \: i\alpha(x_0) > 0 \} \]
  is an adapted positive system with $\Delta^+ \cap \Delta = \Delta_p^+$
  (cf.~\cite[Prop.~VII.2.12]{Ne00}). 
Now $z_0 \in C_{\rm max}^\circ$, and since
  $C_{\rm max}^\circ$ is a connected component of $\comp(\g)^\circ \cap \ft$,
  the convexity of $D_\mu$ implies that
  \begin{equation}
    \label{eq:dmuincl}
 D_\mu \cap \ft \subeq C_{\rm max}^\circ,
 \quad \mbox{   hence that } \quad  D_\mu = \Ad(G)(D_\mu \cap \ft)
 \subeq W_{\rm  max}^\circ.
  \end{equation}

  \nin {\bf Step 4:} Let $\g_1 := \Spann D_\mu$. As $\Ad(G)D_\mu =  D_\mu$,
  this is an ideal of $\g$. 
  Proposition~\ref{prop:1.2}(c) shows that $\g_1$ is admissible.
  It contains the element $z_0 \in D_\mu \cap \fz(\fk)$, and
  $\ft_1 := \ft \cap \g_1$ is a compactly embedded abelian subalgebra
  of $\g_1$.  Since no non-compact root $\alpha$
  vanishes on $z_0$, we obtain
  \[ \g_\C^\alpha = [z_0, \g_\C^\alpha] \subeq \g_{1,\C},\]
  and this implies that the unique maximal compactly embedded subalgebra
  $\fk \subeq \g$ containing $\ft$ (\cite[Prop.~VII.2.5]{Ne00}) satisfies
  \[ \g = \g_1 + \fk.\]
  Since $\fk$ is reductive and $\fk_1 := \fk \cap \g_1$ is an ideal of
  $\fk$, we can write $\fk$ as a direct sum $\fk_1 \oplus \fk_2$ and,
  accordingly, $\ft = \ft_1 \oplus \ft_2$ with $\ft_j := \fk_j \cap \ft$.


As $z_0 \in \ft_1$, we have
  \[ \fz_{\g_1}(\ft_1) \subeq \g_1 \cap \fz_\g(z_0)
    = \g_1 \cap \fk = \fk_1,\]
  and since $\ft_1$ is a Cartan subalgebra of $\fk_1$, it follows that
  $\fz_{\g_1}(\ft_1) = \ft_1$.
  Therefore $\ft_1$ is a compactly embedded Cartan subalgebra of $\g_1$.
 
  \nin {\bf Step 5:} We claim that $C_{\rm min}$ is pointed and contained in
  $C_{\rm max}$. Let $\alpha \in \Delta_p^+ \subeq \Delta(\g,\ft)$.
  Then $\g_\C^\alpha = [\ft_1, \g_\C^\alpha] \subeq \g_{1,\C}$ and, for $\alpha_1 := \alpha\res_{\ft_1}$,
  we have
  \[ \g_{1,\C}^{\alpha_1} = \sum_{\alpha\in \Delta, \alpha\res_{\ft_1} = \alpha_1}
    \g_\C^\alpha.\]
  This implies that
  \[    C_{\rm min} = \oline\cone(\{ i [x_\alpha, x_\alpha^*] \: 
    x_\alpha \in \g_\C^\alpha, \alpha \in \Delta_p^+\})
    \subeq C_{\rm min,\g_1} \subeq \ft_1.\]
  Since $\Delta_p^+\res_{\ft_1}$ are the positive non-compact roots of $\g_1$,
  we also have
  \[ C_{\rm max} \cap \ft_1 = C_{\rm max,\g_1}.\]
  Therefore it suffices to show that $C_{\rm min,\g_1}$ is pointed
  and contained in $C_{\rm max,\g_1}$.

  For any $x \in D_\mu \cap \ft \subeq C_{\rm max,\g_1}$
  (cf.\ \eqref{eq:dmuincl}), it follows from
  the Convexity Theorem for Adjoint Orbits (\cite[Thm.~VIII.1.36]{Ne00}) 
  that
  \begin{equation}
    \label{eq:contheo1}
    x + C_{\rm min,\g_1} \subeq \conv(\Ad(G)x) \subeq D_\mu,
  \end{equation}
  so that
  \[ \cL(\mu)(x + y)
    \leq \sup \cL(\mu)(\Ad(G)x)  = \cL(\mu)(x)  \quad \mbox{ for } \quad
    y \in C_{\rm  min,\g_1}.\]
  Since the function $\cL(\mu)$ is reduced by (a), we must have 
$-y \not\in C_{\rm  min,\g_1}$, i.e.,
  that $C_{\rm min,\g_1}$ is pointed.
  Further, \eqref{eq:contheo1} and
  $D_\mu \cap \ft_1 \subeq C_{\rm  max,\g_1}$ entail that
  $C_{\rm min,\g_1} \subeq C_{\rm max,\g_1}$.
  We thus obtain that $C_{\rm min}$ is pointed and contained
  in $C_{\rm max}$. Finally \cite[Thm.~VII.1.19]{Ne00}
  implies that $\g$ is admissible because it contains a compactly
  embedded Cartan subalgebra, is root reduced, and there exists an
  adapted positive system $\Delta^+$ for which
  $C_{\rm min}$ is pointed and contained in $C_{\rm max}$. 

  \nin {\bf Step 6:} We have already seen in
  \eqref{eq:dmuincl} that
  $D_\mu \subeq W_{\rm  max}^\circ.$ 
The convexity of the $\Ad(G)$-invariant function $\cL(\mu)$ on
$D_\mu$, combined with the relation
\[ W_{\rm min} = \lim(\oline{\conv}(\cO_x)) \quad \mbox{ for } \quad
  x \in W_{\rm max}^\circ \]
(Theorem~\ref{thm:invcon-adm}(b)) shows that
\[ \cL(\mu)(x + y) \leq \cL(\mu)(x) \quad \mbox{ for } \quad
  x \in D_\mu, y \in W_{\rm min}, \]
and this in turn leads with
Lemma~\ref{lem:2.3}(a) to $\supp(\mu)\subeq W_{\rm min}^\star$.
\end{prf}

\section{Symplectic Gibbs ensembles} 
\mlabel{sec:3} 

In this section we introduce some of the key concepts concerning Gibbs 
ensembles associated to a Hamiltonian action of a Lie group 
(cf.\ \cite{Ba16}): geometric temperature, the Gibbs ensemble, thermodynamic 
potential and geometric heat. 

\begin{itemize}
\item 
Let $\sigma \: G \times M \to M$ be a (strongly) Hamiltonian action
of the Lie group $G$ on the symplectic manifold 
$(M,\omega)$ and $\Psi \: M \to \g^*$ the corresponding
equivariant momentum map. 
For the derived action 
\[ \dot\sigma \: \g\to \cV(M), \quad
  \dot\sigma(x)(m) := \frac{d}{dt}\Big|_{t = 0}
  \exp(-tx).m,   \] this implies that 
\[ i_{\dot\sigma(x)} \omega = -\dd H_x \quad \mbox{ for }\quad 
H_x(m) := \Psi(m)(x).\] 
\item We write $\lambda_M$ for the 
Liouville measure on $M$, specified by the volume form 
\[ \frac{\omega^n}{(2\pi)^n n!}, \quad \mbox{ where } \quad 
2n = \dim M.\] 
Then the corresponding push-forward measure on~$\g^*$ 
is denoted $\mu := \Psi_*\lambda_M$. 
\end{itemize}

\begin{ex} Throughout this paper, we shall mostly be concerned
  with the case where
  $M = \cO_\lambda = \Ad^*(G)\lambda \subeq \g^*$ is a coadjoint orbit
  in $\g^*$, endowed with the Kostant--Kirillov--Souriau
  symplectic form, given by
\begin{equation}
  \label{eq:omegaalpha}
 \omega_\alpha(\alpha \circ \ad x, \alpha \circ \ad y)
 := \alpha([x,y]) \quad \mbox{ for } \quad x,y \in \g.
\end{equation}
Here
\[ \dot\sigma(x)(\alpha) = \alpha \circ \ad x, \quad
    H_x(\alpha) = \alpha(x),\]
  and the momentum map is the inclusion
  $\Psi \: \cO_\lambda \to \g^*$.
\end{ex}

\begin{defn} \mlabel{def:3.1} 
{\rm(Geometric temperature of a Hamiltonian action)} 
The {\it geometric temperature} is the 
set $\Omega$ of all elements $x\in \g$ for which the 
Hamiltonian functions $H_y, y \in \g,$ have the property that 
\[  \int_M e^{- H_y(m)} \, d\lambda_M(m) < \infty \]
for all $y$ in a neighborhood of $x$. 
This means that the Laplace transform
\[ Z(x) := \cL(\mu)(x) = \int_{\g^*} e^{-\alpha(x)} \, d\mu(\alpha)
  = \int_M e^{- H_x(m)} \, d\lambda_M(m) \]
is finite on a neighborhood of some $x \in \g$. 
It is smooth on the interior $\Omega := \Omega_\mu$ 
of its domain $D_\mu := \cL(\mu)^{-1}(\R)$ in $\g$  
(cf.\ Lemma~\ref{lem:2.3}). Elements
$x \in \Omega$ are called {\it generalized temperatures}.
For $x \in \Omega$, the measure  
\begin{equation}
  \label{eq:gibbslam}
  \lambda_x:= 
  \frac{e^{-H_x}}{Z(x)} \lambda_M
\end{equation}
is a probability measures on $M$, and 
$\mu_x := \Psi_* \lambda_x$ is a probability measure on~$\g^*$. We write 
\begin{equation}
  \label{eq:qdef}
  Q \: \Omega \to \g^*, \quad
  Q(x) := 
  \int_{\g^*} \alpha\, d\mu_x(\alpha)
  = \int_M \Psi(m)\, d\lambda_x(m) \in \oline\conv(\Psi(M)) \subeq \g^* 
\end{equation}
for the expectation value of the probability measure~$\mu_x$
(see \eqref{eq:exp-val}). 
The following terminology comes from \cite{So97} and \cite{Ba16}.
\begin{itemize}
\item 
The family $(\lambda_x)_{x \in \Omega}$ is called 
the {\it Gibbs ensemble of the dynamical group $G$}, acting on $M$, 
\item the map $-\log Z$ is called the {\it thermodynamic potential}, and 
\item $Q \: \Omega \to \g^*$ is called the {\it geometric heat}. 
\end{itemize}
\end{defn}

\begin{rem} \mlabel{rem:entropy}
  In the relation
  \[    s(x)  =    Q(x)(x) + \log Z(x) \]
  from \eqref{eq:entropie-s} in Definition~\ref{def:entropy},
$Q(x)(x)$ is the mean value of the Hamiltonian
function $H_x$ with respect to the probability measure $\lambda_x$, 
hence is interpreted as ``heat'' in the thermodynamical context. 
All other probability measures on $M$, which are completely continuous with 
respect to the Liouville measure $\lambda_M$ and for which $\Psi$ 
has the same expectation value $Q(x)$, 
have an entropy strictly less than~$s(x)$ by Theorem~\ref{thm:b.3}. So 
$\lambda_x$ maximizes the entropy in this class of measures.
This is in accordance with the 2nd Principle of Thermodynamics
which implies that entropy should be maximal in equilibrium states. 
%
\end{rem}

\begin{rem} (a) The measure $\mu = \Psi_*\lambda_M$ 
  on $\g^*$ is $G$-invariant because $\Psi$ is equivariant and
  $\lambda_M$ is $G$-invariant. Therefore $\cL(\mu)$ 
is an invariant convex function on $\Omega$. 

\nin (b) If $\Omega \not=\eset$, then $\mu$ defines a Radon measure on 
$\g^*$, i.e., compact subsets have finite measure. In fact,
the measures $e^{-H_x}\lambda_M$ are finite and the density is bounded
away from $0$ on every compact subset. 

\nin (c) For $M = \R^{2n}$ and $\omega = \sum_{j = 1}^n \dd p_j \wedge \dd q_j$, we have 
$\frac{\omega^n}{n!} = \dd p_1 \wedge \dd q_1 \wedge \cdots \wedge
\dd p_n \wedge \dd q_n,$ 
the Lebesgue volume form in the coordinates 
$(p_1,q_1,\ldots, p_n,q_n)$.
\end{rem} 

\begin{rem} Following Souriau \cite{So97},
  in \cite{Ma20a}, C.-M.~Marle calls $x \in \g$ a
generalized temperature if there exists an integrable function
$f \: M \to \R^+$ and a neighborhood $U$ of $x$ such that
\[ (\forall y \in U)(\forall m \in M)\ e^{-\Psi(m)(y)} \leq f(m).\]
This clearly implies that $\cL(\mu)(y) < \infty$, so that
$x \in \Omega$ in the sense of Definition~\ref{def:3.1}.
If, conversely,
$x \in \Omega$, then there exist affinely independent elements
$x_0, \ldots, x_n \in \Omega$ with
\[ x = \frac{1}{n}(x_0 + \cdots + x_n), \]
and, for all $y \in \conv(\{x_0, \ldots, x_n\})$ and $m \in M$, we have
\[  e^{-\Psi(m)(y)}
  \leq f(m) := \max_{j =  0,\ldots, n} e^{-\Psi(m)(x_j)}
  \leq  \sum_{j = 0}^n e^{-\Psi(m)(x_j)},\]
so that $f$ is integrable. This shows that our simpler
definition of the geometric temperature $\Omega_\mu$ as
the interior of $D_\mu$ is  consistent with
\cite{So97} and \cite{Ma20a}. 
\end{rem}

\section{Coadjoint orbits} 
\mlabel{sec:5}

In this section we specialize the general setting of symplectic 
Gibbs ensembles from Section~\ref{sec:3} to the case where 
$M = \cO_\lambda := \Ad^*(G)\lambda$ is a coadjoint orbit,
endowed with the Kostant--Kirillov--Souriau symplectic form~\eqref{eq:omegaalpha}.

Let $G$ be a connected Lie group with Lie algebra~$\g$. 
For a coadjoint orbit $\cO_\lambda$, we write 
$\mu_\lambda$ for the Liouville measure on $\cO_\lambda$ and 
consider its geometric temperature
\begin{equation}
  \label{eq:omegalambda}
 \Omega_\lambda := \{  x \in \g \: \cL(\mu_\lambda)(x) < \infty\}^\circ
\end{equation}
(cf.\ Definition~\ref{def:3.1}). We write 
\[ C_\lambda := \oline\conv(\cO_\lambda) \]
for the closed convex hull of $\cO_\lambda$
and $\aff(\cO_\lambda)$ for the affine subspace generated by~$\cO_\lambda$.

\subsection{Generalities}
\mlabel{subsec:5.1}

In view of Theorem~\ref{thm:5.9}, the cases of interest 
arise for admissible Lie algebras $\g$. More precisely, we have
the following corollary to Theorem~\ref{thm:5.9}: 

\begin{cor} \mlabel{cor:4.1} 
Suppose that the coadjoint orbit 
$\cO_\lambda$ spans $\g^*$ and that
$D_{\mu_\lambda} \not=\eset$.
Then the following assertions hold:
\begin{itemize}
\item the convex functions 
$\log \cL(\mu_\lambda)$ and $\cL(\mu_\lambda)$ are reduced,
\item the Lie algebra $\g$ is admissible, and
\item there exists a compactly embedded Cartan subalgebra $\ft
\subeq \g$ and an adapted positive system $\Delta^+$
with $C_{\rm min}$ pointed and contained in $C_{\rm max}$, such that 
\[ \lambda \in W_{\rm min}^\star, \quad \mbox{ and } \quad
  \Omega_\lambda \subeq W_{\rm max}^\circ. \]
Here $\Delta_p^+$ is uniquely determined by $\lambda$. 
\end{itemize}
\end{cor}

\begin{rem} \mlabel{rem:5.2} (Reduction to spanning orbits) 
  We may always assume that $\cO_\lambda$ spans $\g^*$.
  Otherwise $\fn := \cO_\lambda^\bot \trile \g$ is an ideal, and we can factorize
the Hamiltonian action of $G$ to one of a group with Lie algebra $\g/\fn$. 
Then we have an inclusion of Lie algebras 
$\g \into (C^\infty(\cO_\lambda), \{\cdot,\cdot\}),
x \mapsto H_x$. In particular, 
an element $z \in \g$ is central if and only if it defines a 
constant function on $\cO_\lambda$, as follows from
\begin{equation}
  \label{eq:hz}
 H_z(\Ad^*(g)\alpha) = \alpha(\Ad(g)^{-1}z) \quad \mbox{ for }  \quad
 g \in G, \alpha \in \cO_\lambda.
\end{equation}
So 
$\dim \fz(\g) \leq 1$, and $\cO_\lambda$ is contained in a proper 
hyperplane in $\g^*$ if and only if $\fz(\g) \not=\{0\}$. In the latter 
case, $\cO_\lambda \subeq \lambda + \fz(\g)^\bot$.  
\end{rem}

\begin{prop} \mlabel{prop:legendre} If $\cO_\lambda$ spans $\g^*$
  and $\Omega_\lambda \not=\eset$, then the geometric heat 
  \[ Q \: \Omega_\lambda
    = \Big\{ x \in \g \:  \int_{\cO_\lambda} e^{- \alpha(x)}\,
    d\mu_\lambda(\alpha) < \infty \Big\}^\circ  \to \g^*, \quad 
Q(x) = \frac{1}{\cL(\mu_\lambda)(x)} \int_{\cO_\lambda} 
\alpha \cdot e^{-\alpha(x)}\, d\mu_\lambda(\alpha)\] 
has the following properties: 
  \begin{itemize}
  \item[\rm(a)] $\Omega_\lambda + \fz(\g) = \Omega_\lambda$ and 
$Q(x+z) = Q(x)$ for $z \in \fz(\g)$ and $x\in \Omega_\lambda$. 
  \item[\rm(b)] $Q$ factors through a function 
$\oline Q \: \Omega_\lambda/\fz(\g) \to  C_\lambda$ 
which is a diffeomorphism onto an open subset of the
affine space~$\aff(\cO_\lambda)$ generated by $\cO_\lambda$. 
  \end{itemize}
\end{prop}

\begin{prf} (a) follows immediately from the fact that the functions 
  $H_z(\alpha) = \alpha(z)$, $z \in \fz(\g)$, are constant on $\cO_\lambda$ 
  (cf.\ \eqref{eq:hz} in Remark~\ref{rem:5.2}). 

\nin (b) The existence of the factorized function $\oline Q$ follows from (a). 
Since $Q(x)$ is the center of mass of a probability measure on 
$\cO_\lambda$, it is contained in $C_\lambda$.
The remaining assertions follow from Proposition~\ref{prop:I.9}(iii). 
\end{prf}

\begin{ex} \mlabel{ex:sl2nil} 
(A non-closed coadjoint orbit with tempered Liouville measure) \\
For $\g = \fsl_2(\R)$, the coadjoint action is equivalent to the action 
of the group $\SO_{1,2}(\R)_e$ on $3$-dimensional Minkowski space because
the Cartan--Killing form has signature $(1,2)$. Then 
\[ \cO 
:= \{ (x_0, x_1, x_2) \: x_0 := (x_1^2 + x_2^2)^{1/2}, (x_1,x_2) \not=(0,0)\} 
= \{ (x_0, x_1, x_2) \: x_0 > 0, x_0^2 = x_1^2 + x_2^2\} \] 
is a nilpotent orbit (an orbit of a nilpotent element)
and the corresponding Liouville measure is proportional to the 
measure defined by 
\[\int_{\R^3} f(x_0, x_1, x_2)\,  d\mu(x_0, x_1, x_2) := 
\int_{\R^2} f\big((x_1^2 + x_2^2)^{1/2}, x_1, x_2\big)\, 
\frac{dx_1 \, dx_2}{\sqrt{x_1^2 + x_2^2}},\]
because both are invariant under rotations and boosts. 
In polar coordinates, it is plain that this  measure is  tempered.
We conclude that there exist non-closed coadjoint orbits 
whose Liouville measure is tempered. 

For the Laplace transform of this measure, we obtain 
\begin{align*}
 \cL(\mu)(z, s \cos \theta, s \sin \theta)
&=  \int_{\R^2} e^{-z (x_1^2 + x_2^2)^{1/2}} e^{-s \la (\cos \theta, \sin \theta), (x_1,x_2) \ra} 
\frac{dx_1 \, dx_2}{\sqrt{x_1^2 + x_2^2}}\\
&=  \int_0^\infty \int_0^{2\pi} 
e^{-z r} e^{-sr (\cos(\theta) \cos(\phi) + \sin(\theta) \sin(\phi))} \, d\phi\, dr \\
&=  \int_0^\infty \int_0^{2\pi} 
e^{-z r} e^{-sr \cos(\theta - \phi)} \, d\phi\, dr
=  \int_0^\infty  \int_0^{2\pi} e^{-r(z + s\cos(\phi))} \, d\phi\, dr \\
&=  \int_0^\infty e^{-rz} \Big(\int_0^{2\pi} e^{-rs\cos(\phi)} \, d\phi\Big)\, dr.
\end{align*}
The next to last expression for this integral shows that, 
for $0 \leq s < z$, this integral exists. For 
$s = 0 < z$, we obtain in particular 
\begin{align*}
 \cL(\mu)(z, 0,0)
= 2\pi \int_0^\infty e^{-rz} \, dr
= \frac{2\pi}{z}.   
\end{align*}
By the invariance of $\cL(\mu)$, this leads to 
\[  \cL(\mu)(z, s \cos \theta, s \sin \theta)
= \frac{2\pi}{\sqrt{z^2 - s^2}}\quad \mbox{ for } \quad z > s \geq 0.\]
Note that 
\[ \cL(\mu)(rx) = r^{-1} \cL(\mu)(x) \quad \mbox{ for }\quad r > 0, x \in\g. \] 

This example is also discussed explicitly in \cite[\S 4.3]{BDNP23}, thus
correcting invalid claims in \cite[\S 3.3]{Ma21} and
\cite[\S 3.5]{Ma20b},
asserting that this orbit does not have a non-trivial geometric
temperature. 
\end{ex}

\begin{rem} 
The finiteness of $\cL(\mu_\lambda)$ in some point 
$x \in \g$ implies that $\mu_\lambda$ is a Radon measure, 
i.e., finite on compact subsets of $\g^*$. 
By  \cite[Thm.~1.8]{Ch90}, 
the Liouville measure of any closed coadjoint orbit of a 
connected Lie group  is tempered, but 
Example~\ref{ex:sl2nil}  shows that 
the temperedness of $\mu_\lambda$
does not imply that $\cO_\lambda$ is closed.
We shall see in Theorem~\ref{thm:temp} below that $\mu_\lambda$
is always tempered if $D_\mu\not=\eset$. 
\end{rem}

\subsection{The affine action on a symplectic vector space} 
\mlabel{subsec:5.2}

Let $(V,\Omega)$ be a symplectic vector space.
In this subsection we discuss the affine
action of the group $G = \Heis(V,\Omega) \rtimes \Sp(V,\Omega)$ 
on~$V$.  We consider the Lie algebra $\g$ of all functions 
\begin{equation}
  \label{eq:hcwx}
 H_{c,w,x} \: V \to \R, \quad 
 H_{c,w,x}(v) =  c + \Omega(w,v) + \frac{1}{2}\Omega(xv,v), 
\end{equation}
endowed with the Poisson bracket on $(V,\Omega)$. Let $2n = \dim V$. Then 
\[ \g \cong \heis(V,\Omega) \rtimes \sp(V,\Omega),\]
where $\heis(V,\Omega)$  is the $(2n+1)$-dimensional Heisenberg algebra, which
corresponds to the functions $H_{c,w,0}$, $w \in V$, $c \in \R$. 

The linear functional $\lambda = \ev_0 \in \g^*$ given by point 
evaluation in~$0$ takes the form
\[ \lambda(c,w,x) = c = H_{c,w,x}(0).\] 
The action of the Lie algebra $\g$ on $V$ 
integrates to a Hamiltonian action of the corresponding group~$G$,
and the momentum map is given by 
\begin{equation}
  \label{eq:PhionV}
  \Psi \: V \to \g^*, \quad 
  \Psi(v)(f) =  f(v)
\end{equation}
(\cite[Prop.~A.IV.15]{Ne00}). 
It follows in particular that $\Psi(V) = {\cal O}_\lambda\subeq \g^*$
is a coadjoint orbit.

\begin{lem} \mlabel{lem:5.7a} 
For $A\in \Sym_n(\R)$ and $\xi \in \R^n$, we have 
\begin{equation}
  \label{eq:posdef}
\frac{1}{\sqrt{2\pi}^n} \int_{\R^n} e^{-\frac{1}{2}\la Ax,x\ra - \la 
\xi, x \ra} \ dx 
=
\begin{cases}
  \det(A)^{-1/2} \cdot e^{\frac{1}{2} \la A^{-1}\xi,\xi\ra}
  & \text{ for\ } A \text{\ positive definite} \\ 
\infty  & \text{ otherwise}.     
\end{cases}
\end{equation}
\end{lem}

\begin{prf} We may evaluate the integral in coordinates adapted to an
  orthogonal  basis of eigenvectors of~$A$, where it boils down to the
  $1$-dimensional case.   
\end{prf}

We now put this into a symplectic context.
We call a complex structure
$I \in \Sp(V,\Omega)$ {\it positive} if 
\begin{equation}
  \label{eq:scalsymp}
  \la v,w \ra := \Omega(v,Iw)
\end{equation}
is  positive definite. Any positive complex structure determines
a maximal compactly embedded subalgebra $\fk_I \subeq \g$ by
\[ \fk_I = \fz_\g(I) = \R \times \{0\} \times \fk_{I,\fs}, \quad
  \fk_{I,\fs} := \{ x \in \sp(V,\omega) \: [x,I] = 0\}.\]
Now any  $x \in \sp(V,\Omega)$, for which $H_x(v) = H_{0,0,x}(v)$
is positive definite, is a compact element, hence contained in a conjugate
of some $\fk_I$, which means that there exists a complex structure
$I$ with $x \in \fk_I$. 

\begin{lem}  \mlabel{lem:5.7} 
  Let $(V, \Omega)$ be a $2n$-dimensional symplectic vector space
  and $(c,w,x) \in \hsp(V,\Omega)$.
  If $H_{0,0,x}$ is positive definite and $I \in \Sp(V,\Omega)$
  with $x \in \fk_I$, then there exists a constant $c_V$ such that 
\begin{align}\label{eq:gauss-int1}
\int_V e^{-H_{c,w,x}(v)}\ d\lambda_V(v) 
  =   \begin{cases}
        c_V \exp(- H_{c,w,x}(-x^{-1}w)) \det(Ix)^{-\frac{1}{2}}
        & \text{ for\ } H_{0,0,x} \text{\ positive definite}, \\ 
\infty  & \text{ otherwise}.
\end{cases}
\end{align}
\end{lem}

Note that $-x^{-1}w$ is the unique minimum of $H_{c,w,x}$ on $V$.

\begin{prf} This can be derived from Lemma~\ref{lem:5.7a}.
  If $H_{0,0,x}$ is not positive definite, it follows that the
  integral does not exist. So we may assume that this quadratic
  function is positive definite.
  Then there exists a positive complex structure $I\in \Sp(V,\Omega)$ 
 commuting with~$x$  (cf.~\eqref{eq:scalsymp}). 
Then the Liouville measure $\lambda_V$ is given by the volume form 
$ \frac{\Omega^n}{(2\pi)^n n!}$ which is a multiple of
Lebesgue measure with respect to the scalar product. We thus obtain
for a suitable constant $c_V> 0$: 
\begin{align}\label{eq:gauss-int1}
\int_V e^{-H_{c,w,x}(v)}\ d\lambda_V(v) 
&=  e^{-c} \int_V e^{-\Omega(w,v) - {\frac{1}{2}}\Omega(xv,v)}\ d\lambda_V(v) 
=  e^{-c} \int_V e^{-\la Iw,v \ra - {\frac{1}{2}}\la Ixv,v \ra}\ d\lambda_V(v) \notag\\
&=  c_V e^{-c} \det(Ix)^{-{\frac{1}{2}}} e^{{\frac{1}{2}} \la (Ix)^{-1}Iw,Iw\ra} 
=  c_V e^{-c} \det(Ix)^{-{\frac{1}{2}}} e^{-{\frac{1}{2}} \la x^{-1}w, Iw\ra} \notag \\
&=  c_V e^{-c} \det(Ix)^{-{\frac{1}{2}}} e^{{\frac{1}{2}} \Omega(x^{-1}w, w)}. 
\qedhere\end{align}
\end{prf}

\subsection{Admissible coadjoint orbits} 
\mlabel{sec:5.3}

A particular nice class of coadjoint orbits
$\cO_\lambda \subeq \g^*$ are the so-called admissible ones;
they are closed and their convex hull contains 
no affine lines.
In this section we describe the explicit formulas 
for the Laplace transform $\cL(\mu_\lambda)$, $\lambda$ admissible, 
that have been obtained in \cite{Ne96a} with
stationary phase methods for proper momentum maps.

\begin{defn} \mlabel{def:adm-orbit}
  We call a coadjoint orbit $\cO_\lambda$ and the element 
$\lambda \in \g^*$ {\it admissible}, 
if $\cO_\lambda$ is closed and 
its closed convex hull $\oline\conv(\cO_\lambda)$ contains no affine lines 
(\cite[Def.~VII.3.14]{Ne00}). 
\end{defn}

\begin{ex} We consider
the linear functional $\lambda(z,v,x) = z$ on $\hsp(V,\Omega)$,
which corresponds to evaluation in $0$. 
  Let $x \in \sp(V,\Omega) \subeq \hsp(V,\Omega)$ be such that $v \mapsto 
\Omega(xv,v)$ is positive definite. Then 
\[ H_{(0,0,x)} \: V \to \R, \quad H_{0,0,x}(v) 
= {\frac{1}{2}}\Omega(xv,v) \] 
is proper and bounded from below on~$(V,\Omega)$. 
Hence ${\cal O}_\lambda$ is closed in $\hsp(V,\Omega)^*$.
Its convex hull contains no affine lines
because the cone $B(\cO_\lambda)$, which contains all functions
$H_{c,w,x}$ with $H_{0,0,x}$ positive definite, has interior points
(\cite[Prop.~V.1.15]{Ne00}). Therefore ${\cal O}_\lambda$ is admissible. 
\end{ex}

\begin{prop}
  \mlabel{prop:adm-orbits}
  Let $\ft \subeq \g$ a compactly embedded Cartan subalgebra, and $\lambda \in \g^*$. Then the following assertions hold:
  \begin{itemize}
  \item[\rm(a)] If $\cO_\lambda$ is admissible and spans $\g^*$, then 
$\g$ is admissible, $B(\cO_\lambda)^\circ \subeq \comp(\g)$
    {\rm(cf.~\eqref{eq:dualcone})}, and \break $\cO_\lambda \cap \ft^* \not=\eset$,
    where $\ft^* \cong  [\ft,\g]^\bot$. 
Moreover, $B(\cO_\lambda) \subeq
    W_{\rm max}$ for an adapted positive system $\Delta^+
    \subeq \Delta(\g,\ft)$ with
  $C_{\rm min}$ pointed and contained in $C_{\rm max}$.
\item[\rm(b)]  If $\Delta^+$ is adapted with $C_{\rm min} \subeq C_{\rm max}$,
  then $\lambda \in C_{\rm min}^\star \subeq \ft^*$ implies that
  $\cO_\lambda$ is admissible and that $W_{\rm max}^\circ = B(\cO_\lambda)^\circ$.
  \end{itemize}
  
\end{prop}

\begin{prf} (a) That $\g$ is admissible follows from
  the fact that $\g^*$ is spanned by an admissible orbit 
(\cite[Lemma~VII.3.17]{Ne00}), and the ellipticity of
  the cone $B(\cO_\lambda)$ from \cite[Prop.~VIII.1.17(iii)]{Ne00}.
  That $\cO_\lambda$ intersects $\ft^*$ for every compactly embedded
  Cartan subalgebra~$\ft$, follows from \cite[Prop.~VIII.1.4]{Ne00}.

The second assertion now follows from \cite[Thm.~VIII.3.10]{Ne00}.

  \nin (b) follows from (a), and from \cite[Thm.~VIII.1.19]{Ne00}, which
  asserts that $C_{\rm max} \subeq B(\cO_\lambda)$,
  and this in turn entails that $W_{\rm  max}^\circ
  = \Ad(G) C_{\rm max}^\circ \subeq B(\cO_\lambda)$. 
\end{prf}

\cite{Ne96a} contains information 
on Laplace transforms of Liouville measures $\mu_\lambda$
of admissible coadjoint orbits $\cO_\lambda$. 
To explain the formula derived in 
\cite[Thm.~II.10]{Ne96a} for the 
Laplace transform of $\mu_\lambda$, we identify 
the tangent space 
\[ T_\lambda(\cO_\lambda) \cong \lambda \circ \ad \g \cong \g/\g_\lambda, 
\quad \mbox{ where } \quad 
\g_\lambda = \{ y \in \g \: \lambda \circ \ad y = 0\}\] 
is the stabilizer Lie algebra of $\lambda$.
We then write
\[ \Delta_\lambda := \{ \alpha \in \Delta^+ \:
  \g_\C^\alpha \not\subeq (\g_\lambda)_\C \} \]
for those positive roots of the pair
$(\g,\ft)$ that appear in the $\ft$-representation
on the complexified tangent space $T_\lambda(\cO_\lambda)_\C$.
For $\alpha \in \Delta^+$, we write 
\[ m_\alpha^\lambda := \dim_\C \g_\C^\alpha/(\g_{\lambda,\C} \cap \g_\C^\alpha) \]
for the multiplicity of $\alpha$ in this representation, so that
$m_\alpha^\lambda > 0$ if and only if $\alpha \in \Delta_\lambda$.

\begin{defn}
We say that an element $x \in \ft$ 
is {\it $\cO_\lambda$-regular} if, 
for every $w \in \cW_{\fk}$ and $\alpha \in \Delta_\lambda$, we 
have $\alpha(wx) \not=0$. Identifying $\ft^*$ with the subspace
of $[\ft,\g]^\bot \subeq \g^*$, this means that the set 
$\cW_{\fk}\lambda = \cO_\lambda^T$ of 
$T$-fixed points in $\cO_\lambda$
(\cite[Lemma~VIII.1.1]{Ne00}) 
consists of isolated fixed points 
of the one-parameter group $\exp(\R x)$. 
\end{defn}

\begin{thm} \mlabel{thm:II.10} 
  Let $\g$ be admissible, $\ft \subeq \g$ be a compactly embedded Cartan
  subalgebra, $\Delta^+$ an adapted positive system with
  $C_{\rm min} \subeq C_{\rm max}$ and $\lambda \in C_{\rm min}^\star$.
  Then $\lambda$ is admissible and 
  \begin{equation}
    \label{eq:lmulambda}
  \cL(\mu_\lambda)(x) 
  = \sum_{w \in {\cal W}} \frac{e^{-\lambda(wx)}}
  {\prod_{\alpha \in \Delta_\lambda} (i \alpha(wx))^{m_\alpha^\lambda}} 
  = \sum_{w \in {\cal W}} \frac{e^{-\lambda(wx)}}
  {\prod_{\alpha \in \Delta^+} (i \alpha(wx))^{m_\alpha^\lambda}}
  \end{equation}
for every  $x \in C_{\rm max}^\circ$ which is  $\cO_\lambda$-regular.
In particular,
\begin{equation}
  \label{eq:wmaxinomegalambda}
  W_{\rm max}^\circ \subeq \Omega_\lambda.
\end{equation}
\end{thm}

\begin{prf} The admissibility of $\lambda$ follows from
  Proposition~\ref{prop:adm-orbits}(b) and the formula for the
  Laplace transform from \cite[Thm.~II.10]{Ne96a}.

  To verify \eqref{eq:wmaxinomegalambda},
  we note that any $x \in C_{\rm max}^\circ$ on which no root vanishes
  is $\cO_\lambda$-regular, so that
  \eqref{eq:lmulambda} implies that
  $x \in \Omega_\lambda$. Since $\Omega_\lambda$ is convex,
  and $\cO_\lambda$-singular elements are convex combinations of
  $\cO_\lambda$-regular elements,   
  it follows that $C_{\rm max}^\circ \subeq \Omega_\lambda$.
  This entails that
  obtain $W_{\rm max}^\circ = \Ad(G)C_{\rm max}^\circ \subeq \Omega_\lambda$.
\end{prf}

If $\lambda$ is contained in the interior of $C_{\rm min}^\star$, 
then $\Delta_p^+ \subeq \Delta_\lambda$. In fact,
for $0 \not= x_\alpha \in \g_\C^\alpha$ and $\alpha \in \Delta_p^+$,
we have $[x_\alpha, x_\alpha^*] \not=0$ because
$\g$ is admissible (\cite[Thm.~VII.3.10(iv)]{Ne00}).
Since $i[x_\alpha, x_\alpha^*] \in C_{\rm min}$, it follows that
$\lambda(i[x_\alpha, x_\alpha^*]) > 0$, and this implies that
$\g_\C^\alpha \not\subeq (\g_\lambda)_\C$, i.e., $\alpha \in \Delta_\lambda$.
Note that the subset $\Delta_p^+ \subeq \Delta_\lambda$ 
is $\cW_{\fk}$-invariant, so that we obtain the following factorization 
of the right hand side of~\eqref{eq:lmulambda}.

\begin{cor} \mlabel{cor:5.4}
  Let $\cO_\lambda$ be an admissible coadjoint orbit spanning $\g^*$
  with $\lambda \in (C_{\rm min}^\star)^\circ$.
  For $K := \exp \fk$, we write $\mu_\lambda^K$ for the 
Liouville measure of the coadjoint $K$-orbit $\cO_\lambda^K = \Ad^*(K)\lambda 
\subeq \fk^*$. Then 
\begin{equation}
  \label{eq:5.3}
  \cL(\mu_\lambda)(x) 
= \frac{\cL(\mu_\lambda^K)(x) }
{\prod_{\alpha \in \Delta_p^+} (i\alpha(x))^{\dim \g_\C^\alpha}} 
\quad \mbox{ for } \quad x \in C_{\rm max}^\circ. 
\end{equation}
For $N := \sum_{\alpha \in \Delta_p^+} \dim \g_\C^\alpha$, we have 
\begin{equation}
  \label{eq:temp-est1}
  \lim_{t \to 0+} \cL(\mu_\lambda)(tx) t^N  =
  \frac{\vol(\cO_\lambda^K)}{\prod_{\alpha \in \Delta_p^+} (i\alpha(x))^{\dim \g_\C^\alpha}}   < \infty,
\end{equation}
and $\mu_\lambda$ is tempered. 
\end{cor}

\begin{prf} First we apply Theorem~\ref{thm:II.10}
  to obtain
\[   \cL(\mu_\lambda)(x) 
= \frac{1}{\prod_{\alpha \in \Delta_p^+} (i\alpha(x))^{\dim \g_\C^\alpha}} 
\cdot \Big(\sum_{w \in {\cal W}} \frac{e^{-\lambda(wx)}}
{\prod_{\alpha \in \Delta_{k,\lambda}} i \alpha(wx)}\Big)  
\ {\buildrel (*) \over =}\ \frac{\cL(\mu_\lambda^K)(x) }
{\prod_{\alpha \in \Delta_p^+} (i\alpha(x))^{\dim \g_\C^\alpha}} 
\]
for those $x \in C_{\rm max}^\circ$ which are $\cO_\lambda$-regular.
Here $(*)$ follows  by applying Theorem~\ref{thm:II.10} to the
compact Lie algebra~$\fk$. 
Since $\cL(\mu_\lambda)$ is a continuous function on 
$C_{\rm max}^\circ$
(\cite[Prop.~V.3.2]{Ne00})
and $\cL(\mu_\lambda^K)$ is continuous on all of~$\ft$,
we obtain \eqref{eq:5.3} by continuity of both sides on~$C_{\rm max}^\circ$.

The assertion on temperedness now follows from
Proposition~\ref{prop:dom-lapl-temp}(b),
where the estimate \eqref{eq:temp-est1} follows from 
$\lim_{t \to 0+}  \cL(\mu_\lambda^K)(tx) = \vol(\cO_\lambda^K)$.
\end{prf}

\begin{ex} The following $2$-dimensional examples also appear in
 \cite[\S 3.3]{Ma21} and \cite{Neu22}. 

\nin  (a) For $\g = \fsl_2(\R)$, we have $\ft = \fk$ and
  we may fix a basis element
  \begin{equation}
    \label{eq:z0-sl2}
    z_0 := \frac{1}{2} \pmat{0 & 1 \\ - 1& 0} \in \ft.
  \end{equation}
  Then we can chose the positive system in such a way that
  $\Delta^+ = \Delta_p^+ = \{ \alpha\}$ with $i\alpha(z_0) = 1$.
  Then $\lambda \in (C_{\rm min}^\star)^\circ \subeq \ft^*$
  if and only if $\lambda(z_0) > 0$, and in this case $\cO_\lambda$
  is a K\"ahler manifold isomorphic to the complex unit disc/upper
  half plane, whose form is scaled by $\lambda(z_0)$. 
  As $\fk$ is abelian, its coadjoint orbits are trivial, so that
  Corollary~\ref{cor:5.4} yields 
  \[ \cL(\mu_\lambda)(t z_0)
    = \frac{e^{-\lambda(tz_0)}}{i\alpha(tz_0)} 
    = \frac{e^{-t\lambda(z_0)}}{t}.\]

  Note that, for $\lambda(z_0) \to 0$, we obtain
  \begin{equation}
    \label{eq:limit-to-0}
\lim_{\lambda(z_0) \to 0}  \cL(\mu_\lambda)(t z_0) = \frac{1}{t},      
\end{equation}
which is a multiple of the
Laplace transform of the nilpotent orbit to which, on the level
of subsets of~$\g^*$, the orbits $\cO_\lambda$ ``converge''
(Example~\ref{ex:sl2nil}). 

  \nin (b)  For $\g = \su_2(\C)$, we may also take
  $\ft = \R z_0$ with $z_0$ as in \eqref{eq:z0-sl2}.
  We chose the positive system in such a way that
  $\Delta^+ = \Delta_k^+ = \{ \alpha\}$ with $-i\alpha(z_0) = 1$
  (cf.~the definition of compact roots in Subsection~\ref{subsec:1.2}) 
  and note that $\cW_{\fk} = \{ \pm \id_\ft\}$. 
  Then $C_{\rm min} = \{0\}$, $C_{\rm min}^\star = \ft^*$,  and
  $m_\alpha^\lambda = 1$ for $\lambda \not=0$.
  Here $\cO_\lambda$ is a compact K\"ahler manifold isomorphic to $\bS^2$,
  whose symplectic form is scaled by $\lambda(z_0)$, which we assume 
w.l.o.g.~to be $\geq 0$. 

  We thus obtain
  \[ \cL(\mu_\lambda)(t z_0)
    = \frac{e^{-\lambda(tz_0)}}{i\alpha(tz_0)}
    +     \frac{e^{\lambda(tz_0)}}{i\alpha(-tz_0)}
= \frac{e^{-\lambda(tz_0)} -e^{\lambda(tz_0)}}{-t} 
=  2\frac{\sinh(t\lambda(z_0))}{t}.\]

\nin (c) The third $2$-dimensional example, where
$\cO_\lambda \cong \R^2 \cong \C$, with a flat K\"ahler structure,
arises for $\g = \heis(\R^2,\Omega) \rtimes \R z_0$, $z_0$ as above,
and $\ft = \R \bc \oplus \R z_0$.

As $[\g,\g] = \heis(\R^2,\Omega)$ is a hyperplane in $\g$, there
exist non-zero linear functionals $\zeta$ vanishing on $[\g,\g]$,
and these are fixed points of the coadjoint action. We thus have
\begin{equation}
  \label{eq:translat}
  \cO_{\lambda + \zeta} = \zeta + \cO_\lambda, 
\end{equation}
where translation by $\zeta$ is a $G$-equivariant
symplectic isomorphism from $\cO_\lambda$ to  $\cO_{\lambda + \zeta}$. 
Then we can chose the positive system in such a way that
  $\Delta^+ = \Delta_r^+ = \{ \alpha\}$ with $i\alpha(z_0) = \frac{1}{2}$.
Here $\lambda \in (C_{\rm min}^\star)^\circ$ 
if and only if $\lambda(\bc) > 0$. 
Then $\cO_\lambda$ 
is a K\"ahler manifold isomorphic to the complex plane,
whose form is scaled by $\lambda(\bc)$.
Combining (a) with \eqref{eq:translat},  we obtain
  \[ \cL(\mu_\lambda)(s \bc + t z_0)
    = e^{-s\lambda(\bc)} \frac{e^{-\lambda(tz_0)}}{i\alpha(tz_0)} 
    = \frac{2 e^{-s\lambda(\bc)-t \lambda(tz_0)}}{t}.\]
This follows from the discussion in Subsection~\ref{subsec:5.2}. 
\end{ex}

\begin{rem} (a) 
For $\cW_{\fk}$-invariant functionals $\lambda_0 \in \ft^*$ with 
$\lambda$ and $\lambda + \lambda_0 \in (C_{\rm min}^\star)^\circ$,
we obtain in particular from Corollary~\ref{cor:5.4} that 
\[ \cL(\mu_{\lambda + \lambda_0})(x) 
= e^{-\lambda_0(x)} \cL(\mu_{\lambda})(x).\] 

\nin (b)  If $\cO_\lambda$ is admissible and spans $\g^*$,
  we find for   $x \in \partial C_{\rm max}$ and $y \in C_{\rm max}^\circ$
  that 
  \[ \lim_{t \to 0+} \cL(\mu_\lambda)(x + ty) = \infty. \]
  By the continuity of $\cL(\mu)$ on closed rays
  (\cite[Cor.~V.3.3]{Ne00}), this implies that $x \not\in D_{\mu_\lambda}$.
  This also follows from Corollary~\ref{cor:measlem}. 
\end{rem}

\subsection{From reductive to simple Lie algebras}
\mlabel{subsec:5.4}

  Suppose that $\g = \g_0 \oplus \g_1 \oplus \cdots \oplus \g_n$ 
  is a direct sum of Lie algebras, where $\g_0$ is abelian,
  then the coadjoint orbit of
$\lambda = \sum_{j = 0}^n \lambda_j$  with $\lambda_j \in \g_j^*$ is a product
of Hamiltonian $G$-spaces 
\[ \cO_\lambda = \{\lambda_0\} \times \prod_{j = 1}^n \cO_{\lambda_j}.\] 
As the Liouville measure is adapted to this product decomposition
(cf.\ \cite[Thm.~16.98]{So97}), 
\[  \cL(\mu)  = e^{- \lambda_0} \cdot \prod_{j = 1}^n \cL(\mu_j),
  \quad 
  D_{\mu_\lambda} = \g_0 \times \prod_{j = 1}^n D_{\mu_{\lambda_j}},\quad  \mbox{ and }
  \quad 
  \Omega_\lambda = \g_0 \times \prod_{j = 1}^n \Omega_{\lambda_j}.\] 
This observation reduces all questions from the
reductive case to simple Lie algebras. 
If $\g_j$ is compact, then $\cO_{\lambda_j}$ is compact and
$\Omega_{\lambda_j} = \g_j$.

\section{Reduction procedures}
\mlabel{sec:6}

In this section we address the classification problem
for coadjoint orbits $\cO_\lambda$ with non-trivial~$D_{\mu_\lambda}$, 
in general finite-dimensional Lie algebras. What we have seen
so far are admissible orbits (Subsection~\ref{sec:5.3}),
which are rather accessible because we have
an explicit formula for the Laplace transform
$\cL(\mu_\lambda)$. The affine coadjoint orbit
$\cO_\lambda \cong (V,\Omega)$ for the non-reductive
Lie algebra $\g = \hsp(V,\Omega)$
is a very special case (Subsection~\ref{subsec:5.2}).

Our strategy for the classification will be to use a
semidirect decomposition $\g = \fu \rtimes \fl$
as in Subsection~\ref{subsec:1.4} to write any orbit
in $W_{\rm min}^\star$ as a sum
\[ \cO_\lambda= \cO_{\lambda_\fz} + \cO_{\lambda_\fl}, \]
which is actually a symplectic product,\begin{footnote}
  {The symplectic
    product $(M,\omega) = (M_1, \omega_1) \times (M_2, \omega_2)$ is
    defined by the relation $\omega = p_{M_1}^*\omega_1 + p_{M_2}^*\omega_2$.
}\end{footnote}
where
$\cO_{\lambda_\fz}$ is isomorphic to the symplectic vector space
$(V,\Omega)$, where $V = [\fl,\fu]$, and $\cO_{\lambda_\fl}$ is a coadjoint orbit
of the reductive Lie algebra $\fl$.

\subsection{Orbits in $W_{\rm min}^\star$}
\mlabel{subsec:6.1}

This subsection is dedicated to the question when
a linear functional $\lambda = \lambda_\fz + \lambda_\fl \in \g^*$
on a semidirect sum
$\g = \fu \rtimes \fl$, which is admissible, 
is contained in $W_{\rm min}^\star$
(cf.\ Theorem~\ref{thm:spind}).

\begin{lem} \mlabel{lem:wminl} For the projection $p_\fl \: \g = \fu \rtimes \fl \to \fl$, we   have
$p_\fl(W_{\rm min}) = W_{\rm min,\fl} \subeq W_{\rm min}.$ 
\end{lem}

Here we use that $\fl$ is admissible as well, so that $W_{\rm min,\fl}$
is defined by $\Delta_p^+ \cap \Delta_s$ as in Theorem~\ref{thm:invcon-adm}. 

\begin{prf} Let $x \in C_{\rm max}^\circ$. Then
  \[ W_{\rm min} = \lim(\co(x))
    \quad \mbox{ for } \quad
    \co(x) := \oline{\conv(\Ad(G)x)}\]
  by Theorem~\ref{thm:invcon-adm}(b).   We have for $\ft_\fl := \ft \cap \fl$
  the decomposition 
  \[ \ft = \fz(\g) \oplus \ft_\fl \quad \mbox{ with }  \quad
    C_{\rm max} = \fz(\g) \oplus (C_{\rm max} \cap \ft_\fl). \]
  We thus assume below that $x = x_\fl \in \ft_\fl$. 

  We write $G = U \rtimes L$ for the simply connected Lie group
  with Lie algebra $\g$. The projection $p_\fl \: \g\to \fl$ is a homomorphism
  of Lie algebras, hence equivariant for the adjoint action,
  and $U$ acts trivially on $\fl$. We conclude that
  \[ p_\fl(\Ad(G)x) = \Ad(L) x = \cO_x^L.\]
  This implies that
$p_\fl(\co(x)) \subeq \co_\fl(x),$   so that
  \[ p_\fl(W_{\rm min})
   = p_\fl(\lim(\co(x)))  \subeq \lim \co_\fl(x) = W_{\rm min,\fl},\]
 where the last equality follows from
  $x \in C_{\rm max}^\circ \cap \ft_\fl
  \subeq C_{\rm max,\fl}^\circ$, where we
  use that $C_{\rm max,\fl} = \ft_\fl \cap (i\Delta_{p,s}^+)^\star$
  (cf.\ Subsection~\ref{subsec:1.2} and \eqref{eq:maxcone2}). 
  As $\co_\fl(x) \subeq \co(x)$ holds trivially, we also  have
  \[ W_{\rm min,\fl} = \lim(\co_\fl(x)) \subeq  \lim(\co(x)) = W_{\rm min}\]
    (cf.\ Theorem~\ref{thm:invcon-adm}(c)). This proves the asserted equality.   
\end{prf}

Below we shall use the notation
\[ C_{\rm mim,\fz} := C_{\rm min} \cap \fz
= \oline\cone(\{ i [x_\alpha, x_\alpha^*] \: 
x_\alpha \in \g_\C^\alpha, \alpha \in \Delta_r^+\}) \subeq \fz,\]
and $C_{\rm min/max,\fl}$ and $W_{\rm min/max,\fl}$ are the cones
specified by $\Delta_p^+ \cap \Delta_s$ in the
admissible Lie algebra~$\fl$ (cf.\ \eqref{eq:maxcon-into})
and Theorem~\ref{thm:invcon-adm}.

\begin{lem} \mlabel{lem:cmin-div} 
  Suppose that $\g$ is admissible and non-reductive.  
  For $\lambda = \lambda_\fz + \lambda_\fl$
  with $\lambda_\fl \in \fl^* \cong \ft^\bot$, and 
  $0 \not=\lambda_\fz \in \fz^* = (V + \fl)^\bot$,
  the following are equivalent:
  \begin{itemize}
  \item[\rm(a)] $\lambda \in W_{\rm min}^\star$. 
  \item[\rm(b)] $\lambda_\fl \in W_{\rm min,\fl}^\star$
    and $\lambda_\fz \in C_{\rm min,\fz}^\star$. 
  \end{itemize}
\end{lem}

If $\cO_\lambda$ is generating, then $\dim \fz(\g) \leq 1$, so that
$C_{\rm min,\fz} = \R_+ \bc$, and the second condition in (b)  reduces
to $\lambda(\bc) \geq 0$.

\begin{prf} (a) $\Rarrow$ (b):
As $\lambda = \lambda_\fz + \lambda_\fl$ and 
$C_{\rm min,\fz} + W_{\rm min,\fl} \subeq W_{\rm min}$
(Lemma~\ref{lem:wminl}) with $C_{\rm min,\fz} \subeq \fz \subeq \fu$ and
$W_{\rm min,\fl} \subeq \fl$, we immediately obtain (b) from (a).

\nin (b) $\Rarrow$ (a): Lemma~\ref{lem:wminl} 
shows that the dual cones satisfy
\[ W_{\rm min,\fl}^\star
  = p_\fl(W_{\rm min})^\star
  = (\fu + W_{\rm min})^\star 
  = W_{\rm min}^\star \cap \fu^\bot
  = W_{\rm min}^\star \cap \fl^*.\] 
We further have, by definition, 
$p_\ft(W_{\rm min}) = C_{\rm min} = C_{\rm min,\fz} + C_{\rm min,\fl},$ 
and this implies that
\[ C_{\rm min,\fz}^\star = \fz^* \cap C_{\rm min}^\star
  \subeq W_{\rm min}^\star. \] 
Therefore $\lambda = \lambda_\fz + \lambda_\fl \in W_{\rm min}^\star$.
\end{prf}

\subsection{Nilpotent orbits in reductive Lie algebras}
\mlabel{subsec:6.2}

In this subsection we
use Rao's Theorem \cite[Thm.~1]{Rao72} on 
adjoint orbits of nilpotent elements  in reductive Lie algebras.
It implies in particular that the Liouville measure
$\mu_\lambda$ is tempered if $\lambda \in \g^*$ is nilpotent.\\

\begin{thm} \mlabel{thm:rao}  {\rm(Rao's Theorem on Nilpotent Orbits)}
  Let $x$ be a nilpotent element in the reductive Lie algebra $\g$
  and let $(h,x,y)$ be a corresponding $\fsl_2$-triple, i.e.,
\[ [h,x] =  2x, \quad [h,y] = -2y \quad \mbox{ and } \quad
  [x,y] = h.
  \begin{footnote}
    {The existence of such elements follows from the Jacobson--Morozov
    Theorem (\cite[Prop.~13.5.3]{Wa72}).}    
  \end{footnote}
\]
We write $\g_\mu := \g_\mu(h)$ for the $h$-eigenspaces,
\begin{footnote}{Note that Rao's paper contains a misprint
    in the definition of the nilpotent Lie algebras
    $\fn$ and $\fn_2$.}\end{footnote}
\[ \fm := \g_0(h), \quad
  \fn := \sum_{\mu > 0} \g_\mu(h), \quad 
  \fp := \sum_{\mu \geq 0} \g_\mu(h) = \fn \rtimes \fm, \quad 
  \fn_2 := \sum_{\mu > 2} \g_\mu(h).\]
Then $V:= \Ad(M)x$ is an open subset of $\g_2$ and,
for every $f \in C_c(\cO_x)$, we have
\begin{equation}
  \label{eq:rao1}
  \int_{\cO_x} f(z)\, dz = c_1 \int_{V + \fn_2} f^K(z_1 + z_2)\phi(z_1) \,
  dz_1 \, dz_2, 
\end{equation}
where
\begin{itemize}
\item $dz_1$ and $dz_2$, resp.,  are Lebesgue measures on
$\g_2$ and $\fn_2$, respectively, 
\item $\phi(z) = |\det(c(z))|^{1/2}$ for
  $c(x) := \ad x\res_{\g_{-1}} \: \g_{-1} \to \g_1,$ and 
\item $f^K= \int_K f \circ \Ad(k)\, dk$, where
  $dk$ is a normalized Haar measure on the compact group
  $K = e^{\ad \fk}$, for a maximal compactly embedded
  subalgebra $\fk \subeq \g$. 
\end{itemize}
\end{thm}

\begin{thm} \mlabel{thm:nilpotorb}
  The invariant measure $\mu_{\cO_x}$ on a nilpotent
  adjoint orbit $\cO_x$ in a reductive Lie algebra is tempered.   
\end{thm}

\begin{prf} First we note that the function
  $\phi$ in Theorem~\ref{thm:rao}
  is of polynomial growth with degree $\shalf \dim \g_1$, i.e.,
  \[ |\phi(x)| \leq c_2 \|x\|^{\frac{\dim \g_1}{2}}\]
  for some $c_2 > 0$. We  assume that
the Cartan involution $\theta$ with $\fk = \Fix(\theta)$ satisfies
$y = \theta(x)$. This can be achieved because every Cartan involution
of the $\fsl_2$-subalgebra spanned by $(h,x,y)$ can be extended to one
on $\g$, i.e., $\fk$ can be chosen to contain a
maximal compactly embedded subalgebra of $\Spann \{h,x,y\}$
(cf.\ \cite[Lemma~I.2]{HNO94}).

If $\kappa(x,y) = \tr(\ad x \ad y)$ is the non-degenerate
Cartan--Killing form on $\g$, then
$\la x,y \ra := - \kappa(x,\theta(y))$ is an
$\Ad(K)$-invariant scalar product on $\g$, defining a
euclidean norm~$\|\cdot\|$. With respect to this scalar
product, $\ad h$ is a symmetric endomorphism, so that its
eigenspaces are orthogonal. In particular,
$\g_1$~and $\g_2$ are orthogonal.

For the $K$-invariant function $f(x) := (\1 + \|x\|^2)^{-k}$, $k \in \N$,
we then have $f^K = f$, so that 
\begin{align*}
 \int_{\cO_z} f(z)\, dz
  &  = c_1 \int_{V + \fn_2} f(z_1 + z_2)\phi(z_1) \,   dz_1 \, dz_2
    = c_1 \int_{V + \fn_2} \frac{\phi(z_1)}{(1 + \|z_1 + z_2\|^2)^{k}}
        \,   dz_1 \, dz_2 \\
  &  \leq c_1 c_2 \int_{V + \fn_2} \frac{\|z_1\|^{\frac{\dim \g_1}{2}}}
    {(1 + \|z_1\|^2  + \|z_2\|^2)^{k}} \,   dz_1 \, dz_2, 
\end{align*}
and this integral is finite if $k$ is sufficiently large.
Here we use that $V\subseteq \g_2$ is open and $dz_1$ is Lebesgue measure 
on $\g_2$.
\end{prf}

\subsection{Mixed orbits in simple Lie algebras} 
\mlabel{subsec:6.3}

We now consider a subalgebra $\fl \subeq \g$, where $\g$ is an admissible
Lie algebra, 
which is the centralizer of an element
$x_\fl \in \ft$. Then $\ft \subeq \fl \subeq \g$ 
implies that $\fl$ contains a Cartan subalgebra
of $\g$, and that $\fl$ is an admissible Lie algebra
because it contains $\ft$, hence intersects the interior of $W_{\rm max}$. 
We write $G := \Inn(\g) \supeq L := \Inn_\g(\fl)$ for the
corresponding adjoint groups. The identity component
$Z_L = \oline{e^{\ad \fz(\fl)}}$ is a torus,  for which
\begin{equation}
  \label{eq:zl-project}
 p_\fl \: \g \to \fl, \quad
 p_\fl(x) = \int_{Z_L} gx \, dg
\end{equation}
is the fixed point projection, where $dg$ stands for a normalized 
Haar measure on $Z_L$. Here we use that the integral
formula obviously is the fixed point projection $\g \to \Fix(Z_L)$,
and that the fixed point space is $\fz_\g(\fz(\fl)) = \fl$
because $\fz(\fl)$ contains the element $x_\fl$ whose centralizer is $\fl$. 
The integral formula for the projection implies that, for any 
$Z_L$-invariant closed convex subsets $C \subeq \g$, we have
\begin{equation}
  \label{eq:pl=intersect}
  p_\fl(C) = C \cap \fl.
\end{equation}

\begin{thm} {\rm(Convexity Theorem for~$p_\fl$)} \mlabel{thm:conv-ell}
  Let $\g$ be an admissible Lie algebra and 
  $\fl \subeq \g$ the centralizer of some element of~$\ft$.
  For $x \in C_{\rm max}^\circ$, we have
\[ p_\fl(\Ad(G)x) \subeq \oline{\conv(\Ad(L)\cW_{\fk} x)}
  + W_{\rm min} \cap \fl.\]
\end{thm}

\begin{prf}  As $p_\ft \circ p_\fl = p_\ft$ follows from $\ft \subeq \fl$,
  we observe that
  \begin{equation}
    \label{eq:pt-pl}
 p_\ft(p_\fl(\co(x))) 
 = p_\ft(\co(x)) \subeq \conv(\cW_{\fk} x) + C_{\rm min}
 \quad \mbox{ for } \quad
x \in C_{\rm max} 
  \end{equation}
  follows from the Convexity Theorem for Adjoint Orbits
  (\cite[Thm.~VIII.1.36]{Ne00}). 
  Let $p_\ft^\fl := p_\ft\res_\fl \: \fl \to \ft$ and
  $L := \la \exp \fl \ra \subeq G$. 
  We consider the closed convex $\Ad(L)$-invariant subset
\[ C_x^L := \{ y \in \fl \: p_\ft(\Ad(L)y)
  \subeq \conv(\cW_{\fk} x) + C_{\rm min}\}
  = \bigcap_{g \in L} \Ad(g) (p_\ft^\fl)^{-1}
  \big(\conv(\cW_{\fk} x) + C_{\rm min}\big). \]
As the closed convex subset $\conv(\cW_{\fk} x) + C_{\rm min}$ of
$\ft$ is invariant under the Weyl group of $\fk \cap \fl$ and
stable under addition of elements in $C_{\rm min,\fl}$, the
Convexity Theorem for Adjoint Orbits, applied to the Lie algebra $\fl$,
shows that $\conv(\cW_{\fk} x) + C_{\rm min}\subeq C_x^L$. 
Hence 
\[ C_x^T := C_x^L \cap \ft
  \ {\buildrel\eqref{eq:pt-pl} \over=}\ \conv(\cW_{\fk} x) + C_{\rm min}.\]
From \eqref{eq:pt-pl} we derive that
the $\Ad(L)$-invariant convex subset $p_\fl(\co(x))$
is contained in $C_x^L.$ 
Therefore it suffices to  show that
\[  C_x^L \subeq \oline{\conv(\Ad(L)\cW_{\fk} x)}  + W_{\rm min} \cap \fl.\]

Next we observe that $x \in C_{\rm max}^\circ$ implies that
$\conv(\cW_{\fk} x) + C_{\rm min} \subeq C_{\rm max}^\circ
\subeq C_{\rm max,\fl}^\circ$, so that
$C_x^L \subeq W_{\rm max,\fl}^\circ$ 
follows from Theorem~\ref{thm:invcon-adm}(d). 
We therefore have 
\[ C_x^L = \Ad(L) C_x^T
\ {\buildrel {(a)} \over =}\ \Ad(L)(\conv(\cW_{\fk} x) + C_{\rm min})
\subeq \conv(\Ad(L)\cW_{\fk} x) + W_{\rm min} \cap \fl.\]
Here (a) follows from the fact that the closed $L$-invariant convex
subset $C_x^L \subeq\fl$ has dense interior, and that all elements
in its interior are conjugate to elements of~$\ft$.
\end{prf}

The following theorem is the version of the Domain Theorem~\ref{thm:dom}
(in the introduction) for reductive Lie algebras.
\begin{thm} \mlabel{thm:6.15-wmax}  Let $\g$ be reductive admissible 
  and $\lambda \in W_{\rm min}^\star$.
  Then
  \[ W_{\rm max}^\circ \subeq
    \Omega_\lambda = \{  x \in \g \: \cL(\mu_\lambda)(x) < \infty\}^\circ,\]
      and equality holds if $\cO_\lambda$ spans $\g^*$.
\end{thm}

\begin{prf} Let $\lambda = \lambda_s + \lambda_n \in W_{\rm min}^\star$
  (Jordan decomposition of $\lambda$, where we include the central
  component in $\lambda_s$ (cf.~\cite[Prop.~1.3.5.1]{Wa72}),  
and write $\fl := \fg^{\lambda_s}$ for the stabilizer Lie algebra of the
semisimple element~$\lambda_s$. This is a reductive Lie subalgebra
of $\g$ (\cite[Prop.~1.3.5.3]{Wa72}). 
We write $L  = \Inn_\g(\fl) \subeq G = \Inn(\g)$
for the integral subgroup corresponding to~$\fl$.

\nin{\bf Step 1:} Let $\beta \in W_{\rm min}^\star \cap \fl^*$ be a
nilpotent element, $\cO_\beta^L = \Ad^*(L)\beta$ and
$\mu_\beta^L$ the $L$-invariant Liouville measure on this orbit.
With respect to the identification of $\fg$ with $\fg^*$, the cone 
$W_{\rm min}^\star$ corresponds to $W_{\rm max}$, so that the closed convex
hull of $\cO_\beta^L \subeq \cO_\beta$ contains no affine line. 
Further, $\mu_\beta^L$ is tempered by Theorem~\ref{thm:nilpotorb},
so that Proposition~\ref{prop:dom-lapl-temp} shows that its Laplace
transform is defined on the open cone
$B(\cO_\beta^L)^\circ \supeq W_{\rm max}^\circ$.
Here we use that $\beta \in W_{\rm max,\fl}^\star$ follows from 
\cite[Thm.~III.9]{HNO94}, applied to the semisimple Lie algebra~$[\fl,\fl]$. 
We thus have 
\begin{equation}
  \label{eq:lapl615}
\cL(\mu_\beta^L)(x) = \int_{\cO_\beta^L} e^{-\alpha(x)}\,
  d\mu_\beta^L(\alpha)
  \quad \mbox{ for }  \quad x \in W_{\rm max,\fl}^\circ.
\end{equation}

As $C_{\rm max,\fl} \supeq C_{\rm max}$ is a consequence of 
$\Delta_{p,\fl} \subeq \Delta_p$, it follows from
$(W_{\rm max} \cap \fl) \cap \ft = C_{\rm max}$ that
\begin{equation}
  \label{eq:wmaxlandg}
  W_{\rm max,\fl}
  = \oline{\Ad(L) C_{\rm max,\fl}}
  \supeq  \oline{\Ad(L) C_{\rm max}} = W_{\rm max} \cap \fl.
\end{equation}

\nin {\bf Step 2:} The function \eqref{eq:lapl615}
on $W_{\rm max}^\circ \cap \fl$
is decreasing in the direction of
$W_{\rm min} \cap \fl$. In fact, for $x \in W_{\rm max}^\circ \cap \fl$
and $y \in W_{\rm min} \cap \fl$, we have for any
$\alpha \in \cO_\beta^L \subeq \cO_\beta \subeq W_{\rm min}^\star$ that
$\alpha(x + y) \geq \alpha(x),$ so that
\[ \cL(\mu_\beta^L)(x+y)
  = \int_{\cO_\beta^L} e^{-\alpha(x+y)}\,  d\mu_\beta^L(\alpha)
\leq \int_{\cO_\beta^L} e^{-\alpha(x)}\,  d\mu_\beta^L(\alpha)
= \cL(\mu_\beta^L)(x).\]

\nin {\bf Step 3:} For $x \in C_{\rm max}^\circ$, we obtain with
the Convexity Theorem~\ref{thm:conv-ell} 
\[ p_\fl(\Ad(G)x) \subeq \conv(\Ad(L)\cW_{\fk} x) + W_{\rm min} \cap \fl,\]
and hence, identifying $\fl^*$ with the  subspace $[\fz(\fl),\g]^\bot
\subeq \g^*$, 
\begin{equation}
  \label{eq:lapl-decr}
 \cL(\mu_\beta^L)(\Ad(g)x) 
= \cL(\mu_\beta^L)(p_\fl(\Ad(g)x))
\leq \sup \cL(\mu_\beta^L)(\Ad(L)\cW_{\fk} x) 
= \max\cL(\mu_\beta^L)(\cW_{\fk} x).
\end{equation}

\nin {\bf Step 4:} We now apply the preceding discussion to
the nilpotent Jordan component $\beta= \lambda_n$,
which is also contained in $W_{\rm min}^\star$ by \cite[Cor.~B.2]{NO22}.
The invariant measure $\mu_\lambda$ on $\cO_\lambda$ takes by
\cite[p.~510]{Rao72} the form
\[ \mu_\lambda = \int_{G/L} g_*\mu_{\lambda}^L\, d\mu(gL). \]
For its Laplace transform, we find on $x \in C_{\rm  max}^\circ$
with \eqref{eq:lapl-decr} the estimate 
  \begin{align} \label{eq:lapl-prod}
 \cL(\mu_\lambda)(x) &
=  \int_{G/L} e^{-\alpha(x)}\ d(g_*\mu_{\lambda}^L)(\alpha)\, d\mu(gL)
    \notag \\
& = \int_{G/L} e^{-\lambda_s(g^{-1}x)} \cL(\mu_{\lambda_n}^L)(g^{-1}x)\, d\mu(gL)
    \notag \\
& \ {\buildrel \eqref{eq:lapl-decr} \over  \leq}\  \int_{G/L} e^{-\lambda_s(g^{-1}x)}
      \cdot \Big(\max  \cL(\mu_{\lambda_n}^L)(\cW_{\fk} x)\Big)\, d\mu(gL) \notag  \\
    &= \max  \cL(\mu_{\lambda_n}^L)(\cW_{\fk} x)
      \cdot \int_{G/L} e^{-\lambda_s(g^{-1}x)}\, d\mu(gL)
      = \max  \cL(\mu_{\lambda_n}^L)(\cW_{\fk} x)
      \cdot \cL(\mu_{\lambda_s})(x).
  \end{align}
  Since $\lambda_s \in W_{\rm min}^\star$ (cf.\ \cite[Cor.~B.2]{NO22})
  corresponds under the duality $\g \to \g^*$ to an elliptic
  element, it is conjugate to an element in
  $C_{\rm min}^\star \subeq \ft^*$, hence admissible
  (Proposition~\ref{prop:adm-orbits}). 
  From Theorem~\ref{thm:II.10} we know that
  $\cL(\mu_{\lambda_s})(x) < \infty$ because
  $x \in C_{\rm max}^\circ$.  Further, the $\cW_\fk$-invariance of
  $C_{\rm max}$ entails that 
  $\cW_{\fk} x \subeq C_{\rm max}^\circ \subeq C_{\rm max,\fl}^\circ$, 
  and $\cL(\mu_{\lambda_n})$ is finite on $W_{\rm max,\fl}^\circ$
  by \eqref{eq:lapl615}. With \eqref{eq:lapl-prod} it follows with
  Theorem~\ref{thm:invcon-adm}(a) that
  $W_{\rm  max}^\circ =  \Ad(G) C_{\rm max}^\circ 
  \subeq  \Omega_\lambda$. Here the last inclusion follows from
the $G$-invariance of $\Omega_\lambda$. 

    To verify the last statement,
    we write $\g = \g_0 \oplus \g_1$, where
  $\g_0 = \cO_\lambda^\bot$, so that $\cO_\lambda$ spans~$\g_1^*$.
  As $W_{\rm max}^\g = W_{\rm max}^{\g_0} \oplus W_{\rm max}^{\g_1}$
  and $\g_0 \subeq H(\Omega_\lambda)$, we may assume that $\cO_\lambda$
  spans $\g^*$. We have seen above that $W_{\rm max}^\circ
  \subeq \Omega_\lambda$, and the converse follows from
  Corollary~\ref{cor:4.1}.
\end{prf}

\begin{prop} \mlabel{prop:4.8}
  Suppose that $\g$ is reductive admissible, and that  
  $\lambda \in W_{\rm min}^\star$ is such that $\cO_\lambda$ spans~$\g^*$.
Then   $D_{\mu_\lambda}$ is open, hence equal to $\Omega_\lambda$.   
\end{prop}

\begin{prf} In view of the uniqueness of $\Delta_p^+$ 
  in Corollary~\ref{cor:4.1}, we derive from 
Theorem~\ref{thm:6.15-wmax} that
$\Omega_\lambda = W_{\rm max}^\circ$.
On the other hand, Theorem~\ref{thm:5.9}(c) shows that
$D_{\mu_\lambda} \subeq W_{\rm max}^\circ$, hence that
$D_{\mu_\lambda} = \Omega_\lambda$.
\end{prf}

 \subsection{Coadjoint orbits in semidirect sums} 
\mlabel{subsec:6.4}

Consider a semidirect sum $\g = \fu \rtimes \fl$ 
and a corresponding $1$-connected Lie group $G = U \rtimes L$.
We write linear functionals on $\g$ as 
 \begin{equation}
   \label{eq:jordan-dec}
 \lambda = \lambda_\fu + \lambda_\fl \quad \mbox{
   with } \quad \lambda_\fu \in \fu^* \cong \fl^\bot
 \quad \mbox{ and }  \quad
   \lambda_\fl \in \fl^* \cong \fu^\bot.
 \end{equation}
 Then $U$ acts trivially on $\fl^* \cong \fu^\bot  \cong (\g/\fu)^*$
 because $\fu \trile \g$ is an ideal.

 \begin{lem} \mlabel{lem:prod-orbit} Suppose that
   $\g = \fu \rtimes \fl$ and that 
   $\lambda = \lambda_\fu+\lambda_\fl \in \g^*$ is decomposed
   accordingly. We assume that 
 \begin{itemize}
 \item $\lambda_\fu$ is fixed by $L$, and that 
 \item the stabilizer group
   $U^{\lambda_\fu} = \{ u \in U \: \Ad^*(u)\lambda_\fu = \lambda_\fu\}$
   is connected.
 \end{itemize}
 Then the following assertions hold: 
 \begin{itemize}
 \item[\rm(a)] $\cO_\lambda = \cO_{\lambda_\fu} + \cO_{\lambda_\fl}
   = \Ad^*(U) \lambda_\fu +\Ad^*(L)\lambda_\fl$, where
   $U$ acts trivially on $\fl^* \cong (\g/\fu)^*$.
 \item[\rm(b)] The addition map defines a $G$-equivariant
   symplectic diffeomorphism
$\add \: \cO_{\lambda_\fu} \times \cO_{\lambda_\fl} \to \cO_\lambda.$ 
 \item[\rm(c)] $p_{\fu} \: \cO_{\lambda_\fu} \to \fu^*$
   is a diffeomorphism onto a coadjoint orbit of $U$ in $\fu^*$. 
 \end{itemize}
 \end{lem}

\begin{prf} From \cite[Prop.~VIII.1.2]{Ne00} (a)-(c) follow,
  with the exception of the diffeomorphism
  \begin{equation}
    \label{eq:psi-sum}
 \Psi \:  \cO_{\lambda_\fu} \times \cO_{\lambda_\fl} 
 \to \cO_\lambda \subeq \g^*, \quad (\alpha,\beta) \mapsto \alpha + \beta
  \end{equation}
  under (b) being symplectic. 
  To see that $\Psi$ is symplectic,
 we first note that the symplectic product structure is $G$-invariant and that
 the product space is homogeneous. The tangent space
 of $\cO_\lambda$ in  $\lambda$ is
 \[ \lambda \circ \ad \g
= \lambda \circ (\ad \fu + \ad \fl) 
= \lambda_\fu \circ \ad \fu + \lambda_\fl \circ \ad \fl,\] 
and the sum of these two subspaces is direct.
For $x_u, y_u \in \fu$ and $x_l, y_l  \in \fl$, we have
\[ \lambda([x_u + x_l, y_u + y_l])
= \lambda_\fu([x_u + x_l, y_u + y_l])
+ \lambda_\fl([x_u + x_l, y_u + y_l])
= \lambda_\fu([x_u, y_u])
+ \lambda_\fl([x_l, y_l]) \]
(cf.\ \eqref{eq:omegaalpha}), and this shows that
\[ T_\lambda(\cO_\lambda)
  \cong T_{\lambda_\fu}(\cO_{\lambda_\fu}^U)
  \oplus T_{\lambda_\fl}(\cO_{\lambda_\fl}^L) \]
as symplectic vector spaces. This completes the proof of (b).
In particular, $\cO_{\lambda_\fu} \times \cO_{\lambda_\fl}$
is a homogeneous Hamiltonian $G$-space whose momentum map is
given by $\Psi$.
\end{prf}

\begin{lem} \mlabel{lem:lambdav=0}
If $\g = \fu \rtimes \fl  = (\fz +  V) \rtimes \fl$
is admissible and decomposed as in {\rm Theorem~\ref{thm:spind}},
and $\cO_\lambda$ spans~$\g^*$,
then $\cO_\lambda$ contains a functional $\alpha$ vanishing on $V$,
i.e., $\alpha_V = 0$. 
\end{lem}

\begin{prf} If $V = \{0\}$ there is nothing to show.
  So we assume that $V \not=\{0\}$.
  That $\cO_\lambda$ spans $\g^*$ implies in particular
  that $\cO_\lambda \res_\fz = \{\lambda_\fz = \lambda\res_\fz\}$
  (central elements define constant Hamiltonian functions on $\cO_\lambda$) 
  separates points on
  $\fz$, so that $\dim \fz = 1$ and $\lambda_\fz \not=0$.
By admissibility of $\g$, 
\begin{equation}
  \label{eq:symp-form}
  \Omega(v,w) := \lambda_\fz([v,w])
\end{equation}
is a symplectic form on $V$ (see \eqref{eq:omegaf}
in Subsection~\ref{subsec:1.4}), so that 
\begin{equation}
  \label{eq:u-orb}
  \lambda_\fu(e^{\ad v}w)
  = \lambda_V(w) + \frac{1}{2} \lambda_\fz([v,w])
  = \lambda_V(w) + \frac{1}{2} \Omega(v,w) 
\end{equation}
shows that, if we choose $v$ in such a way that
$\lambda_V = - \frac{1}{2} \Omega(v,\cdot),$ 
then we obtain a functional in $\cO_\lambda$ that vanishes on $V$.
\end{prf}

In view of the preceding lemma, we may assume that $\lambda_V = 0$. Then
$\lambda = \lambda_\fz + \lambda_\fl$ with
$\lambda_\fz = \lambda\res_\fz$ and $\lambda_\fl = \lambda\res_\fl$, 
and $L$ fixes $\lambda_\fz$. From \eqref{eq:u-orb}
we derive that, if $u = (z,v) \in U$ fixes $\lambda_\fz$, then $v = 0$
because $\Omega$ is non-degenerate, hence 
\[ U^{\lambda_\fz} = \{ u \in U \:\Ad^*(u)\lambda_\fz = \lambda_\fz \}
  = Z(G)_e \]
is connected. Therefore Lemma~\ref{lem:prod-orbit} applies.

\begin{prop} \mlabel{prop:7.3}
  For $\lambda = \lambda_\fz + \lambda_\fl
  \in W_{\rm min}^\star$, we have
  \begin{itemize}
  \item[\rm(a)] $\mu_\lambda = \mu_{\lambda_\fz} * \mu_{\lambda_\fl}$
    is a convolution product. 
  \item[\rm(b)] $\cL(\mu_\lambda) = \cL(\mu_{\lambda_\fz})\cL(\mu_{\lambda_\fl})$. 
  \item[\rm(c)] $D_{\mu_\lambda}
  = D_{\mu_{\lambda_\fz}} \cap D_{\mu_{\lambda_\fl}}
  = \Omega_{\mu_{\lambda_\fz}} \cap D_{\mu_{\lambda_\fl}}$. 
  \end{itemize}
\end{prop}

\begin{prf}
  (a) follows from the fact that the addition map $\Psi$
  from \eqref{eq:psi-sum} in the proof of
  Lemma~\ref{lem:prod-orbit} is the momentum map of the
  $G$-action on the symplectic product $\cO_{\lambda_\fz} \times \cO_{\lambda_\fl}$.

  \nin (b) follows from (a). 

  \nin (c) follows from (b) and the fact that
  $D_{\mu_{\lambda_\fz}}$ is open since the cone of
  positive definite symmetric matrices is open in the
  space of all symmetric matrices (Lemma~\ref{lem:5.7}). 
\end{prf}

\begin{lem} \mlabel{lem:when-span}
  Let $\g = \fu \rtimes \fl = (\fz \times V) \rtimes \fl$
  be admissible, non-reductive,
  $\lambda \in W_{\rm min}^\star$, and $\dim\fz(\g) = 1$.
  We write $\fl = \fl_0 \oplus \fl_1$ with
  $\fl_1 = \fz_\fl(V)$ and, accordingly,
  \[ \lambda = \lambda_\fz + \lambda_V + \lambda_\fl^0 + \lambda_\fl^1
    \in \fz^* \oplus  V^* \oplus \fl_0^* \oplus  \fl_1^*.\]
  Then $\cO_\lambda$ spans $\g^*$ if and only if
  \begin{equation}
    \label{eq:span-condi}
\lambda_\fz \not=0 \quad \mbox{ and } \quad
\cO_{\lambda_\fl^1} \ \mbox{ spans } \ \fl_1^*.
  \end{equation}
\end{lem}

\begin{prf} Suppose first that $\cO_\lambda$ spans $\g^*$. As
$\g = (\fu \rtimes \fl_0) \oplus \fl_1$ 
  is a direct Lie algebra sum, $\cO_{\lambda_\fl^1}$ spans $\fl_1^*$.
  Further, central elements define constant functions on $\cO_\lambda$,
  and since $\fz \not=\{0\}$,   we must have $\lambda_\fz \not=\{0\}$.

  Suppose, conversely, that these two conditions are satisfied.
  We have to show that any
  \[ (z,w,y_0 + y_1) \in \cO_\lambda^\bot \] 
  vanishes. Here $y = y_0 + y_1 \in \fl$ is the decomposition
  into $\fl_0$ and $\fl_1$-component. 
  As $\lambda_\fz \not=0$, $\lambda \circ e^{\ad V}$ contains an element
  with $\lambda_V =\{0\}$ and the same $\fz$-component
  (Lemma~\ref{lem:lambdav=0}).
  We may therefore assume that $\lambda_V = 0$,
  so that $\lambda_\fu =  \lambda_\fz$.  Now  
\begin{equation}
  \label{eq:van1}
 0     = \la \Ad^*(u)\lambda_\fu + \Ad^*(\ell)\lambda_\fl, (z,w,y) \ra
    = \la \Ad^*(u)\lambda_\fu, (z,w,y) \ra
    + \la \Ad^*(\ell)\lambda_\fl, y\ra
\end{equation}
  for all $u \in U$ and $\ell \in L$. 
  It follows in particular
  that $u \mapsto \la \Ad^*(u)\lambda_\fu, (z,w,y) \ra$ is constant. With
\[   \Omega(v,w) := \lambda_\fz([v,w]),\]
this implies that the
  following expression does not depend on $v$:
  \[ \lambda_\fz(e^{\ad v}(z,w,y))
    = \lambda_\fz(z) +\lambda_\fz([v,w]) + \lambda_V([v,y])
    + \frac{1}{2}\lambda_\fz([v,[v,y]])
    = \lambda_\fz(z)  +\lambda_\fz([v,w]) + \frac{1}{2}\Omega(v,[v,y]).\]
  Here we use that
  \begin{align} \label{heisact}
   e^{\ad v}(z,w,y)
&=    (z,w,y) + [v,(z,w,y)] + \frac{1}{2}[v,[v,(z,w,y)]]\\
&=    (z,w,y) + ([v,w], [v,y],0) + \frac{1}{2}([[v,[v,y]],0,0)\notag \\
&  =   \Big(z + [v,w],  \frac{1}{2}[[v,[v,y]],w + [v,y], y\Big).
\notag  \end{align}
This is a polynomial in $v$, for which the summands in
\eqref{heisact} are of degree
$0$, $1$ and $2$, respectively. 
Since it is constant, the 
  homogeneous terms of degree $1$ and $2$ vanish, i.e., 
 \[ \Omega(V,w)  = \{0\}\quad \mbox{ and } \quad 
   \Omega(v,[v,y]) = 0 \quad \mbox{ for } \quad v \in V.\]
 We conclude that $w = 0$ because $\Omega$ is non-degenerate, 
 and also by polarization that $\Omega(w,[v,y]) = 0$ for all $w,v \in V$,
 so that $[y,V] = \{0\}$, i.e., $y \in \fl_1$.

 The relation \eqref{eq:van1} now reduces to
 \[  0   = \lambda_\fz(z) +  \la \Ad^*(\ell)\lambda_\fl^1, y\ra
   \quad \mbox{ for } \quad \ell \in L.\]
 As $\Ad^*(L) \lambda_\fl^1$ spans $\fl_1$,
 the fact that the Hamiltonian function $H_y$ is constant on this orbit
 implies that $y \in \fz(\fl_1) \subeq \fz(\g) \subeq \fu$, hence that
 $y = 0$. Finally $\lambda_\fz(z) = 0$ entails $z = 0$
 because $\fz$ is $1$-dimensional and $\lambda_\fz \not=0$.
 \end{prf}

 \begin{ex} For the non-reductive admissible Lie algebra
   $\g = \hsp(V,\Omega)$ we have $\fl = \fl_0$.
   For the functional  $\lambda = \lambda_\fz = \ev_0$,
   the orbit $\cO_\lambda$ spans $\g^*$, but $\lambda_\fl = 0$.    
 \end{ex}

 \subsection{The general case}
\mlabel{subsec:6.5}
 
 For the following theorem, we recall that, replacing $\g$ by
 $\g/\fn$ for $\fn := \cO_\lambda^\bot$, 
 we may always reduce to the situation where  $\cO_\lambda$ spans $\g^*$.
 
 \begin{thm} \mlabel{conj:4.4}
   Let $\cO_\lambda \subeq \g$ be a coadjoint orbit spanning $\g^*$.
   Then the following assertions hold:
   \begin{itemize}
   \item[\rm(a)] If $D_{\mu_\lambda}\not=\eset$, then
     $\g$ is admissible   and there exists an adapted positive system
     $\Delta^+$ with $C_{\rm min}$ pointed and contained in $C_{\rm max}$
     for which $\lambda \in W_{\rm min}^\star$.
   \item[\rm(b)] Suppose that $\ft \subeq \g$ is a compactly embedded
     Cartan subalgebra and $\Delta^+$ an adapted positive system 
     with $C_{\rm min}$ pointed and contained in $C_{\rm max}$.
Then, for $\lambda \in W_{\rm min}^\star$, 
   \begin{itemize}
   \item[\rm(1)] $D_{\mu_\lambda} = \Omega_\lambda = W_{\rm max}^\circ$.
   \item[\rm(2)] $Q \: \Omega_\lambda \to C_\lambda,
  Q(x) =  \frac{1}{Z_\lambda(x)}
  \int_{\g^*} \alpha e^{-\alpha(x)}\, d\mu_\lambda(\alpha),$ defines a diffeomorphism
     from $\Omega_\lambda/\fz(\g)$ onto $C_\lambda^\circ= \conv(\cO_\lambda)^\circ$. 
   \end{itemize}
   \end{itemize}
\end{thm}

\begin{prf} (a) follows from Corollary~\ref{cor:4.1}. 

\nin (b.1)   As in Lemma~\ref{lem:when-span}, we write
  \[ \g = (\fu \rtimes \fl_0) \oplus \fl_1 \quad \mbox{ with } \quad
    \fl_1 = \fz_\fl(\fu).\]
  If $V = [\ft,\fu]\not=\{0\}$, then $\fz = [V,V] \not=\{0\}$
  because the bracket $V \times V  \to \fz$ is a non-degenerate
  vector-valued alternating form (Theorem~\ref{thm:spind}).
  Since $\cO_\lambda$ spans $\g^*$ and restricts to a singleton
  on $\fz$, it follows that $\dim \fz = 1$ and
  $\lambda_\fz = \lambda\res_\fz \not=0$.
  In view of Lemma~\ref{lem:lambdav=0}, we may assume that $\lambda_V = 0$,
  so that
  \[ \lambda = \lambda_\fz + \lambda_{\fl_0} + \lambda_{\fl_1}.\] 
  Proposition~\ref{prop:7.3} shows that
  \begin{equation}
    \label{eq:intersect1}
    D_{\mu_\lambda} = \Omega_{\lambda_\fz} \cap D_{\mu_{\lambda_\fl}},
  \end{equation}
    and $\fl = \fl_0 \oplus \fl_1$ entails that
    \[ D_{\mu_{\lambda_\fl}}
      = D_{\mu_{\lambda_{\fl_0}}} \cap D_{\mu_{\lambda_{\fl_1}}}.\]
          Lemma~\ref{lem:when-span} further implies that
  $\cO_{\lambda_\fl^1}$ spans $\fl_1^*$.
  With the ideal $\fl_{0,0} := \cO_{\lambda_{\fl_0}}^\bot \cap\fl_0
  = \cO_{\lambda_\fl}^\bot \cap \fl \trile \fl$ and a
  complementary ideal $\fl_{0,1}$, we now have
  $\fl_0 = \fl_{0,0} \oplus \fl_{0,1}$.
  We thus obtain 
  the direct sum decomposition $\fl = \fl_{0,0} \oplus \fl_{0,1} \oplus \fl_1.$
  Then $\ft_\fl = \fl \cap \ft$ and the minimal and maximal cones
  are adapted to this decomposition and $\cO_{\lambda_\fl}$ spans
  the dual of the ideal $\fl_{0,1} \oplus \fl_1$ 
   of the admissible Lie algebra   $\fl$, which also is admissible. 
   From Lemma~\ref{lem:cmin-div} we derive that $\lambda_\fl \in
   W_{\rm min,\fl}^\star$, so that Theorem~\ref{thm:6.15-wmax} yields 
   $W_{\rm max,\fl}^\circ \subeq \Omega_{\lambda_\fl}$ 
  because the cone $W_{\rm max,\fl}$ is adapted to the decomposition of~$\fl$.
  Thus
  \[ \Omega_{\lambda_\fl}
   \supeq \fl_{0,0}
    \oplus  W_{\rm max,\fl_{0,1}}^\circ 
    \oplus  W_{\rm max,\fl_{1}}^\circ,\]
  and with Proposition~\ref{prop:4.8}, applied
  to the ideal $\fl_{0,1} \oplus \fl_1$,   this further leads to 
$D_{\mu_{\lambda_\fl}} = \Omega_{\lambda_\fl}$.
In view of \eqref{eq:intersect1}, it follows that $D_{\mu_\lambda}$ 
is open.

\nin (b.2) Since $\lambda_\fz \in C_{\rm min}^\star$
(Lemma~\ref{lem:cmin-div}) is admissible, 
Theorem~\ref{thm:II.10} implies that
$W_{\rm  max}^\circ \subeq \Omega_{\lambda_\fz}$.
Further, $\fu$~acts trivially on $\fl^*$ and
the projection $\g \to \fl$ maps $W_{\rm max}$ into $W_{\rm max,\fl}$, 
so that 
\[ \Omega_{\lambda_\fl}\supeq \fu +  W_{\rm max,\fl}^\circ 
  \supeq W_{\rm max}^\circ.\]
We thus find with \eqref{eq:intersect1} that 
$W_{\rm max}^\circ \subeq \Omega_\lambda$.
As  $\cO_\lambda$ spans $\g^*$, Theorem~\ref{thm:5.9} implies
that $D_{\mu_\lambda} \subeq W_{\rm max}^\circ$, so that we actually have
the equality $D_{\mu_\lambda} = \Omega_\lambda = W_{\rm max}^\circ$.
In particular, $D_{\mu_\lambda}$ is open, so that
(b) follows from Theorem~\ref{thm:conv-lapl}.
\end{prf}

The preceding theorem brings us full circle in the classification
of coadjoint orbits $\cO_\lambda$ for which $D_{\mu_\lambda} \not=\eset$,
and we have actually seen that (after some reduction), this
$D_{\mu_\lambda} \not=\eset$.
We had already seen above, that, factorizing the ideal
$\cO_\lambda^\bot$, we may always assume that the orbit spans
$\g^*$. Then Theorem~\ref{conj:4.4}(a) tells us 
where these functionals $\lambda$ can be found, namely in some
$W_{\rm min}^\star$, and part (b) shows that all these functionals
actually satisfy $\Omega_\lambda \not=\eset$.

At this point one should note that the positive system
$\Delta^+_p$ is uniquely determined by $\lambda$,
and that, given $\g$, there are only finitely many such systems. 

\subsection{Temperedness of the Liouville measures}

In this subsection we show that, whenever $D_{\mu_\lambda} \not=\eset$,
the Liouville measure on $\cO_\lambda$ is tempered, i.e.,
defines a tempered distribution on~$\g^*$.

\begin{thm} \mlabel{thm:temp}
  Let $\cO_\lambda \subeq \g^*$ be a coadjoint orbit for which
  $D_{\mu_\lambda} \not=\eset$. Then
  the Liouville measure on $\cO_\lambda$ is tempered.
\end{thm}

\begin{prf} We first observe that we may assume that $\cO_\lambda$
  spans $\g^*$, so that Theorem~\ref{thm:5.9} entails that
  $\g$ is admissible and that $\lambda \in W_{\rm min}^\star$
  for an adapted positive system for which
  $C_{\rm min}$ is pointed and contained in $C_{\rm max}$. 

  In view of Lemma~\ref{lem:lambdav=0},
  Proposition~\ref{prop:7.3} provides a decomposition
  $\lambda = \lambda_\fz + \lambda_\fl$ for which
  \[ \cL(\mu_\lambda) = \cL(\mu_{\lambda_\fz}) \cL(\mu_{\lambda_\fl}).\]
  Since $\cO_{\lambda_\fz}$ is admissible, there exists an $N \in \N$
  such that, for $x \in C_{\rm max}^\circ$, we have 
  \[ c := \lim_{t \to 0+} \cL(\mu_{\lambda_\fz})(tx) t^{-N}  \]
  exists.

  From Step 4 in the proof of Theorem~\ref{thm:6.15-wmax},
  we further obtain a Jordan decomposition
  $\lambda_\fl= \lambda_s + \lambda_n$ such that
  \[  \cL(\mu_\lambda)(x)
      = \max  \cL(\mu_{\lambda_n}^L)(\cW_{\fk} x)
      \cdot \cL(\mu_{\lambda_s})(x). \] 
    As $\cO_{\lambda_s}$ is admissible, $\mu_{\lambda_s}$
    is tempered by  Corollary~\ref{cor:5.4}.
    Further $\mu_{\lambda_n}$ is tempered by
    Theorem~\ref{thm:nilpotorb}.
    So Proposition~\ref{prop:dom-lapl-temp} yields a $k \in \N$ for which
    \[ \limsup_{t \to 0+} \cL(\mu_{\lambda_\fl})(tx) t^k < \infty.\] 
Therefore 
\[ \limsup_{t \to 0+} \cL(\mu_{\lambda})(tx) t^{k+N} < \infty,\]
so that Proposition~\ref{prop:dom-lapl-temp} shows that
$\mu_\lambda$ is tempered.   
\end{prf}

\begin{rem}
In \cite{Ch90, Ch96}, Charbonnel shows that, for any connected Lie group $G$,
the Liouville measure on a closed coadjoint orbit is tempered.
This is already claimed in \cite[Thm.~1.8]{Ch90},
but the argument in \cite{Ch90} only worked under 
the assumption that the Lie algebra
$\ad \g$ is stable under Jordan decomposition.
This gap was filled in \cite{Ch96}.
For the connection between the Fourier transforms
of closed coadjoint orbits and characters of unitary representations,
we refer to \cite{BV83} and \cite{Ne96a}.

It is quite plausible that, for a reductive Lie algebra,
all Liouville measures are tempered. For nilpotent orbits we saw this in
Theorem~\ref{thm:nilpotorb}, and for orbits with non-trivial geometric
temperatures, it follows from Theorem~\ref{thm:temp}.
We expect
that the methods developed in \cite{dCl91} can be used
to prove that this is true; as suggested in an email
from Yoshiki Oshima.
\end{rem}

\section{Disintegration of invariant measures}
\mlabel{sec:disint}

In this section we take a closer look at the $\Ad^*(G)$-invariant
measures $\mu$ on $\g^*$ that arise from general Hamiltonian
$G$-actions with non-trivial geometric temperature.
We know already from Theorem~\ref{thm:5.9} that
we may assume that $\g$ is admissible and that 
$\Psi(M) \subeq W_{\rm min}^\star$ holds for an adapated 
positive system $\Delta^+$ of roots with respect to a 
compactly embedded Cartan subalgebra $\ft$, for which
$C_{\rm min}$ is pointed and contained in $C_{\rm max}$.

Our strategy is to use results on algebraic groups, which
is based on the following observation. 

\begin{lem} \mlabel{lem:8.1}
We consider the  action of the closure $G_c := \oline{\Ad(G)}$
on the corresponding invariant cone $W_{\rm min}^\star \subeq \g^*$.
Then the following assertions hold:
\begin{itemize}
\item[\rm(a)] Write $\g = \fu \rtimes \fl$ with
  $\ft = \fz(\g) \oplus \ft_\fl$, so that $\fz(\fl) \subeq \ft$.
  The group $Z_L := \oline{e^{\ad \fz(\fl)}}$ is a torus
  and $G_c = \Ad(G) Z_L$ is the identity component,
  with respect to the Lie group topology, of an algebraic group,
namely the Zariski closure of $\Ad(G)$. 
\item[\rm(b)]   For $\lambda \in W_{\rm min}^\star$, the coadjoint
  orbit $\cO_\lambda$ is also invariant under $G_c$.
\end{itemize}
\end{lem}

\begin{prf} (a) Since $\ft$ is compactly embedded, 
  $Z_L$ is a compact group, hence in particular
  algebraic. Let $a(\ad \g)$ denote the Lie algebra
  of the Zariski closure of $\Ad(G)$, i.e., the algebraic hull of
  $\ad \g$. In view of \cite[Prop.~I.6(iii)]{Ne94},
  $\ad([\g,\g]) = [\ad \g,\ad \g]$ is algebraic and we have seen above
  that $\L(Z_L)$ is also algebraic. We also have
  \[ \g =  [\g,\g] + \ft = [\g,\g] + \fz(\fl) \]
  because
  \[ \ft = \fz(\g) + \fz(\fl) +  \ft \cap [\fl,\fl]
    \subeq [\g,\g]  + \fz(\fl).\] 
  Therefore
$\ad([\g,\g]) + \L(Z_L)$ is the Lie algebra of an algebraic group 
(\cite[Prop.~I.6(ii)]{Ne94}), and since this is the Lie algebra of
$\Ad(G) Z_L = \Ad(G) \oline{e^{\ad \ft}}= \oline{\Ad(G)}$ 
(cf. \cite[Thm.~14.5.3(ii)]{HN12}), the assertion follows.

\nin  (b) Using that $\cO_\lambda \subeq (\g/\fn)^*$ for
  $\fn = \cO_\lambda^\bot$, we may assume that $\cO_\lambda$ spans $\g^*$.
  Then $\dim \fz(\g) \leq 1$ and  Lemma~\ref{lem:lambdav=0}, combined
  with   Proposition~\ref{prop:7.3}, provides a decomposition
  $\lambda = \lambda_\fz + \lambda_\fl$, according to
  $\g = \fu \rtimes \fl$. Here $\lambda_\fz \in \ft^*$ is fixed by $L$, and 
  $\lambda_\fl \in \fl^*$, where $\fl$ is reductive.
  Therefore $Z_L$ fixes $\lambda_\fl$. This shows that
  $G_c = \Ad(G) Z_L$ leaves $\cO_\lambda = \cO_{\lambda_\fz} + \cO_{\lambda_\fl}$ 
  invariant.   
\end{prf}

\begin{thm} \mlabel{thm:disint} {\rm(Disintegration Theorem)} 
  Let $\mu$ be an $\Ad^*(G)$-invariant measure on the closed
  convex cone $C := W_{\rm min}^\star$ associated to an adapted positive
  system with $C_{\rm min}$ pointed and contained in $C_{\rm max}$. 
We assume that there exists an $x \in \g$ with 
  $\cL(\mu)(x) < \infty$. Then there exists a
  measure $\nu$ on the Borel quotient $C/G$ for which
\[ \mu = \int_{C/G} \mu_\lambda\,  d\nu([\lambda]).\]
\end{thm}

\begin{prf}  {\bf Step 1:} As $\cL(\mu)(x) < \infty$, the measure
  $\tilde \mu := e^{-H_x} \mu$ is finite, so that $\mu$ is a Radon measure, i.e.,
  finite on compact sets. Therefore the same argument as in the proof
  of Theorem~\ref{thm:measlemb}(a) shows that the stabilizer of $\mu$
  in $\GL(\g^*)$ is closed, hence contains $G_c$ from
  Lemma~\ref{lem:8.1}.

  \nin {\bf Step 2:} (Chevalley's Theorem) Let $H$ be an affine algebraic group
acting regularly on an affine algebraic variety $X$
and write $H_e$ for the identity component in the Lie group topology.
Then the Borel space $X/H_e$ is countably separated, i.e., the $\sigma$-algebra
of $H_e$-invariant Borel sets is countably generated. 
This result was
never published by Chevalley himself, but a sketch
of the proof and corresponding references are given
on page 183 of \cite{Dix66}; see also the introduction
of \cite{Dix57} and \cite[Thm.~VI.10]{Fa00}.  

Applying Pukanszky's Theorem \cite[p.~50]{Pu72} to the
action of the Zariski closure $H$ of $\Ad^*(G)$ on~$\g^*$,
considered as the unitary dual of the additive group $(\g,+)$, 
it implies that the orbit space $\g^*/H_e = \g^*/G_c$ is countably separated,
so that $S := C/G = C/G_c$ (Lemma~\ref{lem:8.1}) is also countably separated.
Thus \cite[Thm.~VI.11]{Fa00} implies the existence of a
Borel cross section. We may thus consider $S$ as a subset of $C$,
meeting every $G$-orbit exactly once. We write
\[ q \: C \to S \quad \mbox{ with } \quad 
  q(\cO_\lambda) = \{\lambda\}, \quad \lambda \in S, \]
for the corresponding quotient map. 

\nin {\bf Step 3:} The measure $\mu$ on $C$ is Radon, hence in
particular $\sigma$-finite and equivalent to the finite measure
$\tilde\mu$ from above. We also note that $\tilde \mu$ is quasi-invariant
under $\Ad^*(G)$. 

Now $\tilde\nu := q_*\tilde\mu$ is a finite positive Borel measure
on $S$ and the Disintegration Theorem \cite[Thm.~I.27]{Fa00}
implies the existence of a family of finite measures 
$(\tilde\mu_\lambda)_{\lambda \in S}$ such that 
  \begin{itemize}
\item[\rm(1)] For each Borel set $E \subeq C$,  the map
    $S \to [0,\infty], \lambda \mapsto \tilde\mu_\lambda(E)$ is measurable and
    \begin{equation}
      \label{eq:tildemu-dec}
      \tilde\mu(E) = \int_{S} \tilde\mu_\lambda(E)\, d\tilde\nu(\lambda).
    \end{equation}
\item[\rm(2)] The function $\lambda \mapsto \tilde\mu_\lambda$
  is unique $\tilde\nu$ almost   everywhere. 
\item[\rm(3)] $\tilde\mu_\lambda(C \setminus \cO_\lambda)= 0$ for
  $\tilde\nu$ almost every $\lambda \in S$. 
  \end{itemize}

  \nin {\bf Step 4:} For $g \in G$, the relation $g_*\mu = \mu$ implies that
  \[ g_* \tilde \mu = c_g \tilde \mu \quad \mbox{ for }\quad
    c_g = e^{H_x - H_x \circ \Ad^*(g)^{-1}},\]
  resp., $\tilde \mu = c_g^{-1} g_*\mu$.  Writing
  \eqref{eq:tildemu-dec} as
  \[ \tilde\mu = \int_S \tilde\mu_\lambda\, d\tilde\nu(\lambda),\]
  we thus obtain 
  \[ \int_S c_g\cdot \tilde\mu_\lambda\, d\tilde\nu(\lambda)
=     c_g \tilde\mu
= g_*\tilde\mu = \int_S g_*\tilde\mu_\lambda\, d\tilde\nu(\lambda).\]

Property (3) implies that, for
  almost every $\lambda \in S$,  the measure $\tilde\mu_{\lambda}$ 
  is a Borel measure on the coadjoint orbit $\cO_\lambda$. 
  Let $\Gamma \subeq G$ be a dense countable subgroup.
  Then the uniqueness property (2) implies that, for almost every 
  $\lambda\in S$, we have
  \begin{equation}
    \label{eq:Gamma-quasiinv}
 g_*\tilde\mu_\lambda = c_g \cdot \tilde\mu_\lambda
 \quad \mbox{ for } \quad g \in \Gamma.
  \end{equation}
We may thus assume w.l.o.g.\  that this is the case for every $\lambda \in S$.

In view of  \cite[Thm.~VI.10]{Fa00},  the natural map
  $G/G^\lambda \into C, gG^\lambda \mapsto \Ad^*(g)\lambda$ is a topological
  embedding. The regularity of the measure
  $\tilde\mu_{\lambda}$ on $\cO_\lambda$ thus follows from 
  \cite[Thm.~2.18]{Ru86}, so that it is a Radon measure on $\cO_\lambda$.
  Now \eqref{eq:Gamma-quasiinv} implies that this relation holds for
  every $g \in G$. Therefore the measure
  $e^{H_x} \tilde\mu_\lambda$ on $\cO_\lambda$ is $G$-invariant,
  hence of the form $c_\lambda \mu_\lambda$,
  where $\mu_\lambda$ is the $G$-invariant
  Liouville measure on $\cO_\lambda$.

  \nin {\bf Step 5:} This leads to 
  \[ \mu = e^{H_x} \tilde \mu
    = \int_{S} e^{H_x} \tilde\mu_\lambda\, d\tilde\nu(\lambda)
    = \int_{S}  c_\lambda \mu_\lambda\, d\tilde\nu(\lambda),\]
  which is the desired disintegration for 
  $d\nu(\lambda) = c_\lambda d\tilde\nu(\lambda)$. 
\end{prf} 

At this point one may wonder which measures $\mu$ on $\g^*$
occur naturally for Hamiltonian $G$-actions and
$\mu = \Psi_*\lambda_M$, where $\lambda_M$ is the Liouville measure
on~$M$. A particularly interesting class of examples arises as follows.

\subsection*{Open domains in $T^*(\Gamma\backslash G)$} 

Let $\Omega \subeq \g$ be an open convex set on which
we have a smooth convex function $f \: \Omega \to \R$ that is
strictly convex and has a closed epigraph.
Then $\dd f \:  \Omega \to \g^*$ maps
$\Omega$ diffeomorphically onto an open $\Ad^*(G)$-invariant
subset $\cC \subeq \g^*$. We thus obtain an open subset
\[ \cC_G := G \times \cC \subeq G \times \g^* \cong T^*(G) \]
of the symplectic manifold $T^*(G)$, on which $G$ acts 
by right translations in a Hamiltonian fashion with momentum map
\[  \cC_G \to \g^*, \quad (g,\alpha) \mapsto \alpha\]
(cf.\ \cite[\S III]{Ne00b}). Let $\Gamma \subeq G$ be a lattice,
i.e., a discrete subgroup for which $\Gamma\backslash G$ has finite
volume. Then
\[ M := \Gamma\backslash \cC_G \subeq T^*(\Gamma \backslash G) \]
is an open $G$-right-invariant subset, the $G$-right action is Hamiltonian,
and the momentum set takes the form
\[ \Psi \: M \to \g^*, \quad (\Gamma g,\alpha) \mapsto \alpha.\]
As $\vol(\Gamma\backslash G) < \infty$, the Liouville measure
$\lambda_M$ projects onto a multiple of Lebesgue measure
$\lambda_{\g^*}$, restricted to $\cC$. Since $G$ is unimodular
by \cite[Thm.~VII.1.8]{Ne00}, the coadjoint action preserves
any Lebesgue measure on $\g^*$.

We conclude that measures of the form
$\mu := \lambda_{\g^*}\res_{\cC}$ occur as the image
of the Liouville measure for a Hamiltonian $G$-action. If $\cC$
contains no affine lines, the temperedness of Lebesgue measure
implies that
\[ f(x) := \log \cL(\mu)(x)
  = \log \int_{\cC} e^{-\alpha(x)}\, d\lambda_{\g^*}(\alpha) \]
is finite on the open cone $B(\cC)^\circ$ 
(Proposition~\ref{prop:dom-lapl-temp}). 
If $\cC$ is a cone, this is the logarithm of the
Koecher--Vinberg characteristic function of the
cone~$C$.

If $x \in \partial B(\cC)$, then there exsits
$\alpha \in \lim(\cC)$ with $\alpha(x) = 0$. For any open subset
$O \subeq \cC$ we then have $O + \R_+ \alpha \subeq \cC$
and the Lebesgue measure of this set is infinite. This implies that
$\cL(\mu)(x) = \infty$. So $D_\mu =  B(\cC)^\circ$ and
\[ \dd f \: B(\cC)^\circ \to C_\mu^\circ = \cC\]
is a diffeomorphism by Theorem~\ref{thm:conv-lapl}.

\begin{ex} (a) If $G = (\g,+)$ is abelian, then
  $\g \cong \R^n$ and $\Gamma = \Z^n$ is a lattice in $G$.

\nin (b) The Theorem of Borel--Harish-Chandra \cite[p.2]{Zi84}
  (see also \cite[Thm.~7.8]{BHC62}, \cite[Thm.~14.1]{Ra72}),
  combined with Chevalley's Theorem on the existence of
  $\Z$-basis in simple real Lie algebras, implies that every
  connected semisimple Lie group $G$ contains a lattice~$\Gamma$.
  If $G$ is the identity component of an algebraic group defined over $\Q$,
  then the $\Z$-points of $G$ are such a lattice.
  We thus obtain in particular the lattice $\Gamma = \Sp_{2n}(\Z)
  \subeq G = \Sp_{2n}(\R)$.

  \nin (c) In the Jacobi group
  \[ G = \Heis(\R^{2n},\omega) \rtimes \Sp_{2n}(\R)
    = \R \times \R^{2n} \rtimes \Sp_{2n}(\R) \]
  we have the lattice 
  \[ \Gamma = \Z \times \Z^{2n} \rtimes \Sp_{2n}(\Z) \]
  (cf.\ \cite[Thm.~9.4]{BHC62}).

  \nin (d) If $G$ contains a lattice $\Gamma$, then
  $\Ad(G)$ is closed by \cite[Thm.~2]{GG66}.
  For an admissible Lie algebra $\g = \fu \rtimes \fl$
  and $\ft = \fz(\g) \oplus \ft_\fl$, \cite[Prop.~VII.1.4]{Ne00}
  implies that $\Ad(G)$ is closed if and only if
  $e^{\ad \fz(\fl)}$ is closed. It is easy to construct examples
  where $\fl = \fz(\fl)$ is abelian and this is not the case.
  The simplest ones are of the form
  \[ \g = \Heis(\R^4,\Omega) \rtimes \R D,\]
  where $D \in \sp_4(\R)$ is of the form
  \[ D = \pmat{
      0 & 1 & 0 & 0 \\
-1 & 0 & 0 & 0 \\
0 & 0 & 0 & \sqrt{2}  \\
0 & 0 & -\sqrt{2} & 0}.\]
In this case the closure of $\exp(\R D)$ is a $2$-dimensional torus.
\end{ex}

\section{Non-strongly Hamiltonian actions}
\mlabel{sec:non-strong}

As already noted in the introduction, one may also consider
symplectic actions $\sigma \: G \times M \to M$ of a connected
Lie group $G$ on a connected symplectic manifold that are
Hamiltonian in the  sense that all vector fields
$\dot\sigma(x)$ on $M$ are Hamiltonian, but the homomorphism
$\dot\sigma \: \g \to \Ham(M,\omega)$ may not lift to a homomorphism
to $(C^\infty(M), \{\cdot,\cdot\})$. These actions are not strongly
Hamiltonian. As
\[ \R 1 \into C^\infty(M) \onto \Ham(M,\omega) \]
is a central extension of Lie algebras, this obstruction can always
be overcome by replacing $\g$ by a central extension
\[ \g^\sharp = \R \oplus_\beta \g
  \quad \mbox{ with } \quad
  [(t,x),(t',x')] = (\beta(x,x'), [x,x']).\] 
Then the corresponding simply connected Lie group $G^\sharp$ is a
central extension of $G$ that acts on $M$ with an equivariant momentum
map 
\[ \Psi^\sharp \: M \to \{1\} \times \g^* \subeq  (\g^\sharp)^*
  \cong \R \times \g^*.\]
The coadjoint action of $G^\sharp$ on $\g^\sharp$ factors through an action
of $G$ that leaves the affine hyperplane $\{1\} \times \g^*$ invariant.
So $\Psi^\sharp$ can be considered as a map $M \to \g^*$ that is equivariant
with respect to an action of $G$ on $\g^*$ by affine maps.

Having this in mind, one may always translate between Hamiltonian
actions of $G$ with a momentum map equivariant for an affine action and
strongly Hamiltonian actions of a central extension $G^\sharp$.
As we throughout adopted the latter perspective, we briefly
discuss this translation in the thermodynamic context.

As before, we assume that $M$ is connected and that
$\Psi(M)^\bot = \{0\}$, i.e., that the Lie algebra $\g$ acts effectively
on $M$. Then $\Psi(M)$ spans $\g^*$ and one of the following two cases occurs:
\begin{itemize}
\item[(A)] Affine type: Then $\Psi(M)$ is contained in a proper affine
  hyperplane of $\g^*$. Then $\g$ contains central elements with non-zero
  constant Hamiltonian function, so that $\fz(\g) \not=\{0\}$ is $1$-dimensional. Thus $\g$ is a central extension of $\g^\flat := \g/\fz(\g)$
  and the corresponding quotient group $G^\flat := G/Z(G)_e$ acts
  on $M$ with a momentum map that is equivariant for an affine action. 
\item[(L)] Linear type: Then $\Psi(M)$ is not contained in a proper affine 
  hyperplane of $\g^*$. Since central elements of $\g$ are constant on
  $\Psi(M)$, it follows that $\fz(\g) = \{0\}$.
\end{itemize}

Recall from Theorem~\ref{thm:spind} that admissible Lie algebras can always be written as
\[ \g = (\fz \oplus V)  \rtimes \fl \quad \mbox{ with } \quad
  [V,V] \subeq \fz\]
and $\fl$ reductive.
So $\fz(\g) \not=\{0\}$ always holds if $\g$ is not reductive, and then
\[ \g/\fz(\g) \cong V \rtimes \fl.\]
But $\g$ may also be reductive, i.e., $\g = \fl$, with non-trivial
center. Then $\g$ is a trivial central extension of the semisimple
Lie algebra $[\g,\g]$, so that in this case the affine action
of $G^\flat$ always has a fixed point, hence can be linearized.

This discussion shows that the non-reductive Lie groups
$G^\flat$ that may possess (non-strongly) Hamiltonian actions with
non-trivial geometric temperatures have Lie algebras of the form
\[ \g^\flat \cong V \rtimes_\delta \fl,\]
where $V$ is an abelian ideal carrying a symplectic form $\Omega$ for which
$\delta(\fl) \subeq \sp(V,\Omega)$ contains elements with a positive
definite Hamiltonian function
(cf.~Theorem~\ref{thm:spind}).

\section{Perspectives}
\mlabel{sec:7} 

In this final section, we collect some references
and possibly interesting connections with other areas. 

\subsection{Non-commutative relatives}

In \cite{St96, St99, NS99}, Nencka and Streater consider
non-commutative statistical manifolds obtained from a 
unitary Lie group representation $U \:  G \to \U(\cH)$.
Let $\partial U(x)$ be the skew-adjoint infinitesimal generator of
the unitary one-parameter group $(U(\exp tx))_{t \in \R}$, $x \in \g$.
We call 
\[ \Omega_U := \{ x \in \g \:  \tr(e^{-i\partial U(x)}) < \infty \} \]
the corresponding {\it trace class domain}. 
In \cite[Thm.~2.3.1]{Si23}, T.~Simon proves the following  result,
which is a ``non-commutative'' analog of our Domain Theorem~\ref{thm:dom}: 

\begin{thm} If $(U,\cH)$ is irreducible and $\ker U$ discrete, then 
  $\g$ is admissible and
  there exists an adapted positive system~$\Delta^+$ with $C_{\rm min}
  \subeq C_{\rm max}$ such that $\Omega_U = W_{\rm max}^\circ.$.
\end{thm}

The case of reducible representations is more complicated,
but the requirement $\Omega_U^\circ \not=\eset$ implies that the representation
decomposes as a countable direct sum of irreducible representations
(\cite[Prop.~III.3.18]{Ne00}).

The function
\[ Z \: \Omega_U^\circ \to \R, \quad
  Z(x) := \tr(e^{-i \partial U(x)}) \]
is the non-commutative/quantum analog of the parition function
from thermodynamics. It is also analytic, $G$-invariant and convex,
and even strictly convex if $U$ has discrete kernel
(cf.\ \cite{Ne96a}). 

In this context the natural Riemannian metric on $\Omega_U$, 
specified by  the second derivative of $\log Z$, 
is called the {\it Bogoliubov--Kubo--Mori metric}
(cf.~also \cite{Ta06}). In this context Balian's paper \cite{Ba05}
is particularly interesting, where, for finitely many
selfadjoint operators $H_1, \ldots, H_n$,
Gibb's ensembles are parametrized by
\[ \Omega := \Big\{ x \in \R^n \: \tr\exp\Big(-\sum_{j = 1}^n x_j H_j\Big) < \infty \Big\},\]
the corresponding Gibbs states are of the form
\[ \exp\Big(-z(x)\1 -\sum_{j = 1}^n x_j H_j\Big) \]
and characterized by maximizing a suitable entropy,
so that the situations very much resembles the geometry
of Theorem~\ref{thm:b.3}.

It would be very interesting to understand the precise relation
between the geometric temperature 
$\Omega_\lambda = W_{\rm max}^\circ$ associated to a coadjoint orbit
$\cO_\lambda$, which for unitary
highest weight representations, coincides with the corresponding
trace class domain by Simon's Theorem. But the
``commutative'' and the ``non-commutative'' partition functions
do not coincide in general. We refer to \cite{Ne96a} for
a detailed discussion of examples.
This leaves the question how they are related on a conceptual level.
A natural key could
be the Duistermaat--Heckman formulas for the holomorphic character
in terms of admissible coadjoint orbits, as described in \cite{Ne96a}.

\subsection{Coherent state orbits and trace class operators}

Let $(U,\cH)$ be a unitary lowest weight representation
of an admissible Lie group $G$ and $[v_\lambda] \in \bP(\cH)$
the lowest weight ray, where $v_\lambda$ is a unit vector
of lowest weight~$\lambda$
(\cite{Ne96a}, \cite{Ne00}). Then the momentum map
\[ \Psi \: \bP(\cH^\infty) \to \g^*, \quad
  \Psi([v])(x) := -i \frac{ \la v, \dd U(x) v}{\la v,v\ra} \]
is $G$-equivariant and maps the complex manifold
$M := G.[v_\lambda]$ diffeomorphically onto the admissible coadjoint
orbit $\cO_\lambda$ (see \cite[Ch.~XV]{Ne00} for
coherent state representations).

Since $\Omega_\lambda = W_{\rm max}^\circ$ is non-trivial, we obtain
on $M$ a family of probability measure~$\mu_x$,
parametrized by $x \in W_{\rm max}^\circ$.
Using the $G$-equivariant embedding
\[ \bP(\cH) \into B_1(\cH), \quad  [v] \mapsto P_v, \quad
  P_v(w) := \frac{\la v, w \ra}{\|v\|^2}v \]
we obtain a $G$-equivariant injection
\[ \Psi \:  \cO_\lambda \to B_1(\cH), \quad
  \Psi(\Ad^*(g)\lambda) = U(g) P_{v_\lambda} U(g)^{-1}.\]
Then
\[ A_x := \int_{\cO_\lambda}\, \Psi(\alpha)\,  d\mu_x(\alpha) \]
defines a positive trace class operator with $\tr(A_x) = 1$. 
The map $\Psi$ is continuous because $G$ acts continuously on $B_1(\cH)$.
Therefore the symbol map
\[  \Psi^\vee \: B(\cH) \to C(\cO_\lambda), \quad
  \Psi^\vee(A)(\alpha) = \tr(A \Psi(\alpha)) \]
is a linear $G$-equivariant map with
\begin{align*}
 \Psi^\vee(A)(\Ad^*(g)\lambda) 
&  = \tr(A U(g) P_{v_\lambda} U(g)^{-1})
  = \tr(U(g)^{-1}A U(g) P_{v_\lambda}) \\
&  = \la v_\lambda, U(g)^{-1}A U(g) v_\lambda\ra 
  = \la U(g) v_\lambda, A U(g) v_\lambda\ra.
\end{align*}
Therefore the map $\Psi^\vee$ may be viewed as a
dequantization or a symbol map, turning operators
into functions. This correspondence is of particular interest
for representations which are square integrable modulo the center,
resp., which can be realized in holomorphic $L^2$-sections
of line bundles; see in particular~\cite{Ne96c, Ne97, Ne00}.

\subsection{Infinite dimensions}

Symplectic manifolds also make sense in infinite dimensions,
but not the Liouville measure. However, measures on infinite-dimensional
spaces make good sense. If, for instance, $\mu$ is a Borel measure
on the dual $V^*$ of the real vector space $V$,
endowed with the smallest $\sigma$-algebra making all evaluations
measurable, then
\[ \cL(\mu) \:  V \to \R \cup \{\infty\}, \quad
  \cL(\mu)(v) = \int_{V^*} e^{-\alpha(v)}\, d\mu(\alpha) \]
is finite, and one can study measures for which
it is finite on a non-empty open subset. Interesting examples
appear in \cite{NO02} on domains in the space of
Hilbert--Schmidt operators. Here the major sources
are Gaussian measure and their images under non-linear maps. 

\begin{ex} To see infinite dimensional examples that are closer to the applications 
in physics, one may also consider Lie algebras 
of the form $\g = \su_2(\C)^{(\N)}$ (countable direct sum), whose dual 
space is the full sequence space $\g^* \cong \su_2(\C)^\N$. 
This space carries many invariant probability measures.
We refer to \cite{NR24} for a discussion of possibly related
unitary representations of infinite-dimensional Hilbert--Lie groups. 
\end{ex}

For information geometry in the infinite-dimensional context
of diffeomorphism groups, we refer to the recent survey \cite{KMM24}.
Results concerning infinite-dimensional convex functions
can be found in \cite{Mi08}, \cite{Bou07} and \cite[\S 3]{Ro74}. \\

In \cite{Fr91} Friedrich's discussion of the Fisher--Rao metric on the
space of probability measures is infinite-dimensional in spirit.
For a probability space $(X,\fS,\mu)$, he considers the set $\cA$ of
all probability measures of the form $f \mu$, with the
tangent space in $\mu$ given by 
\[ T_\mu(\cA) = \Big\{ f \in L^2(X,\mu) \: \int_X f \, d\mu = 0\Big\}, \] 
endowed with the Riemannian metric inherited from $L^2(X,\mu)$.
For the case where $\mu$ comes from an $n$-form $\lambda$ on a manifold $M$,
Friedrich even associates to each vector field preserving $\lambda$ a
Poisson structure on the corresponding manifold $\cA$ of probability
measures with smooth densities. For $M = \bS^1$, this leads to the symplectic
structure corresponding to identifying $\cA$ with a coadjoint orbit
of the infinite-dimensional group $\Diff(\bS^1)_+$
(\cite[Bem.~2]{Fr91}).

\subsection{Weinstein's modular automorphisms} 

Let $(M,\omega)$ be a symplectic manifold, $\mu$ its Liouville measure and 
$H \: M \to \R$ a smooth function for which 
$e^{-H}\mu$ is a finite measure. For a Hamiltonian vector field 
$X_F$ with $X_F G = \{G,F\}$ for $G \in C^\infty(M)$, we then have 
\[ \cL_{X_F} (e^{-H}\mu) 
= -X_F(H)(e^{-H}\mu) = \{F,H\}(e^{-H}\mu).\]
Therefore 
\[ \div_{e^{-H}\mu}(X_F) = \{F,H\} = X_H(F).\] 
This shows that the modular flow corresponding to the 
``KMS state'' $e^{-H}\mu$  in the sense of \cite{We97}
coincides with the flow of the Hamiltonian vector field 
$X_H$ on $M$. We refer to \cite{We97} for a discussion of
KMS states in the context of Poisson- and symplectic manifolds.
More recent results in this context can be found in \cite{DW23}.
This paper also contains for a connected symplectic manifold $M$
a characterization of the measures of the form
$e^{-H}\lambda_M$ as the KMS functionals corresponding to the flow
generated by the Hamiltonian function~$H$. Finiteness of
these measures is only discussed in \cite{DW23} for the trivial
case where $M$ is compact. A corresponding result in
the context of deformation quantization
is stated in \cite[Thm.~4.1]{BRW98}, characterizing KMS states
as Gibbs states.

\begin{ex}
As the context of Weinstein's paper is Poisson manifolds, 
one may also consider open domains $M \subeq \g^*$, where  
$\g$ is a finite dimensional Lie algebra. Here the case 
where $M$ is the interior of a cone $W_{\rm min}^\star$,
or the interior of the convex hull of an orbit
$\cO_\lambda$ with $\Omega_\lambda\not=\eset$ provide
interesting examples, connecting with information geometry.
\end{ex}

\subsection{Locally symmetric spaces}

Let $G$ be a linear semisimple Lie group,
$K \subeq G$ be a maximal compact subgroup and $G/K$
the corresponding non-compact Riemannian symmetric space.
If $\Gamma \subeq G$ is a torsionfree lattice,
$X := \Gamma\backslash G/K$ is called a {\it locally symmetric space}.
Then $\vol(X) < \infty$ and the Liouville measure
$\lambda_M$ on the symplectic manifold $M := T^*(X)$ has
strong finiteness properties. For example the energy function
\[ H \: T^*(X) \to \R, \quad H(\alpha) = \frac{1}{2}\|\alpha\|^2 \] 
is the Hamiltonian function of the geodesic flow on $T^*(X) \cong T(X)$.
Since $X$ has finite volume, it follows that 
\[ Z(\beta) := \int_M e^{-\beta H}\, d\lambda_M
  = \int_X \Big(\int_{T^*_p(X)} e^{-\frac{\beta}{2}\|\alpha\|^2}\, d\alpha\Big)
  \, d\mu_X(p) < \infty\]
for every $\beta > 0$. 

\begin{ex} For $G = \PSL_2(\R)$ and $K = \PSO_2(\R)$,
  $G/K$ is the hyperbolic plane, resp., the open unit disc in $\C$,
  and the group $G$ acts transitively on the
  level sets of the energy function in $T^*(G/K)$. Factorization 
  of a lattice $\Gamma$, leads to submanifolds of finite volume.

  In this case the geodesic flow on $X$ can be implemented by
  the subgroup $A \cong \PSO_{1,1}(\R) \cong \R$. We thus obtain
  a Gibbs measure on $T^*(X)$ for the action of a hyperbolic
  one-parameter group which acts ergodically on the level sets of~$H$
  (cf.\ also \cite[p.~386]{We97}).
\end{ex}

For groups of rank $r> 1$, one has
Hilgert's Ergodic Arnold--Liouville Theorem 
(\cite[Thm.~8.3(v)]{Hi05}) which specifies a
Poisson commuting set $C_1, \ldots, C_r$ of smooth functions on $T^*(X)$ that
are obtained from $G$-invariant functions on $T^*(G/K)$ by factorization.
Any finite-dimensional linear subspace $\fh$ of the algebra $\cA$ generated 
by  these functions that leads to complete Hamiltonian vector fields
defines a Hamiltonian action of $H = \R^r$ on $T^*(X)$ and one may expect
that suitable choices even lead to a non-trivial geometric
temperature, as for the geodesic flow and $r= 1$.


\end{document}